\documentclass[11pt,reqno]{amsart}

\usepackage[text={160mm,240mm},centering]{geometry}            
\usepackage{amssymb,amsmath,setspace,geometry,indentfirst,changepage,inputenc,amsthm}
\usepackage{mathrsfs,amsfonts,savesym,graphicx,bm,color}

\usepackage{graphicx}
\usepackage{amssymb}
\usepackage{epstopdf}
\usepackage{dsfont}
\usepackage{placeins}
\usepackage{pifont}
\usepackage{romannum}

\usepackage{lscape}

\usepackage{tikz}
\usetikzlibrary{positioning}
\usepackage{xcolor}

\usepackage{float}

\usepackage[all,2cell]{xy}    
\UseAllTwocells               

\usepackage{amsmath,amsthm}

\usepackage[mathscr]{eucal}
\usepackage[all]{xy}
\usepackage{mathrsfs}
\usepackage{hyperref}
\usepackage{color}

\newtheorem{definition-lemma}[theorem]{Definition-Lemma}
\newtheorem{definition-theorem}[theorem]{Definition-Theorem}

\theoremstyle{definition}

\allowdisplaybreaks

\title[]{Hilbert Series of Second Order Jets of  Determinantal Varieties}

\author{Yifan Chen}
\address{
	Department of Mathematical Sciences,
	Tsinghua University,
	Beijing, 100084, P. R. China.}
\email{c-yf20@tsinghua.org.cn}

\author{Yongxin Xu}
\address{
	Zhili College,
	Tsinghua University,
	Beijing, 100084, P. R. China.}
\email{xuyongxi22@mails.tsinghua.edu.cn}

\author{Huaiqing Zuo}
\address{Department of Mathematical Sciences,
	Tsinghua University,
	Beijing, 100084, P. R. China.}
\email{hqzuo@mail.tsinghua.edu.cn}

\begin{document}
	
	\setcounter{page}{1}
	\pagenumbering{arabic}
	
	\maketitle
	
	\parskip=.3em
	
	\setcounter{tocdepth}{1}

	\begin{abstract}
		In this paper, we will investigate the jet schemes of determinantal varieties.  It is quite often the case that the geometric information concerning the jet schemes of an algebraic variety can be described, but the more refined algebraic information is quite mysterious. For example, it is known that computing the Hilbert function associated to a natural grading on these jet schemes is a very hard problem. The present paper handles a few such computations. It succeeds in computing the Hilbert functions of the second order jet schemes in the case of maximal minors of a $2\times n$ matrix.
		
		
	\end{abstract}
	
	\section{Introduction}
	\par In recent years, the theory of jet schemes has been widely studied and applied since it was introduced by Nash in his paper \cite{c5}. Ein and Mustaţă studied some invariants of singularities using jet schemes, with results collected in \cite{c6}. We will discuss jet schemes of determinantal varieties in this paper.
	\par Let $K$ be an algebraically closed field of characteristic zero and $\mathbf{A}^{mn}_{K}$ be the affine space of dimension $mn$ with $m\le n$. For a fixed $m \times n$ matrix $(x_{i,j})_{1\le i\le m, 1\le j\le n}$, all $r \times r$ minors of the matrix generate an ideal $I^{m,n}_{r}$ in $K[x_{i,j}|1\le i\le m, 1\le j\le n]$, which determines an algebraic variety $\mathcal{L}^{m,n}_{r}$. We say $\mathcal{L}^{m,n}_{r}$ is a determinantal variety of $m \times n$ matrix.
	\par Let $x_{i,j}(t)=\underset{s=0}{\overset{\infty}{\sum}} x^{(s)}_{i,j}t^{s}$ and $X(t)=(x_{i,j}(t))_{1 \le i \le m,1 \le j \le n}$ an $m \times n$ matrix. Every $r \times r$ minor of $X(t)$ is of the form $\underset{l=0}{\overset{\infty}{\sum}}f_{l}(x^{(s)}_{i,j})t^{l}$, where $f_{l}$ is a polynomial in $K[x_{i,j}^{(s)}|s\in \mathbb{Z}_{\ge 0},1\le i\le m,1\le j\le n]$. We denote $I^{m,n}_{r,k}$ to be the ideal generated by $f_{0},f_{1},\cdots,f_{k}$ from all $r \times r$ minors. The variety determined by $I^{m,n}_{r,k}$ is the $k$th jet scheme of $\mathcal{L}^{m,n}_{r}$, denoted by $\mathcal{L}^{m,n}_{r,k}$. The $0$th jet scheme of $\mathcal{L}^{m,n}_{r}$ is just itself.
	\par The Hilbert series of the classical determinantal varieties was obtained by Alde Conca and Jürgen Herzog in their paper \cite{c9}.  \\\\
	\textbf{Theorem 1.1. }(\cite{c9}, Corollary 1) \textit{The Hilbert series of $\mathcal{L}^{m,n}_{r+1}$ is}
	\begin{equation}
		\centering
		H_{\mathcal{L}^{m,n}_{r+1}}(z)= \frac{det\big( \underset{k\ge 0}{\sum}\binom{m-i}{k} \binom{n-j}{k} z^{k}\big)_{1\le i\le r,1\le j\le r}}{z^{\binom{r}{2} }(1-z)^{r(m+n-r)}}.
	\end{equation}\\
	\par It is complicated to compute the Hilbert series for jets of arbitrary order of determinantal varieties. The simplest nontrivial case is when $(r,k)=(2,1)$, and this case has been addressed in  three papers. In \cite{c3}, the authors gave a Gröbner basis for $I^{m,n}_{2,1}$. Then Boyan Jonov characterized the Stanley-Reisner simplicial complex of the leading ideal of $I^{m,n}_{2,1}$ and showed its shellability (see Definition 2.11) in \cite{c4}. Finally in \cite{c8}, the authors computed the Hilbert series of $\mathcal{L}^{m,n}_{2,1}$, obtaining the following result:
	\begin{equation}
		\centering
		H_{\mathcal{L}^{m,n}_{2,1}}(z)=\Bigg( 
		\frac{\underset{k=0}{\overset{m-1}{\sum}}\binom{m-1}{k} \binom{n-1}{k} z^{k}}{(1-z)^{m+n-1}}
		\Bigg)^{2}.
	\end{equation}
	\par This study extends the work of these three papers. The computation of Gröbner bases of $\mathcal{L}^{m,n}_{r,k}$ is quite difficult when $r$ and $k$ are larger. Hence we deal with cases for small $r$ and $k$. The computation is involved and highly non-trivial even in these cases. In \cite{c10}, Chen and Zuo have obtained some results on the Hilbert series of $\mathcal{L}^{3,n}_{3,1}$, with result as follows:
	\begin{equation}
		H_{\mathcal{L}^{3,n}_{3,1}}(z)=\big(  \frac{1+(n-2)z+\frac{(n-1)(n-2)}{2}z^{2}}{(1-z)^{2n+2}}\big)^{2}.
	\end{equation}
	\par In this paper we will compute the Hilbert series of $\mathcal{L}^{2,n}_{2,2}$. For convenience, we will denote $x^{(0)}_{i,j}+x^{(1)}_{i,j}t+x^{(2)}_{i,j}t^{2}$ by $x_{i,j}+y_{i,j}t+z_{i,j}t^{2}$. The form of $X(t)$ in this case is as follows:
	\begin{equation}
		X(t)=\left(
		\begin{array}{cccc}
			x_{1,1}+y_{1,1}t+z_{1,1}t^{2} & x_{1,2}+y_{1,2}t+z_{1,2}t^{2} & \cdots & x_{1,n}+y_{1,n}t+z_{1,n}t^{2} \\ x_{2,1}+y_{2,1}t+z_{2,1}t^{2} & x_{2,2}+y_{2,2}t+z_{2,2}t^{2} & \cdots & x_{2,n}+y_{2,n}t+z_{2,n}t^{2}
		\end{array}
		\right).
	\end{equation}
	\par We state our main result here: \\\\
	\textbf{Main Theorem. }\textit{Let $L(I_{2,2}^{2,n})$ be the leading ideal of $I^{2,n}_{2,2}$. Suppose  $\Delta_{L(I^{2,n}_{2,2})}$ is the Stanley-Reisner simplicial complex corresponding to $K[x_{i,j},y_{i,j},z_{i,j}|1\le i\le 2,1\le j\le n]/L(I_{2,2}^{2,n})$. Then $\Delta_{L(I^{2,n}_{2,2})}$ is shellable and the Hilbert series of $\mathcal{L}^{2,n}_{2,2}$ is  }
	\begin{center}
		$H_{\mathcal{L}^{2,n}_{2,2}}(z)=\big( \frac{1+(n-1)z}{(1-z)^{n+1}}\big) ^{3}$.
	\end{center}
	\par We will give a short introduction to Gröbner bases, Stanley-Reisner complexes and shellable complexes in the second section. Then we will find a Gröbner basis for $I_{2,2}^{2,n}$ in the third section and find all facets of $\Delta_{L(I^{2,n}_{2,2})}$ in the fourth section. Finally we will show $\Delta_{L(I^{2,n}_{2,2})}$ is shellable in the fifth section and compute the Hilbert series of $\mathcal{L}^{2,n}_{2,2}$ in the sixth section. Indeed, we also propose a conjecture about the relation between the Hilbert series of determinantal varieties and their $k$th order jets, which is included in the sixth section. We introduce it here first.\\\\
	\textbf{Conjecture. }The Hilbert series of $\mathcal{L}^{m,n}_{r,k}$ is the $(k+1)$-th power of the Hilbert series of $\mathcal{L}^{m,n}_{r}$.  \\
	
	\section{Preliminaries}
	\subsection{Gröbner Bases}
	\par $\quad$
	\par In this section, $K$ is always a field and $K[x_{1},\cdots,x_{n}]$ is the polynomial ring in $n$ indeterminates. A polynomial of the form
	\begin{center}
		$t=x_{1}^{e_{1}}x_{2}^{e_{2}} \cdots x_{n}^{e_{n}}$($e_{i}\in \mathbb{N}$)
	\end{center}
	\par \noindent is called a monomial.\\
	\textbf{Definition 2.1. }Let $M$ be the set of all monomials in $K[x_{1},\cdots,x_{n}]$. A monomial ordering on $K[x_{1},\cdots,x_{n}]$ is a total ordering "$\le$"on $M$ satisfying the following properties:
	\par (1)$1\le t$ for all $t\in M$.
	\par (2)if $s$, $t_{1}$, $t_{2}\in M$ and $t_{1}\le t_{2}$, then $st_{1}\le st_{2}$.\\
	\par If $f\in K[x_{1},\cdots,x_{n}]$ is a nonzero polynomial written as $f=c_{1}t_{1}+c_{2}t_{2}+\cdots+c_{m}t_{m}$ with $t_{1}> \cdots > t_{m}$ under a fixed monomial ordering, then we define $lm(f)=t_{1}$ (leading monomial), $lc(f)=c_{1}$ (leading coefficient), and $lt(f)=c_{1}t_{1}$ (leading term). Let $S\subset K[x_{1},\cdots ,x_{n}]$ be a subset of polynomials. The ideal generated by all leading monomials in $S$ is called the leading ideal of $S$, denoted by $L(S)$ .
	\par For convenience, we fix a monomial ordering $\le$ on $K[x_{1},...,x_{n}]$ for following discussion.\\\\
	\textbf{Definition 2.2. }Let $I\subset K[x_{1},\cdots ,x_{n}]$ be an ideal. A finite subset $G\subset I$ is called a Gröbner basis (with respect to the chosen monomial ordering "$\le$") of $I$ if $L(I)=L(G)$.\\\\
	\textbf{Proposition 2.3. }(\cite{c1}, Corollary 9.10) \textit{Let $G$ be a Gröbner basis of an ideal $I\subset K[x_{1},\cdots ,x_{n}]$. Then $I$ is equal to the ideal generated by $G$.}\\
	\par It is not easy to compute Gröbner bases of an ideal. Here we introduce Buchberger's criterion to identify whether a finite set is a Gröbner basis of the chosen ideal or not. \\\\
	\textbf{Definition 2.4. }Let $S=\{g_{1},\cdots ,g_{r}\}$ be a finite subset of $K[x_{1},\cdots ,x_{n}]$ and $f\in K_[x_{1},\cdots ,x_{n}]$. We say $h\in K[x_{1},\cdots ,x_{n}]$ is a normal form of $f$ with respect to $S$ if 
	\par (1)no monomial of $h$ is divisible by $lm(g_{i})$ for any $i=1,\cdots ,r$.
	\par (2)there exist polynomials $h_{1},\cdots ,h_{r}\in K[x_{1},\cdots ,x_{n}]$ such that $f-h=\underset{i=1}{\overset{r}{\sum}}g_{i}h_{i}$ and $lm(h_{i}g_{i})\le lm(f)$ for all $i$.
	\par \noindent We write $f\to_{S}0$ if $f\in K[x_{1},\cdots ,x_{n}]$ satisfies (2) in Definition 2.4.\\\\
	\textbf{Theorem 2.5. }(\cite{c1}, Theorem 9.9) \textit{Let $G$ be a Gröbner basis of an ideal $I\subset K[x_{1},\cdots ,x_{n}]$. Every $f\in K[x_{1},\cdots ,x_{n}]$ has precisely one normal form with respect to $G$.}\\
	\par For $f,g\in K[x_{1},\cdots ,x_{n}]$ nonzero, let $t$ be the gcd of $lm(f)$ and $lm(g)$. Then 
	\begin{center}
		$S(f,g)=\frac{lt(g)}{t}f-\frac{lt(f)}{t}g$
	\end{center} 
	\par \noindent is called the s-polynomial of $f$ and $g$. \\\\
	\textbf{Theorem 2.6. }(\cite{c1}, Theorem 9.12) \textit{Let $I$ be an ideal in $K[x_{1},\cdots ,x_{n}]$ and $G=\{g_{1},\cdots ,g_{r}\}$ be a finite set generating $I$. Then $G$ is a Gröbner basis of $I$ if and only if $S(g_{i},g_{j})\to_{G}0$ for all $ 1\le i<j\le r$.}\\
	\par For the condition $S(g_{i},g_{j})\to_{G}0$, we have convenient methods to check in some special cases. One common proposition is as follows:\\\\
	\textbf{Proposition 2.7. }\textit{Let $G$ be a finite subset of $K[x_{1},\cdots ,x_{n}]$. We have $S(f,g)\to_{G}0$ for polynomials $f,g\in K[x_{1},\cdots ,x_{n}]$ if $lm(f)$ and $lm(g)$ are coprime.}\\
	\par We will use the graded reverse lexicographic ordering, as defined below: \\\\
	\textbf{Definition 2.8. }Let $t=x_{1}^{e_{1}}\cdots x^{e_{n}}_{n}$ and $t^{'}=x_{1}^{e_{1}^{'}}\cdots x^{e_{n}^{'}}_{n}$ be monomials. The graded reverse lexicographic ordering on $K[x_{1},x_{2},\cdots ,x_{n}]$ with $x_{1}>x_{2}>\cdots >x_{n}$ is defined by saying that $t\le t^{'}$ if $t=t^{'}$ or deg($t$) $<$ deg($t^{'}$), or deg($t$) $=$ deg($t^{'}$) and $e_{i}>e^{'}_{i}$ for the largest index $i$ with $e_{i} \neq e_{i}^{'}$.\\
	\par For example, we have $x_{1}x_{3}<x_{2}^{2}$ if $x_{1}>x_{2}>x_{3}$.
	\par For $\mathcal{L}^{2,n}_{2,2}$, we use the graded reverse lexicographic ordering on $K[x_{i,j},y_{i,j},z_{i,j}]$ induced by the ordering of variables:
	\begin{center}
		$z_{1,1}>z_{1,2}>\cdots >z_{2,1}>\cdots >z_{2,n}>y_{1,1}>\cdots >y_{2,n}>x_{1,1}>\cdots >x_{2,n}$
	\end{center}
	\subsection{Stanley-Reisner Complexes and Shellable Complexes}
	\par $\quad$\\\\
	\textbf{Definition 2.9. }Let $S=\{1,2,\cdots ,n\}$. An abstract simplicial complex $\Delta$ on $S$ is a collection of subsets of $S$ closed under taking subsets, that is, if $X\in \Delta$ and $Y\subset X$, then $Y\in \Delta$. $\Delta$ is a partially ordered set under the inclusion relation. The maximal elements are called facets of $\Delta$. The minimal elements are called vertices of $\Delta$. The dimension of a facet $X$ is $|X|-1$. The dimension of $\Delta$ is the maximal dimension of its facets, denoted by $dim\Delta$.\\
	\par An abstract simplicial complex corresponds to a ring, which is called its Stanley-Reisner ring, as defined below:\\\\
	\textbf{Definition 2.10. }Let $\Delta$ be an abstract simplicial complex with $n$ vertices. The Stanley-Reisner ring of $\Delta$ is $R_{\Delta}=K[x_{1},x_{2},\cdots ,x_{n}]/I_{\Delta}$, where the Stanley-Reisner ideal $I_{\Delta}$ is generated by $\{x_{i_{1}}x_{i_{2}}\cdots x_{i_{s}}|\{ i_{1},i_{2},\cdots ,i_{s}\}\notin \Delta \}$.\\
	\par Let $I$ be an ideal of $K[x_{1},x_{2},\cdots ,x_{n}]$. If $I$ is radical, then there exists a unique abstract simplicial complex $\Delta_{I}$, such that its Stanley-Reisner ideal is $I$.\\\\
	\textbf{Definition 2.11. }Let $\Delta$ be an abstract simplicial complex. We say $\Delta$ is pure if all facets of $\Delta$ have the same dimension. We say $\Delta$ is shellable, if it is pure and its facets can be ordered as $F_{1}<F_{2}<\cdots <F_{e}$ such that there exist $v\in F_{j}-F_{i}$ and $k<j$ satisfying $F_{j}-F_{k}=\{v\}$ for all $1\le i\le j\le e$. This ordering on facets is called a shelling of $\Delta$.\\\\
	\textbf{Definition 2.12. }Let $\Delta$ be a shellable simplicial complex with a shelling $F_{1}<F_{2}<\cdots <F_{e}$. We define $c(F_{t})=\{v\in F_{t}|$ $F_{t}-F_{s}=\{v\}$ for some $s<t\}$. If $T$ is a family of facets of $\Delta$, we define $h_{j}(T)=|\{F_{t}\in T||c(F_{t})|=j\}|$ and $h_{j}=h_{j}(S)$ where $S$ is the family of all facets of $\Delta$.\\
	\par The following proposition about shelling is important for our computation. \\\\
	\textbf{Proposition 2.13. }(\cite{c2}, Theorem 6.3) \textit{Let $\Delta$ be a shellable simplicial complex with a shelling $F_{1}<F_{2}<\cdots <F_{e}$. Then:}
	\par (1)\textit{the Stanley-Reisner ring $R_{\Delta}$ is a Cohen-Macaulay ring of dimension $d=1+dim\Delta$.}
	\par (2)\textit{the Hilbert series of $R_{\Delta}$ is}\\
	\begin{equation}
		H_{R_{\Delta}}(z)=\frac{\underset{j\ge 0}{\sum}h_{j}z^{j}}{(1-z)^{d}}.
	\end{equation}

	\section{Gröbner bases for $I^{2,n}_{2,2}$}
	\par In this section, we present a Gröbner basis $\Gamma$ for the ideal $I^{2,n}_{2,2}$ and identify leading monomials of polynomials in $\Gamma$.
	\par By computing all second order minors of (4), one gets generators of $I^{2,n}_{2,2}$ as follows:
	\begin{center}
		$a_{p,q}=det\left(
		\begin{array}{cc}
			x_{1,p} & x_{1,q}\\x_{2,p} & x_{2,q}
		\end{array}
		\right),1 \le p<q \le n$,
	\end{center}
	\begin{center}
		$b_{p,q}=det\left(
		\begin{array}{cc}
			x_{1,p} & x_{1,q}\\y_{2,p} & y_{2,q}
		\end{array}
		\right)+det\left(
		\begin{array}{cc}
			y_{1,p} & y_{1,q}\\x_{2,p} & x_{2,q}
		\end{array}
		\right),1 \le p<q \le n$,
	\end{center}	
	\begin{center}
		$c_{p,q}=det\left(
		\begin{array}{cc}
			x_{1,p} & x_{1,q}\\z_{2,p} & z_{2,q}
		\end{array}
		\right)+det\left(
		\begin{array}{cc}
			z_{1,p} & z_{1,q}\\x_{2,p} & x_{2,q}
		\end{array}
		\right)+det\left(
		\begin{array}{cc}
			y_{1,p} & y_{1,q}\\y_{2,p} & y_{2,q}
		\end{array}
		\right),1 \le p<q \le n$.
	\end{center}
	\par We introduce the following polynomials:
	\begin{center}
		$d_{p,q,r}=det\left(
		\begin{array}{ccc}
			x_{2,p} & x_{2,q} & x_{2,r}\\y_{1,p} & y_{1,q} & y_{1,r}\\z_{2,p} & z_{2,q} & z_{2,r}
		\end{array}
		\right)+det\left(
		\begin{array}{ccc}
			x_{2,p} & x_{2,q} & x_{2,r}\\z_{1,p} & z_{1,q} & z_{1,r}\\y_{2,p} & y_{2,q} & y_{2,r}
		\end{array}
		\right),p<q<r$,
	\end{center} 
	\begin{center}
		$e_{p,q,r}=det\left(
		\begin{array}{ccc}
			x_{1,p} & x_{1,q} & x_{1,r}\\y_{1,p} & y_{1,q} & y_{1,r}\\z_{2,p} & z_{2,q} & z_{2,r}
		\end{array}
		\right)+det\left(
		\begin{array}{ccc}
			z_{1,p} & z_{1,q} & z_{1,r}\\y_{1,p} & y_{1,q} & y_{1,r}\\x_{2,p} & x_{2,q} & x_{2,r}
		\end{array}
		\right),p<q<r$,
	\end{center} 
	\begin{center}
		$f_{l,p,q,r}=det\left(
		\begin{array}{cccc}
			0 & x_{2,p} & x_{2,q} & x_{2,r}\\x_{2,l} & y_{2,p} & y_{2,q} & y_{2,r}\\y_{1,l} & z_{1,p} & z_{1,q} & z_{1,r}\\y_{2,l} & z_{2,p} & z_{2,q} & z_{2,r}
		\end{array}
		\right),1 \le l \le p<q<r \le n$,
	\end{center} 
	\begin{center}
		$g_{l,p,q,r}=det\left(
		\begin{array}{cccc}
			x_{2,l} & x_{2,p} & x_{2,q} & x_{2,r}\\y_{2,l} & y_{2,p} & y_{2,q} & y_{2,r}\\z_{2,l} & z_{2,p} & z_{2,q} & z_{2,r}\\z_{1,l} & z_{1,p} & z_{1,q} & z_{1,r}
		\end{array}
		\right),1 \le l<p<q<r \le n$.
	\end{center}
	\par The leading monomials of these polynomials are
	\begin{center}
		$lm(a_{p,q})=x_{1,q}x_{2,p}$, \\ $lm(b_{p,q})=x_{1,p}y_{2,q}$, \\
		$lm(c_{p,q})=y_{1,q}y_{2,p}$,\\
		$lm(d_{p,q,r})=x_{2,p}y_{1,q}z_{2,r}$, \\ $lm(e_{p,q,r})=x_{1,p}y_{1,q}z_{2,r}$, \\
		$lm(f_{l,p,q,r})=x_{2,p}y_{1,l}y_{2,q}z_{2,r}$, \\
		$lm(g_{l,p,q,r})=x_{2,l}y_{2,p}z_{1,r}z_{2,q}$.\\
	\end{center}
	\textbf{Proposition 3.1. }\textit{The polynomials $d_{p,q,r}$, $e_{p,q,r}$, $f_{l,p,q,r}$, $g_{l,p,q,r}$ are belong to $I^{2,n}_{2,2}$.}
	\begin{proof}[Proof.]
		We have following equations:
		\begin{center}
			$d_{p,q,r}=-z_{2,p}b_{q,r}+z_{2,q}b_{p,r}-z_{2,r}b_{p,q}-y_{2,p}c_{q,r}+y_{2,q}c_{p,r}-y_{2,r}c_{p,q}$,\\ $e_{p,q,r}=-y_{1,p}c_{q,r}+y_{1,q}c_{p,r}-y_{1,r}c_{p,q}$,\\
			$f_{l,p,q,r}=z_{2,l}z_{2,p}a_{q,r}-z_{2,l}z_{2,q}a_{p,r}+z_{2,l}z_{2,r}a_{p,q}-z_{2,q}x_{2,p}c_{l,r}+z_{2,r}x_{2,p}c_{l,q}-z_{2,r}x_{2,q}c_{l,p}+y_{2,l}d_{p,q,r}$,\\
			$g_{l,p,q,r}=(z_{2,q}y_{2,r}-y_{2,q}z_{2,r})c_{l,p}-(z_{2,p}y_{2,r}-y_{2,p}z_{2,r})c_{l,q}+(z_{2,l}y_{2,r}-y_{2,l}z_{2,r})c_{p,q}+(z_{2,p}y_{2,q}-y_{2,p}z_{2,q})c_{l,r}-(z_{2,l}y_{2,q}-y_{2,l}z_{2,q})c_{p,r}+(z_{2,l}y_{2,p}-y_{2,l}z_{2,p})c_{q,r}$.\\
		\end{center}
		Since $a_{p,q}$, $b_{p,q}$, $c_{p,q}\in I_{2,2}^{2,n}$, it follows that $d_{p,q,r}$, $e_{p,q,r}$, $f_{l,p,q,r}$, $g_{l,p,q,r}\in I^{2,n}_{2,2}$. 
	\end{proof}	
	\par Let $\Gamma$ be the set of all $a_{p,q}$, $b_{p,q}$, $c_{p,q}$, $d_{p,q.r}$, $e_{p,q,r}$, $f_{l,p,q,r}$, $g_{l,p,q,r}$. We aim to prove that $\Gamma$ forms a Gröbner basis for $I^{2,n}_{2,2}$. By Buchberger's criterion (Theorem 2.6), it suffices to show that for every pair of polynomials $\alpha$, $\beta \in \Gamma$, the s-polynomial $S(\alpha,\beta)$ can be expressed as $\underset{\theta \in \Gamma}{\sum} f_{\theta} \theta$ with $lm(f_{\theta} \theta) \le lm(S(\alpha,\beta))$ for all $\theta\in \Gamma$. By the following lemma from \cite{c3}, we can simplify the checking work.\\ \\
	\textbf{Lemma 3.2. }(\cite{c3}, Lemma 2.5) \textit{Let $S(\alpha,\beta)=\underset{\theta \in \Gamma}{\sum}f_{\theta}\theta$ satisfying $lm(f_{\theta} \theta) \le lm(S(\alpha,\beta))$, and all the variables of $\alpha$, $\beta$, $\theta$, $f_{\theta}$ belong to some fixed submatrix $T= \{ (i_{s},j_{s})|i_{1}<\cdots <i_{p}$, $j_{1}<\cdots <j_{q} \}$, $p$, $q \in \mathbb{N}$. If we replace $T$ with another submatrix $T^{'}= \{ (i^{'}_{s},j^{'}_{s})|i^{'}_{1}<\cdots <i^{'}_{p}$, $j^{'}_{1}<\cdots <j^{'}_{q} \} $ of the same size,then we will get new polynomials $\alpha^{'}$, $\beta^{'}$, $\theta^{'}$, $g_{\theta^{'}}$ and $S(\alpha^{'},\beta^{'})=\underset{\theta^{'} \in \Gamma}{\sum}$, $lm(g_{\theta^{'}} \theta^{'}) \le lm(S(\alpha^{'},\beta^{'}))$.}\\
	\par For each $\alpha \in \Gamma$, variables of $\alpha$ comes from at most a $2 \times 4$ matrix. So a $2 \times 8$ matrix can contain all variables of $\alpha$, $\beta \in \Gamma$. If $lm(\alpha)$, $lm(\beta)$ are coprime, then $S(\alpha,\beta)\to_{\Gamma}0$ must hold. So it suffices to consider cases that submatrices of $\alpha$ and $\beta$ have coincident columns. In conclusion, we only have to prove that $\Gamma$ is a Gröbner basis for the case $4 \le n \le 7$. We computed $L(I^{2,n}_{2,2})$ for $n=4,5,6,7$ by the computer algebra system $\mathtt{SINGULAR}$. The results confirmed that $L(I^{2,n}_{2,2})$ is the ideal generated by all leading monomials of polynomials in $\Gamma$ in each case.
	
	\section{Leading Ideal of $I^{2,n}_{2,2}$ and Facets of $\Delta_{L(I^{2,n}_{2,2})}$}
	\par In this section, we describe the facets of the Stanley-Reisner complex associated to the leading ideal $L(I^{2,n}_{2,2})$.
	\par The leading ideal $L(I^{2,n}_{2,2})$ is generated by leading monomials of polynomials in $\Gamma$. We respectively denote the families of sets corresponding to leading monomials of $a_{p,q}, b_{p,q}, c_{p,q}, d_{p,q.r}, e_{p,q,r}, f_{l,p,q,r}, g_{l,p,q,r}$ as $A, B, C, D, E, F^{*}, G$, which can be respectively depicted in the following diagrams:\\
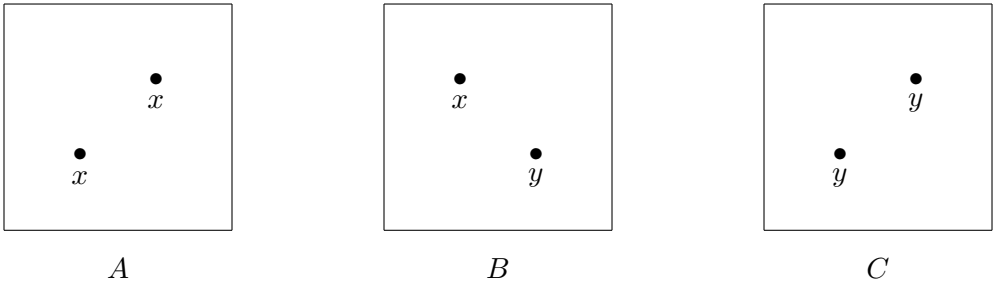
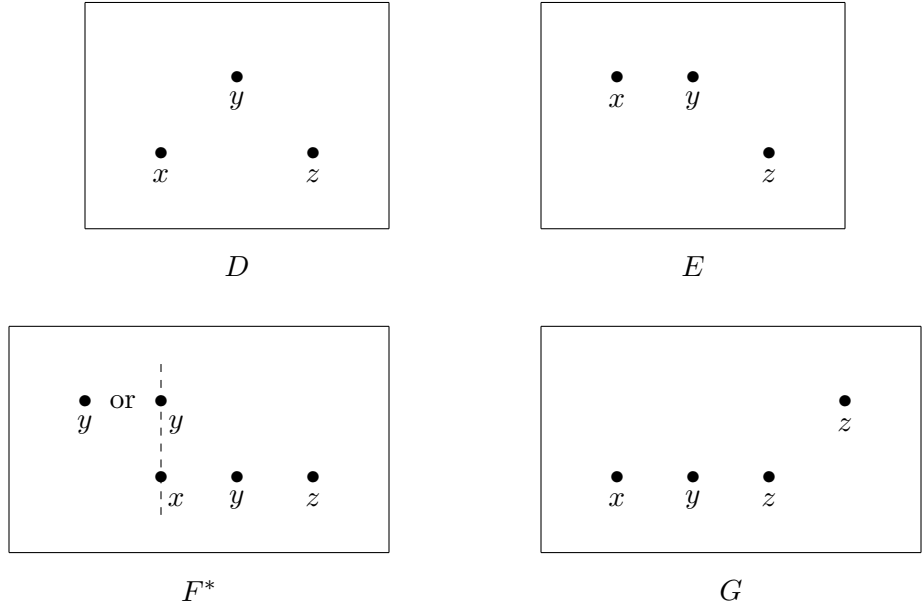
	\begin{figure}[h]
		\centering
		\begin{tikzpicture}
			\coordinate (a) at (0,0);
			\coordinate (b) at (3,0);
			\coordinate (c) at (0,3);
			\coordinate (d) at (3,3);
			\draw(a)--(b);
			\draw(b)--(d);
			\draw(d)--(c);
			\draw(c)--(a);
			\node(small-dot) at (1,1){$\bullet$};
			\node(small-dot) at (1,0.7){$x$};
			\node(small-dot) at (2,2){$\bullet$};   
			\node(small-dot) at (2,1.7){$x$};		
			\coordinate (a1) at (5,0);
			\coordinate (b1) at (8,0);
			\coordinate (c1) at (5,3);
			\coordinate (d1) at (8,3);
			\draw(a1)--(b1);
			\draw(b1)--(d1);
			\draw(d1)--(c1);
			\draw(c1)--(a1);
			\node(small-dot) at (6,2){$\bullet$};
			\node(small-dot) at (6,1.7){$x$};
			\node(small-dot) at (7,1){$\bullet$};   
			\node(small-dot) at (7,0.7){$y$};
			\coordinate (a2) at (10,0);
			\coordinate (b2) at (13,0);
			\coordinate (c2) at (10,3);
			\coordinate (d2) at (13,3);
			\draw(a2)--(b2);
			\draw(b2)--(d2);
			\draw(d2)--(c2);
			\draw(c2)--(a2);
			\node(small-dot) at (11,1){$\bullet$};
			\node(small-dot) at (11,0.7){$y$};
			\node(small-dot) at (12,2){$\bullet$};   
			\node(small-dot) at (12,1.7){$y$};
			\node(small-dot) at (1.5,-0.5){$A$};
			\node(small-dot) at (6.5,-0.5){$B$};
			\node(small-dot) at (11.5,-0.5){$C$};
		\end{tikzpicture}
	\end{figure}	\FloatBarrier
	\begin{figure}[h]
		\centering
		\begin{tikzpicture}
			\coordinate (a) at (0,0);
			\coordinate (b) at (4,0);
			\coordinate (c) at (0,3);
			\coordinate (d) at (4,3);
			\draw(a)--(b);
			\draw(b)--(d);
			\draw(d)--(c);
			\draw(c)--(a);
			\node(small-dot) at (1,1){$\bullet$};
			\node(small-dot) at (1,0.7){$x$};
			\node(small-dot) at (2,2){$\bullet$};
			\node(small-dot) at (2,1.7){$y$};
			\node(small-dot) at (3,1){$\bullet$};
			\node(small-dot) at (3,0.7){$z$};
			\coordinate (a1) at (6,0);
			\coordinate (b1) at (10,0);
			\coordinate (c1) at (6,3);
			\coordinate (d1) at (10,3);
			\draw(a1)--(b1);
			\draw(b1)--(d1);
			\draw(d1)--(c1);
			\draw(c1)--(a1);
			\node(small-dot) at (7,2){$\bullet$};
			\node(small-dot) at (7,1.7){$x$};
			\node(small-dot) at (8,2){$\bullet$};
			\node(small-dot) at (8,1.7){$y$};
			\node(small-dot) at (9,1){$\bullet$};
			\node(small-dot) at (9,0.7){$z$};
			\node(small-dot) at (2,-0.5){$D$};
			\node(small-dot) at (8,-0.5){$E$};
		\end{tikzpicture}
	\end{figure}\FloatBarrier
	\begin{figure}[h]
		\centering
		\begin{tikzpicture}
			\coordinate (a) at (0,0);
			\coordinate (b) at (5,0);
			\coordinate (c) at (0,3);
			\coordinate (d) at (5,3);
			\draw(a)--(b);
			\draw(b)--(d);
			\draw(d)--(c);
			\draw(c)--(a);
			\node(small-dot) at (1,2){$\bullet$};
			\node(small-dot) at (1,1.7){$y$};
			\node(small-dot) at (2,2){$\bullet$};
			\node(small-dot) at (2.2,1.7){$y$};
			\node(small-dot) at (2,1){$\bullet$};
			\node(small-dot) at (2.2,0.7){$x$};
			\node(small-dot) at (3,1){$\bullet$};
			\node(small-dot) at (3,0.7){$y$};
			\node(small-dot) at (4,1){$\bullet$};
			\node(small-dot) at (4,0.7){$z$};
			\node(small-dot) at (1.5,2){or};
			\coordinate (aa1) at (2,2.5);
			\coordinate (aa2) at (2,0.5);
			\draw[dashed](aa1)--(aa2);
			\coordinate (a1) at (7,0);
			\coordinate (b1) at (12,0);
			\coordinate (c1) at (7,3);
			\coordinate (d1) at (12,3);
			\draw(a1)--(b1);
			\draw(b1)--(d1);
			\draw(d1)--(c1);
			\draw(c1)--(a1);
			\node(small-dot) at (8,1){$\bullet$};
			\node(small-dot) at (8,0.7){$x$};
			\node(small-dot) at (9,1){$\bullet$};
			\node(small-dot) at (9,0.7){$y$};
			\node(small-dot) at (10,1){$\bullet$};
			\node(small-dot) at (10,0.7){$z$};
			\node(small-dot) at (11,2){$\bullet$};
			\node(small-dot) at (11,1.7){$z$};
			\node(small-dot) at (2.5,-0.5){$F^{*}$};
			\node(small-dot) at (9.5,-0.5){$G$};
		\end{tikzpicture}
	\end{figure}	\FloatBarrier
	\par The ideal $L(I^{2,n}_{2,2})$ is radical since leading terms of polynomials in $\Gamma$ are square-free. Hence it corresponds to a unique abstract simplicial complex $\Delta_{L(I^{2,n}_{2,2})}$, denoted by $\Delta _{0}$ for convenience. By definition, conditions of facets are restricted by $A, B, C, D, E, F^{*}, G$, that is, a facet contains none of the sets in $A,\cdots , G$. 
	\par Let $F$ be a facet of $\Delta_{0}$. We will consider $F$ as a polynomial, a set, and a diagram. As a polynomial, we have $F=F_{x}F_{y}F_{z}$ where $F_{x}, F_{y}, F_{z}$ are respectively composed of $x_{i,j}, y_{i,j}, z_{i,j}$.\\\\
	\textbf{Example 4.1. }We describe $F_{x}$-part, $F_{y}$-part, and $F_{z}$-part for $F=\{x_{1,1}, x_{1,2}, x_{2,2}, x_{2,3}, y_{1,3}, z_{2,1}, z_{2,2}, z_{2,3}\}$ with following diagram:\\
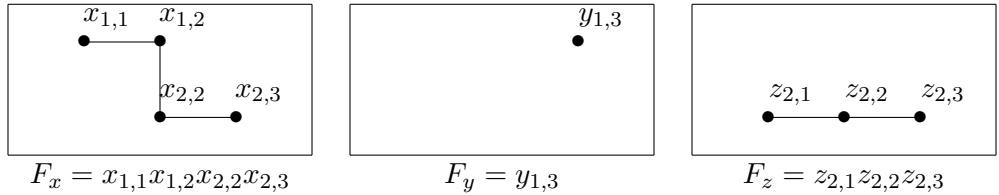
	\begin{figure}[h]
		\centering
		\begin{tikzpicture}
			\coordinate (a) at (0,0.5);
			\coordinate (b) at (4,0.5);
			\coordinate (c) at (0,2.5);
			\coordinate (d) at (4,2.5);
			\draw(a)--(b);
			\draw(b)--(d);
			\draw(d)--(c);
			\draw(c)--(a);
			\coordinate (p1) at (1,2);
			\coordinate (p2) at (2,2);
			\coordinate (p3) at (2,1);
			\coordinate (p4) at (3,1);
			\draw(p1)--(p2);
			\draw(p2)--(p3);
			\draw(p3)--(p4);
			\node(small-dot) at (1,2){$\bullet$};
			\node(small-dot) at (2,2){$\bullet$};
			\node(small-dot) at (2,1){$\bullet$};
			\node(small-dot) at (3,1){$\bullet$};
			\node(small-dot) at (1.3,2.3){$x_{1,1}$};
			\node(small-dot) at (2.3,2.3){$x_{1,2}$};
			\node(small-dot) at (2.3,1.3){$x_{2,2}$};
			\node(small-dot) at (3.3,1.3){$x_{2,3}$};
			\node(small-dot) at (2,0.2){$F_{x}=x_{1,1}x_{1,2}x_{2,2}x_{2,3}$};
			\coordinate (a1) at (4.5,0.5);
			\coordinate (b1) at (8.5,0.5);
			\coordinate (c1) at (4.5,2.5);
			\coordinate (d1) at (8.5,2.5);
			\draw(a1)--(b1);
			\draw(b1)--(d1);
			\draw(d1)--(c1);
			\draw(c1)--(a1);
			\node(small-dot) at (7.5,2){$\bullet$};
			\node(small-dot) at (7.8,2.3){$y_{1,3}$};
			\node(small-dot) at (6.5,0.2){$F_{y}=y_{1,3}$};
			\coordinate (a2) at (9,0.5);
			\coordinate (b2) at (13,0.5);
			\coordinate (c2) at (9,2.5);
			\coordinate (d2) at (13,2.5);
			\draw(a2)--(b2);
			\draw(b2)--(d2);
			\draw(d2)--(c2);
			\draw(c2)--(a2);
			\coordinate (ppp1) at (10,1);
			\coordinate (ppp2) at (12,1);
			\draw(ppp1)--(ppp2);
            \node(small-dot) at (10,1){$\bullet$};
            \node(small-dot) at (11,1){$\bullet$};
            \node(small-dot) at (12,1){$\bullet$};
            \node(small-dot) at (10.3,1.3){$z_{2,1}$};
            \node(small-dot) at (11.3,1.3){$z_{2,2}$};
            \node(small-dot) at (12.3,1.3){$z_{2,3}$};
			\node(small-dot) at (11,0.2){$F_{z}=z_{2,1}z_{2,2}z_{2,3}$};
		\end{tikzpicture}
	\end{figure}	\FloatBarrier
	\par We will use such diagrams as above later in this paper, so detailed explanation is necessary. A line segment along one row is closed, that is, it contains all points including the endpoints on itself. A line segment along one column contains two points, but it is mainly used to express "paths walking only to the right or down" as stated below.
	\par Restriction of $A$ shows that $F_{x}$ is a subset of a path walking continuously from $x_{1,1}$ to $x_{2,n}$ only to the right or down. Similarly, considering $F_{y}$ and $C$, one see that $F_{y}$ is a subset of a path walking continuously from $y_{1,1}$ to $y_{2,n}$ only to the right or down.
	\par We first consider the case that the first row of $F_{x}$ is nonempty. The result is as follows: \\\\
	\textbf{Proposition 4.2. }\textit{Let $F$ be a facet of $\Delta_{0}$. If the first row of $F_{x}$ is nonempty, all possible diagrams of $F$ are as follows:}\\
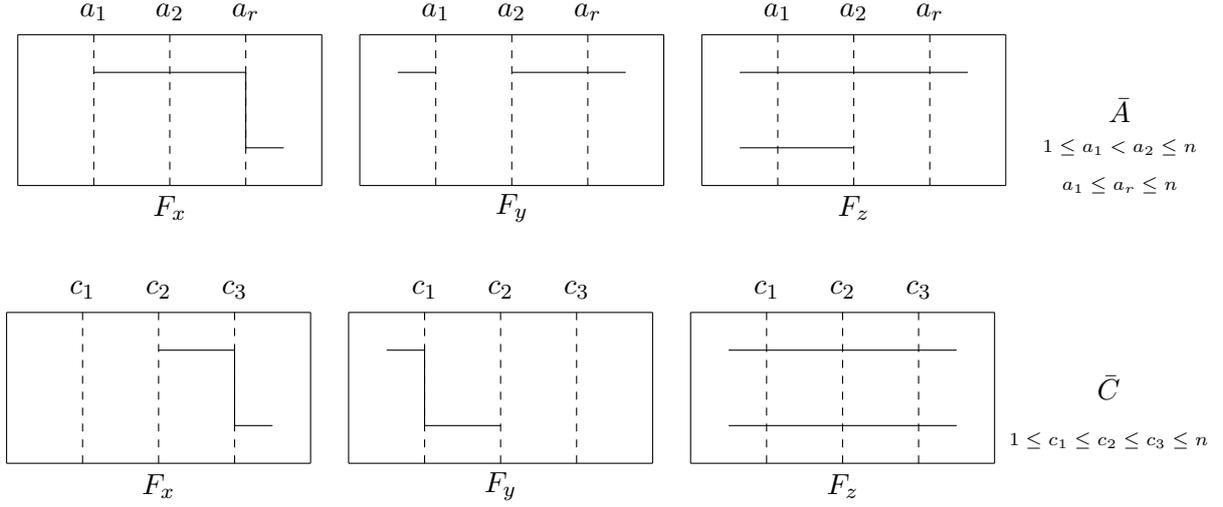
	\begin{figure}[h]
		\centering
		\begin{tikzpicture}
			\coordinate (a) at (0,0.5);
			\coordinate (b) at (4,0.5);
			\coordinate (c) at (0,2.5);
			\coordinate (d) at (4,2.5);
			\draw(a)--(b);
			\draw(b)--(d);
			\draw(d)--(c);
			\draw(c)--(a);
			\coordinate (x1) at (1,0.5);
			\coordinate (y1) at (1,2.5);
			\coordinate (x2) at (2,0.5);
			\coordinate (y2) at (2,2.5);
			\coordinate (x3) at (3,0.5);
			\coordinate (y3) at (3,2.5);
			\node(small-dot) at (1,2.8){$a_{1}$};
			\node(small-dot) at (2,2.8){$a_{2}$};
			\node(small-dot) at (3,2.8){$a_{r}$};
			\draw[dashed](x1)--(y1);
			\draw[dashed](x2)--(y2);
			\draw[dashed](x3)--(y3);
			\coordinate (p1) at (1,2);
			\coordinate (p2) at (3,2);
			\coordinate (p3) at (3,1);
			\coordinate (p4) at (3.5,1);
			\draw(p1)--(p2);
			\draw(p2)--(p3);
			\draw(p3)--(p4);
			\node(small-dot) at (2,0.2){$F_{x}$};
			\coordinate (a1) at (4.5,0.5);
			\coordinate (b1) at (8.5,0.5);
			\coordinate (c1) at (4.5,2.5);
			\coordinate (d1) at (8.5,2.5);
			\draw(a1)--(b1);
			\draw(b1)--(d1);
			\draw(d1)--(c1);
			\draw(c1)--(a1);
			\coordinate (x4) at (5.5,0.5);
			\coordinate (y4) at (5.5,2.5);
			\coordinate (x5) at (6.5,0.5);
			\coordinate (y5) at (6.5,2.5);
			\coordinate (x6) at (7.5,0.5);
			\coordinate (y6) at (7.5,2.5);
			\node(small-dot) at (5.5,2.8){$a_{1}$};
			\node(small-dot) at (6.5,2.8){$a_{2}$};
			\node(small-dot) at (7.5,2.8){$a_{r}$};
			\draw[dashed](x4)--(y4);
			\draw[dashed](x5)--(y5);
			\draw[dashed](x6)--(y6);
			\coordinate (pp1) at (5,2);
			\coordinate (pp2) at (5.5,2);
			\coordinate (pp3) at (6.5,2);
			\coordinate (pp4) at (8,2);
			\draw(pp1)--(pp2);
			\draw(pp3)--(pp4);
			\node(small-dot) at (6.5,0.2){$F_{y}$};
			\coordinate (a2) at (9,0.5);
			\coordinate (b2) at (13,0.5);
			\coordinate (c2) at (9,2.5);
			\coordinate (d2) at (13,2.5);
			\draw(a2)--(b2);
			\draw(b2)--(d2);
			\draw(d2)--(c2);
			\draw(c2)--(a2);
			\coordinate (x7) at (10,0.5);
			\coordinate (y7) at (10,2.5);
			\coordinate (x8) at (11,0.5);
			\coordinate (y8) at (11,2.5);
			\coordinate (x9) at (12,0.5);
			\coordinate (y9) at (12,2.5);
			\node(small-dot) at (10,2.8){$a_{1}$};
			\node(small-dot) at (11,2.8){$a_{2}$};
			\node(small-dot) at (12,2.8){$a_{r}$};
			\draw[dashed](x7)--(y7);
			\draw[dashed](x8)--(y8);
			\draw[dashed](x9)--(y9);
			\coordinate (ppp1) at (9.5,2);
			\coordinate (ppp2) at (12.5,2);
			\coordinate (ppp3) at (9.5,1);
			\coordinate (ppp4) at (11,1);
			\draw(ppp1)--(ppp2);
			\draw(ppp3)--(ppp4);
			\node(small-dot) at (11,0.2){$F_{z}$};
			\node(small-dot) at (14.5,1.5){$\bar{A}$};
			\node(small-dot) at (14.5,1){\tiny{$1\le a_{1}< a_{2}\le n$}};
			\node(small-dot) at (14.5,0.5){\tiny{$a_{1}\le a_{r}\le n$}};
		\end{tikzpicture}
	\end{figure}	\FloatBarrier
	\begin{figure}[h]
		\centering
		\begin{tikzpicture}
			\coordinate (a) at (0,0.5);
			\coordinate (b) at (4,0.5);
			\coordinate (c) at (0,2.5);
			\coordinate (d) at (4,2.5);
			\draw(a)--(b);
			\draw(b)--(d);
			\draw(d)--(c);
			\draw(c)--(a);
			\coordinate (x1) at (1,0.5);
			\coordinate (y1) at (1,2.5);
			\coordinate (x2) at (2,0.5);
			\coordinate (y2) at (2,2.5);
			\coordinate (x3) at (3,0.5);
			\coordinate (y3) at (3,2.5);
			\node(small-dot) at (1,2.8){$c_{1}$};
			\node(small-dot) at (2,2.8){$c_{2}$};
			\node(small-dot) at (3,2.8){$c_{3}$};
			\draw[dashed](x1)--(y1);
			\draw[dashed](x2)--(y2);
			\draw[dashed](x3)--(y3);
			\coordinate (p1) at (2,2);
			\coordinate (p2) at (3,2);
			\coordinate (p3) at (3,1);
			\coordinate (p4) at (3.5,1);
			\draw(p1)--(p2);
			\draw(p2)--(p3);
			\draw(p3)--(p4);
			\node(small-dot) at (2,0.2){$F_{x}$};
			\coordinate (a1) at (4.5,0.5);
			\coordinate (b1) at (8.5,0.5);
			\coordinate (c1) at (4.5,2.5);
			\coordinate (d1) at (8.5,2.5);
			\draw(a1)--(b1);
			\draw(b1)--(d1);
			\draw(d1)--(c1);
			\draw(c1)--(a1);
			\coordinate (x4) at (5.5,0.5);
			\coordinate (y4) at (5.5,2.5);
			\coordinate (x5) at (6.5,0.5);
			\coordinate (y5) at (6.5,2.5);
			\coordinate (x6) at (7.5,0.5);
			\coordinate (y6) at (7.5,2.5);
			\node(small-dot) at (5.5,2.8){$c_{1}$};
			\node(small-dot) at (6.5,2.8){$c_{2}$};
			\node(small-dot) at (7.5,2.8){$c_{3}$};
			\draw[dashed](x4)--(y4);
			\draw[dashed](x5)--(y5);
			\draw[dashed](x6)--(y6);
			\coordinate (pp1) at (5,2);
			\coordinate (pp2) at (5.5,2);
			\coordinate (pp3) at (5.5,1);
			\coordinate (pp4) at (6.5,1);
			\draw(pp1)--(pp2);
			\draw(pp2)--(pp3);
			\draw(pp3)--(pp4);
			\node(small-dot) at (6.5,0.2){$F_{y}$};
			\coordinate (a2) at (9,0.5);
			\coordinate (b2) at (13,0.5);
			\coordinate (c2) at (9,2.5);
			\coordinate (d2) at (13,2.5);
			\draw(a2)--(b2);
			\draw(b2)--(d2);
			\draw(d2)--(c2);
			\draw(c2)--(a2);
			\coordinate (x7) at (10,0.5);
			\coordinate (y7) at (10,2.5);
			\coordinate (x8) at (11,0.5);
			\coordinate (y8) at (11,2.5);
			\coordinate (x9) at (12,0.5);
			\coordinate (y9) at (12,2.5);
			\node(small-dot) at (10,2.8){$c_{1}$};
			\node(small-dot) at (11,2.8){$c_{2}$};
			\node(small-dot) at (12,2.8){$c_{3}$};
			\draw[dashed](x7)--(y7);
			\draw[dashed](x8)--(y8);
			\draw[dashed](x9)--(y9);
			\coordinate (ppp1) at (9.5,2);
			\coordinate (ppp2) at (12.5,2);
			\coordinate (ppp3) at (9.5,1);
			\coordinate (ppp4) at (12.5,1);
			\draw(ppp1)--(ppp2);
			\draw(ppp3)--(ppp4);
			\node(small-dot) at (11,0.2){$F_{z}$};
			\node(small-dot) at (14.5,1.5){$\bar{C}$};
			\node(small-dot) at (14.5,0.8){\tiny{$1\le c_{1}\le c_{2}\le c_{3}\le n$}};
		\end{tikzpicture}
	\end{figure}	\FloatBarrier
	\begin{proof}[Proof.]
		It is straightforward to verify that $F$ with diagram of the type $\bar{A}$ or $\bar{C}$ is indeed a facet, that is, $F$ contains no forbidden sets in $A,\cdots ,G$ and admits no other points. Suppose that $s$ is the maximal column index such that $x_{1,s}\in F_{x}$. Consider the following two diagrams:\\
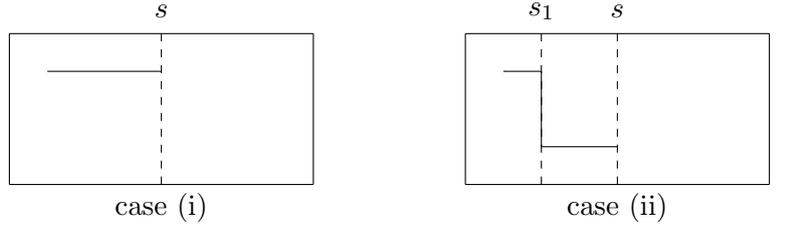
		\begin{figure}[h]
			\centering
			\begin{tikzpicture}
				\coordinate (a) at (0,0.5);
				\coordinate (b) at (4,0.5);
				\coordinate (c) at (0,2.5);
				\coordinate (d) at (4,2.5);
				\draw(a)--(b);
				\draw(b)--(d);
				\draw(d)--(c);
				\draw(c)--(a);
				\coordinate (p1) at (0.5,2);
				\coordinate (p2) at (2,2);
				\coordinate (a1) at (2,0.5);
				\coordinate (b1) at (2,2.5);
				\draw[dashed](a1)--(b1);
				\draw(p1)--(p2);
				\node(small-dot) at (2,2.8){$s$};
				\node(small-dot) at (2,0.2){case (i)};
				\coordinate (a2) at (6,0.5);
				\coordinate (b2) at (10,0.5);
				\coordinate (c2) at (6,2.5);
				\coordinate (d2) at (10,2.5);
				\draw(a2)--(b2);
				\draw(b2)--(d2);
				\draw(d2)--(c2);
				\draw(c2)--(a2);
				\coordinate (p3) at (6.5,2);
				\coordinate (p4) at (7,2);
				\coordinate (p5) at (7,1);
				\coordinate (p6) at (8,1);
				\coordinate (a3) at (8,0.5);
				\coordinate (b3) at (8,2.5);
				\coordinate (a4) at (7,0.5);
				\coordinate (b4) at (7,2.5);
				\draw[dashed](a3)--(b3);
				\draw[dashed](a4)--(b4);
				\draw(p3)--(p4);
				\draw(p4)--(p5);
				\draw(p5)--(p6);
				\node(small-dot) at (8,2.8){$s$};
				\node(small-dot) at (7,2.8){$s_{1}$};
				\node(small-dot) at (8,0.2){case (ii)};
				\node(small-dot) at (10.5,1.5){,};
			\end{tikzpicture}
		\end{figure}    \FloatBarrier
		\par \noindent The $\tiny{\le} s$-part of $F_{y}$ must be a subset of one set corresponding to case (i) or case (ii). By restriction of $B$, the second row of $\tiny{>}s$-part of $F_{y}$ is empty.
		\par In case (i), $F_{y}$ is empty in the second row. Therefore $F$ cannot contain any set in $F^{*}$ or $G$. 
		\par If the $\tiny{>}s$-part of $F_{y}$ is nonempty, we suppose that $t$ is the minimal column index such that $y_{1,t}\in F_{y}$ and $s<t$. Let $t_{1}$ be the maximal column index such that $z_{2,t_{1}}\in F_{z}$.\\\\
		\textit{Claim }1. If $t_{1}\le s$ and $s\neq n$, then $t=s+1$. 
		\begin{proof}[Proof of Claim]
			It is obvious that $F$ cannot contain any set in $D$ or $E$ when $t_{1}\le s$. Suppose $t>s+1$ and $F_{0}$ is a facet for $s,t,t_{1}$ in such condition. Let $F^{'}=F_{0}y_{1,s+1}$. Then $F^{'}\in \Delta_{0}$ and this yields a contradiction. 
		\end{proof}
		\par \noindent \textit{Claim }2. If $t_{1}>s$, then $t=t_{1}$.
		\begin{proof}[Proof of Claim]
			The subset $\{x_{1,s},y_{1,t},z_{2,t_{1}}\} \subset F$ cannot be of the form as diagram $E$. Hence $t\ge t_{1}$ and then $F$ cannot contain any set in $D$ or $E$. Suppose $t>t_{1}$ and $F_{0}$ is a facet for $s,t,t_{1}$ in such condition. Let $F^{'}=F_{0}z_{2,t}$. Then $F^{'}\in \Delta_{0}$ and this yields a contradiction.
		\end{proof}
		\par Let $t_{2}$ be the minimal column index such that $x_{1,t_{2}}\in F_{x}$. Then we have two possible diagrams corresponding to the case $t_{1}\le s$ and the case $t_{1}>s$:\\
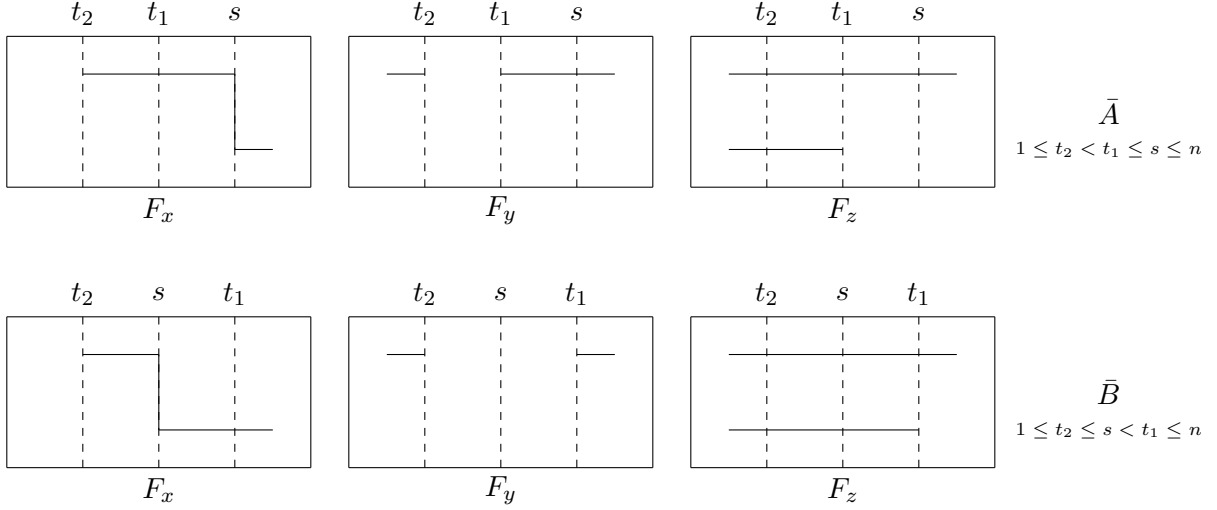
		\begin{figure}[h]
			\centering
			\begin{tikzpicture}
				\coordinate (a) at (0,0.5);
				\coordinate (b) at (4,0.5);
				\coordinate (c) at (0,2.5);
				\coordinate (d) at (4,2.5);
				\draw(a)--(b);
				\draw(b)--(d);
				\draw(d)--(c);
				\draw(c)--(a);
				\coordinate (x1) at (1,0.5);
				\coordinate (y1) at (1,2.5);
				\coordinate (x2) at (2,0.5);
				\coordinate (y2) at (2,2.5);
				\coordinate (x3) at (3,0.5);
				\coordinate (y3) at (3,2.5);
				\node(small-dot) at (1,2.8){$t_{2}$};
				\node(small-dot) at (2,2.8){$t_{1}$};
				\node(small-dot) at (3,2.8){$s$};
				\draw[dashed](x1)--(y1);
				\draw[dashed](x2)--(y2);
				\draw[dashed](x3)--(y3);
				\coordinate (p1) at (1,2);
				\coordinate (p2) at (3,2);
				\coordinate (p3) at (3,1);
				\coordinate (p4) at (3.5,1);
				\draw(p1)--(p2);
				\draw(p2)--(p3);
				\draw(p3)--(p4);
				\node(small-dot) at (2,0.2){$F_{x}$};
				\coordinate (a1) at (4.5,0.5);
				\coordinate (b1) at (8.5,0.5);
				\coordinate (c1) at (4.5,2.5);
				\coordinate (d1) at (8.5,2.5);
				\draw(a1)--(b1);
				\draw(b1)--(d1);
				\draw(d1)--(c1);
				\draw(c1)--(a1);
				\coordinate (x4) at (5.5,0.5);
				\coordinate (y4) at (5.5,2.5);
				\coordinate (x5) at (6.5,0.5);
				\coordinate (y5) at (6.5,2.5);
				\coordinate (x6) at (7.5,0.5);
				\coordinate (y6) at (7.5,2.5);
				\node(small-dot) at (5.5,2.8){$t_{2}$};
				\node(small-dot) at (6.5,2.8){$t_{1}$};
				\node(small-dot) at (7.5,2.8){$s$};
				\draw[dashed](x4)--(y4);
				\draw[dashed](x5)--(y5);
				\draw[dashed](x6)--(y6);
				\coordinate (pp1) at (5,2);
				\coordinate (pp2) at (5.5,2);
				\coordinate (pp3) at (6.5,2);
				\coordinate (pp4) at (8,2);
				\draw(pp1)--(pp2);
				\draw(pp3)--(pp4);
				\node(small-dot) at (6.5,0.2){$F_{y}$};
				\coordinate (a2) at (9,0.5);
				\coordinate (b2) at (13,0.5);
				\coordinate (c2) at (9,2.5);
				\coordinate (d2) at (13,2.5);
				\draw(a2)--(b2);
				\draw(b2)--(d2);
				\draw(d2)--(c2);
				\draw(c2)--(a2);
				\coordinate (x7) at (10,0.5);
				\coordinate (y7) at (10,2.5);
				\coordinate (x8) at (11,0.5);
				\coordinate (y8) at (11,2.5);
				\coordinate (x9) at (12,0.5);
				\coordinate (y9) at (12,2.5);
				\node(small-dot) at (10,2.8){$t_{2}$};
				\node(small-dot) at (11,2.8){$t_{1}$};
				\node(small-dot) at (12,2.8){$s$};
				\draw[dashed](x7)--(y7);
				\draw[dashed](x8)--(y8);
				\draw[dashed](x9)--(y9);
				\coordinate (ppp1) at (9.5,2);
				\coordinate (ppp2) at (12.5,2);
				\coordinate (ppp3) at (9.5,1);
				\coordinate (ppp4) at (11,1);
				\draw(ppp1)--(ppp2);
				\draw(ppp3)--(ppp4);
				\node(small-dot) at (11,0.2){$F_{z}$};
				\node(small-dot) at (14.5,1.5){$\bar{A}$};
				\node(small-dot) at (14.5,1){\tiny{$1\le t_{2}< t_{1}\le s\le n$}};
			\end{tikzpicture}
		\end{figure}	\FloatBarrier
		\begin{figure}[h]
			\centering
			\begin{tikzpicture}
				\coordinate (a) at (0,0.5);
				\coordinate (b) at (4,0.5);
				\coordinate (c) at (0,2.5);
				\coordinate (d) at (4,2.5);
				\draw(a)--(b);
				\draw(b)--(d);
				\draw(d)--(c);
				\draw(c)--(a);
				\coordinate (x1) at (1,0.5);
				\coordinate (y1) at (1,2.5);
				\coordinate (x2) at (2,0.5);
				\coordinate (y2) at (2,2.5);
				\coordinate (x3) at (3,0.5);
				\coordinate (y3) at (3,2.5);
				\node(small-dot) at (1,2.8){$t_{2}$};
				\node(small-dot) at (2,2.8){$s$};
				\node(small-dot) at (3,2.8){$t_{1}$};
				\draw[dashed](x1)--(y1);
				\draw[dashed](x2)--(y2);
				\draw[dashed](x3)--(y3);
				\coordinate (p1) at (1,2);
				\coordinate (p2) at (2,2);
				\coordinate (p3) at (2,1);
				\coordinate (p4) at (3.5,1);
				\draw(p1)--(p2);
				\draw(p2)--(p3);
				\draw(p3)--(p4);
				\node(small-dot) at (2,0.2){$F_{x}$};
				\coordinate (a1) at (4.5,0.5);
				\coordinate (b1) at (8.5,0.5);
				\coordinate (c1) at (4.5,2.5);
				\coordinate (d1) at (8.5,2.5);
				\draw(a1)--(b1);
				\draw(b1)--(d1);
				\draw(d1)--(c1);
				\draw(c1)--(a1);
				\coordinate (x4) at (5.5,0.5);
				\coordinate (y4) at (5.5,2.5);
				\coordinate (x5) at (6.5,0.5);
				\coordinate (y5) at (6.5,2.5);
				\coordinate (x6) at (7.5,0.5);
				\coordinate (y6) at (7.5,2.5);
				\node(small-dot) at (5.5,2.8){$t_{2}$};
				\node(small-dot) at (6.5,2.8){$s$};
				\node(small-dot) at (7.5,2.8){$t_{1}$};
				\draw[dashed](x4)--(y4);
				\draw[dashed](x5)--(y5);
				\draw[dashed](x6)--(y6);
				\coordinate (pp1) at (5,2);
				\coordinate (pp2) at (5.5,2);
				\coordinate (pp3) at (7.5,2);
				\coordinate (pp4) at (8,2);
				\draw(pp1)--(pp2);
				\draw(pp3)--(pp4);
				\node(small-dot) at (6.5,0.2){$F_{y}$};
				\coordinate (a2) at (9,0.5);
				\coordinate (b2) at (13,0.5);
				\coordinate (c2) at (9,2.5);
				\coordinate (d2) at (13,2.5);
				\draw(a2)--(b2);
				\draw(b2)--(d2);
				\draw(d2)--(c2);
				\draw(c2)--(a2);
				\coordinate (x7) at (10,0.5);
				\coordinate (y7) at (10,2.5);
				\coordinate (x8) at (11,0.5);
				\coordinate (y8) at (11,2.5);
				\coordinate (x9) at (12,0.5);
				\coordinate (y9) at (12,2.5);
				\node(small-dot) at (10,2.8){$t_{2}$};
				\node(small-dot) at (11,2.8){$s$};
				\node(small-dot) at (12,2.8){$t_{1}$};
				\draw[dashed](x7)--(y7);
				\draw[dashed](x8)--(y8);
				\draw[dashed](x9)--(y9);
				\coordinate (ppp1) at (9.5,2);
				\coordinate (ppp2) at (12.5,2);
				\coordinate (ppp3) at (9.5,1);
				\coordinate (ppp4) at (12,1);
				\draw(ppp1)--(ppp2);
				\draw(ppp3)--(ppp4);
				\node(small-dot) at (11,0.2){$F_{z}$};
				\node(small-dot) at (14.5,1.5){$\bar{B}$};
				\node(small-dot) at (14.5,1){\tiny{$1\le t_{2}\le s< t_{1}\le n$}};
			\end{tikzpicture}
		\end{figure}	\FloatBarrier
		\par We can reduce $\bar{A}$ and $\bar{B}$ to the same condition $\bar{A}$ as stated in the proposition.
		\par If the $\tiny{>}s$-part of $F_{y}$ is empty, then $F$ cannot contain any set in $D$. Let $t$ be the maximal column index such that $y_{1,t}\in F_{y}$. Suppose $F_{0}$ is a facet for $s,t$ in such condition. Let $F^{'}=F_{0}y_{2,t}$. Then $F^{'}\in \Delta_{0}$ and this yields a contradiction.
		\par In case (ii), we suppose that the second row of $F_{y}$ is nonempty, otherwise the condition can be reduced to case (i). Any set in $D, F^{*}, G$ cannot be contained in $F$ because $F_{x}$ walks down after $s$. Let $s_{1}$ be the minimal column index such that $y_{2,s_{1}}\in F_{y}$ and $t$ be the maximal column index such that $z_{2,t}\in F_{z}$. The first row of $\tiny{<}s_{1}$-part of $F_{x}$ is empty by restriction of $B$. Hence any set in $E$ cannot be contained in $F$. We have not given any restriction to $t$ and restrictions from $A,B,C$ have no relation with $F_{z}$. Hence $t=n$ because $F$ is maximal in $\Delta_{0}$. Let $t_{1}$ be the minimal column index such that $x_{1,t_{1}}\in F_{x}$. We have diagrams for this case as follows:\\
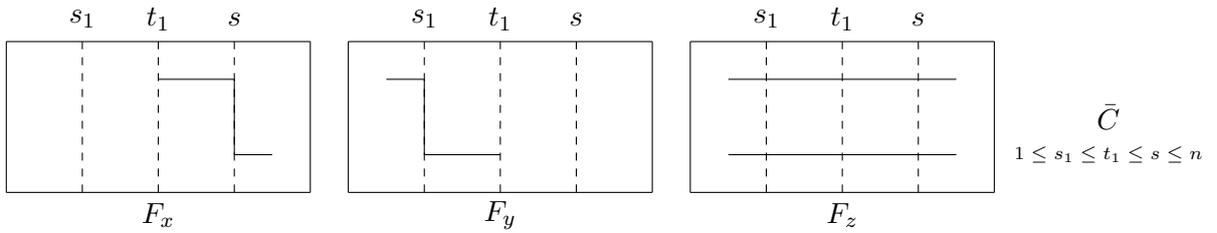
		\begin{figure}[!h]
			\centering
			\begin{tikzpicture}
				\coordinate (a) at (0,0.5);
				\coordinate (b) at (4,0.5);
				\coordinate (c) at (0,2.5);
				\coordinate (d) at (4,2.5);
				\draw(a)--(b);
				\draw(b)--(d);
				\draw(d)--(c);
				\draw(c)--(a);
				\coordinate (x1) at (1,0.5);
				\coordinate (y1) at (1,2.5);
				\coordinate (x2) at (2,0.5);
				\coordinate (y2) at (2,2.5);
				\coordinate (x3) at (3,0.5);
				\coordinate (y3) at (3,2.5);
				\node(small-dot) at (1,2.8){$s_{1}$};
				\node(small-dot) at (2,2.8){$t_{1}$};
				\node(small-dot) at (3,2.8){$s$};
				\draw[dashed](x1)--(y1);
				\draw[dashed](x2)--(y2);
				\draw[dashed](x3)--(y3);
				\coordinate (p1) at (2,2);
				\coordinate (p2) at (3,2);
				\coordinate (p3) at (3,1);
				\coordinate (p4) at (3.5,1);
				\draw(p1)--(p2);
				\draw(p2)--(p3);
				\draw(p3)--(p4);
				\node(small-dot) at (2,0.2){$F_{x}$};
				\coordinate (a1) at (4.5,0.5);
				\coordinate (b1) at (8.5,0.5);
				\coordinate (c1) at (4.5,2.5);
				\coordinate (d1) at (8.5,2.5);
				\draw(a1)--(b1);
				\draw(b1)--(d1);
				\draw(d1)--(c1);
				\draw(c1)--(a1);
				\coordinate (x4) at (5.5,0.5);
				\coordinate (y4) at (5.5,2.5);
				\coordinate (x5) at (6.5,0.5);
				\coordinate (y5) at (6.5,2.5);
				\coordinate (x6) at (7.5,0.5);
				\coordinate (y6) at (7.5,2.5);
				\node(small-dot) at (5.5,2.8){$s_{1}$};
				\node(small-dot) at (6.5,2.8){$t_{1}$};
				\node(small-dot) at (7.5,2.8){$s$};
				\draw[dashed](x4)--(y4);
				\draw[dashed](x5)--(y5);
				\draw[dashed](x6)--(y6);
				\coordinate (pp1) at (5,2);
				\coordinate (pp2) at (5.5,2);
				\coordinate (pp3) at (5.5,1);
				\coordinate (pp4) at (6.5,1);
				\draw(pp1)--(pp2);
				\draw(pp2)--(pp3);
				\draw(pp3)--(pp4);
				\node(small-dot) at (6.5,0.2){$F_{y}$};
				\coordinate (a2) at (9,0.5);
				\coordinate (b2) at (13,0.5);
				\coordinate (c2) at (9,2.5);
				\coordinate (d2) at (13,2.5);
				\draw(a2)--(b2);
				\draw(b2)--(d2);
				\draw(d2)--(c2);
				\draw(c2)--(a2);
				\coordinate (x7) at (10,0.5);
				\coordinate (y7) at (10,2.5);
				\coordinate (x8) at (11,0.5);
				\coordinate (y8) at (11,2.5);
				\coordinate (x9) at (12,0.5);
				\coordinate (y9) at (12,2.5);
				\node(small-dot) at (10,2.8){$s_{1}$};
				\node(small-dot) at (11,2.8){$t_{1}$};
				\node(small-dot) at (12,2.8){$s$};
				\draw[dashed](x7)--(y7);
				\draw[dashed](x8)--(y8);
				\draw[dashed](x9)--(y9);
				\coordinate (ppp1) at (9.5,2);
				\coordinate (ppp2) at (12.5,2);
				\coordinate (ppp3) at (9.5,1);
				\coordinate (ppp4) at (12.5,1);
				\draw(ppp1)--(ppp2);
				\draw(ppp3)--(ppp4);
				\node(small-dot) at (11,0.2){$F_{z}$};
				\node(small-dot) at (14.5,1.5){$\bar{C}$};
				\node(small-dot) at (14.5,1){\tiny{$1\le s_{1}\le t_{1}\le s\le n$}};
			\end{tikzpicture}
		\end{figure}	\FloatBarrier
	\end{proof} 
	\par If $F_{x}$ is empty in the first row, then $F$ cannot contain any set in $B$ or $E$. We only have to consider restriction from $D, F^{*}, G$. The result is as follows:\\\\
	\textbf{Proposition 4.3. }\textit{Let $F$ be a facet of $\Delta_{0}$. If the first row of $F_{x}$ is empty, all possible diagrams of $F$ are as follows:}\\
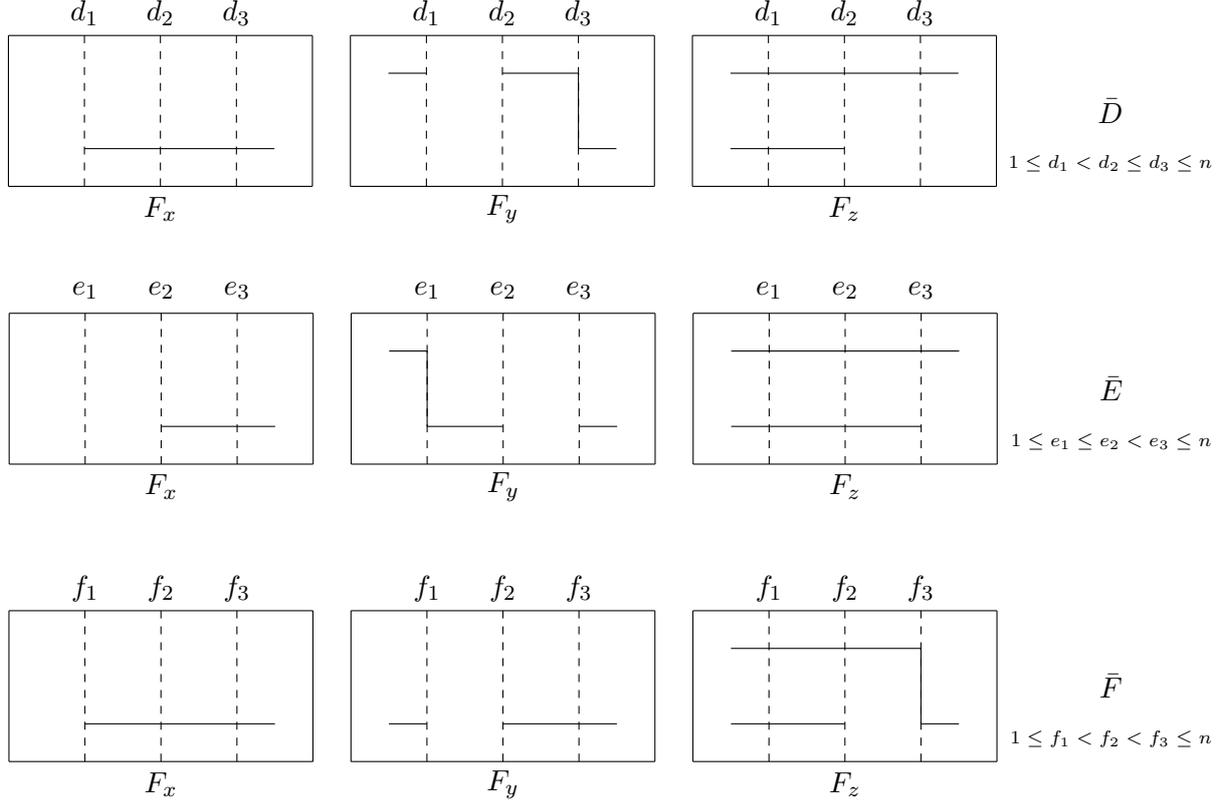
	\begin{figure}[!h]
		\centering
		\begin{tikzpicture}
			\coordinate (a) at (0,0.5);
			\coordinate (b) at (4,0.5);
			\coordinate (c) at (0,2.5);
			\coordinate (d) at (4,2.5);
			\draw(a)--(b);
			\draw(b)--(d);
			\draw(d)--(c);
			\draw(c)--(a);
			\coordinate (x1) at (1,0.5);
			\coordinate (y1) at (1,2.5);
			\coordinate (x2) at (2,0.5);
			\coordinate (y2) at (2,2.5);
			\coordinate (x3) at (3,0.5);
			\coordinate (y3) at (3,2.5);
			\node(small-dot) at (1,2.8){$d_{1}$};
			\node(small-dot) at (2,2.8){$d_{2}$};
			\node(small-dot) at (3,2.8){$d_{3}$};
			\draw[dashed](x1)--(y1);
			\draw[dashed](x2)--(y2);
			\draw[dashed](x3)--(y3);
			\coordinate (p1) at (1,1);
			\coordinate (p2) at (3.5,1);
			\draw(p1)--(p2);
			\node(small-dot) at (2,0.2){$F_{x}$};
			\coordinate (a1) at (4.5,0.5);
			\coordinate (b1) at (8.5,0.5);
			\coordinate (c1) at (4.5,2.5);
			\coordinate (d1) at (8.5,2.5);
			\draw(a1)--(b1);
			\draw(b1)--(d1);
			\draw(d1)--(c1);
			\draw(c1)--(a1);
			\coordinate (x4) at (5.5,0.5);
			\coordinate (y4) at (5.5,2.5);
			\coordinate (x5) at (6.5,0.5);
			\coordinate (y5) at (6.5,2.5);
			\coordinate (x6) at (7.5,0.5);
			\coordinate (y6) at (7.5,2.5);
			\node(small-dot) at (5.5,2.8){$d_{1}$};
			\node(small-dot) at (6.5,2.8){$d_{2}$};
			\node(small-dot) at (7.5,2.8){$d_{3}$};
			\draw[dashed](x4)--(y4);
			\draw[dashed](x5)--(y5);
			\draw[dashed](x6)--(y6);
			\coordinate (pp1) at (5,2);
			\coordinate (pp2) at (5.5,2);
			\coordinate (pp3) at (6.5,2);
			\coordinate (pp4) at (7.5,2);
			\coordinate (pp5) at (7.5,1);
			\coordinate (pp6) at (8,1);
			\draw(pp1)--(pp2);
			\draw(pp3)--(pp4);
			\draw(pp4)--(pp5);
			\draw(pp5)--(pp6);
			\node(small-dot) at (6.5,0.2){$F_{y}$};
			\coordinate (a2) at (9,0.5);
			\coordinate (b2) at (13,0.5);
			\coordinate (c2) at (9,2.5);
			\coordinate (d2) at (13,2.5);
			\draw(a2)--(b2);
			\draw(b2)--(d2);
			\draw(d2)--(c2);
			\draw(c2)--(a2);
			\coordinate (x7) at (10,0.5);
			\coordinate (y7) at (10,2.5);
			\coordinate (x8) at (11,0.5);
			\coordinate (y8) at (11,2.5);
			\coordinate (x9) at (12,0.5);
			\coordinate (y9) at (12,2.5);
			\node(small-dot) at (10,2.8){$d_{1}$};
			\node(small-dot) at (11,2.8){$d_{2}$};
			\node(small-dot) at (12,2.8){$d_{3}$};
			\draw[dashed](x7)--(y7);
			\draw[dashed](x8)--(y8);
			\draw[dashed](x9)--(y9);
			\coordinate (ppp1) at (9.5,2);
			\coordinate (ppp2) at (12.5,2);
			\coordinate (ppp3) at (9.5,1);
			\coordinate (ppp4) at (11,1);
			\draw(ppp1)--(ppp2);
			\draw(ppp3)--(ppp4);
			\node(small-dot) at (11,0.2){$F_{z}$};
			\node(small-dot) at (14.5,1.5){$\bar{D}$};
			\node(small-dot) at (14.5,0.8){\tiny{$1\le d_{1}< d_{2}\le d_{3}\le n$}};
		\end{tikzpicture}
	\end{figure}	\FloatBarrier
	\begin{figure}[!h]
		\centering
		\begin{tikzpicture}
			\coordinate (a) at (0,0.5);
			\coordinate (b) at (4,0.5);
			\coordinate (c) at (0,2.5);
			\coordinate (d) at (4,2.5);
			\draw(a)--(b);
			\draw(b)--(d);
			\draw(d)--(c);
			\draw(c)--(a);
			\coordinate (x1) at (1,0.5);
			\coordinate (y1) at (1,2.5);
			\coordinate (x2) at (2,0.5);
			\coordinate (y2) at (2,2.5);
			\coordinate (x3) at (3,0.5);
			\coordinate (y3) at (3,2.5);
			\node(small-dot) at (1,2.8){$e_{1}$};
			\node(small-dot) at (2,2.8){$e_{2}$};
			\node(small-dot) at (3,2.8){$e_{3}$};
			\draw[dashed](x1)--(y1);
			\draw[dashed](x2)--(y2);
			\draw[dashed](x3)--(y3);
			\coordinate (p1) at (2,1);
			\coordinate (p2) at (3.5,1);
			\draw(p1)--(p2);
			\node(small-dot) at (2,0.2){$F_{x}$};
			\coordinate (a1) at (4.5,0.5);
			\coordinate (b1) at (8.5,0.5);
			\coordinate (c1) at (4.5,2.5);
			\coordinate (d1) at (8.5,2.5);
			\draw(a1)--(b1);
			\draw(b1)--(d1);
			\draw(d1)--(c1);
			\draw(c1)--(a1);
			\coordinate (x4) at (5.5,0.5);
			\coordinate (y4) at (5.5,2.5);
			\coordinate (x5) at (6.5,0.5);
			\coordinate (y5) at (6.5,2.5);
			\coordinate (x6) at (7.5,0.5);
			\coordinate (y6) at (7.5,2.5);
			\node(small-dot) at (5.5,2.8){$e_{1}$};
			\node(small-dot) at (6.5,2.8){$e_{2}$};
			\node(small-dot) at (7.5,2.8){$e_{3}$};
			\draw[dashed](x4)--(y4);
			\draw[dashed](x5)--(y5);
			\draw[dashed](x6)--(y6);
			\coordinate (pp1) at (5,2);
			\coordinate (pp2) at (5.5,2);
			\coordinate (pp3) at (5.5,1);
			\coordinate (pp4) at (6.5,1);
			\coordinate (pp5) at (7.5,1);
			\coordinate (pp6) at (8,1);
			\draw(pp1)--(pp2);
			\draw(pp2)--(pp3);
			\draw(pp3)--(pp4);
			\draw(pp5)--(pp6);
			\node(small-dot) at (6.5,0.2){$F_{y}$};
			\coordinate (a2) at (9,0.5);
			\coordinate (b2) at (13,0.5);
			\coordinate (c2) at (9,2.5);
			\coordinate (d2) at (13,2.5);
			\draw(a2)--(b2);
			\draw(b2)--(d2);
			\draw(d2)--(c2);
			\draw(c2)--(a2);
			\coordinate (x7) at (10,0.5);
			\coordinate (y7) at (10,2.5);
			\coordinate (x8) at (11,0.5);
			\coordinate (y8) at (11,2.5);
			\coordinate (x9) at (12,0.5);
			\coordinate (y9) at (12,2.5);
			\node(small-dot) at (10,2.8){$e_{1}$};
			\node(small-dot) at (11,2.8){$e_{2}$};
			\node(small-dot) at (12,2.8){$e_{3}$};
			\draw[dashed](x7)--(y7);
			\draw[dashed](x8)--(y8);
			\draw[dashed](x9)--(y9);
			\coordinate (ppp1) at (9.5,2);
			\coordinate (ppp2) at (12.5,2);
			\coordinate (ppp3) at (9.5,1);
			\coordinate (ppp4) at (12,1);
			\draw(ppp1)--(ppp2);
			\draw(ppp3)--(ppp4);
			\node(small-dot) at (11,0.2){$F_{z}$};
			\node(small-dot) at (14.5,1.5){$\bar{E}$};
			\node(small-dot) at (14.5,0.8){\tiny{$1\le e_{1}\le e_{2}< e_{3}\le n$}};
		\end{tikzpicture}
	\end{figure}	\FloatBarrier
	\begin{figure}[!h]
		\centering
		\begin{tikzpicture}
			\coordinate (a) at (0,0.5);
			\coordinate (b) at (4,0.5);
			\coordinate (c) at (0,2.5);
			\coordinate (d) at (4,2.5);
			\draw(a)--(b);
			\draw(b)--(d);
			\draw(d)--(c);
			\draw(c)--(a);
			\coordinate (x1) at (1,0.5);
			\coordinate (y1) at (1,2.5);
			\coordinate (x2) at (2,0.5);
			\coordinate (y2) at (2,2.5);
			\coordinate (x3) at (3,0.5);
			\coordinate (y3) at (3,2.5);
			\node(small-dot) at (1,2.8){$f_{1}$};
			\node(small-dot) at (2,2.8){$f_{2}$};
			\node(small-dot) at (3,2.8){$f_{3}$};
			\draw[dashed](x1)--(y1);
			\draw[dashed](x2)--(y2);
			\draw[dashed](x3)--(y3);
			\coordinate (p1) at (1,1);
			\coordinate (p2) at (3.5,1);
			\draw(p1)--(p2);
			\node(small-dot) at (2,0.2){$F_{x}$};
			\coordinate (a1) at (4.5,0.5);
			\coordinate (b1) at (8.5,0.5);
			\coordinate (c1) at (4.5,2.5);
			\coordinate (d1) at (8.5,2.5);
			\draw(a1)--(b1);
			\draw(b1)--(d1);
			\draw(d1)--(c1);
			\draw(c1)--(a1);
			\coordinate (x4) at (5.5,0.5);
			\coordinate (y4) at (5.5,2.5);
			\coordinate (x5) at (6.5,0.5);
			\coordinate (y5) at (6.5,2.5);
			\coordinate (x6) at (7.5,0.5);
			\coordinate (y6) at (7.5,2.5);
			\node(small-dot) at (5.5,2.8){$f_{1}$};
			\node(small-dot) at (6.5,2.8){$f_{2}$};
			\node(small-dot) at (7.5,2.8){$f_{3}$};
			\draw[dashed](x4)--(y4);
			\draw[dashed](x5)--(y5);
			\draw[dashed](x6)--(y6);
			\coordinate (pp1) at (5,1);
			\coordinate (pp2) at (5.5,1);
			\coordinate (pp3) at (6.5,1);
			\coordinate (pp4) at (8,1);
			\draw(pp1)--(pp2);
			\draw(pp3)--(pp4);
			\node(small-dot) at (6.5,0.2){$F_{y}$};
			\coordinate (a2) at (9,0.5);
			\coordinate (b2) at (13,0.5);
			\coordinate (c2) at (9,2.5);
			\coordinate (d2) at (13,2.5);
			\draw(a2)--(b2);
			\draw(b2)--(d2);
			\draw(d2)--(c2);
			\draw(c2)--(a2);
			\coordinate (x7) at (10,0.5);
			\coordinate (y7) at (10,2.5);
			\coordinate (x8) at (11,0.5);
			\coordinate (y8) at (11,2.5);
			\coordinate (x9) at (12,0.5);
			\coordinate (y9) at (12,2.5);
			\node(small-dot) at (10,2.8){$f_{1}$};
			\node(small-dot) at (11,2.8){$f_{2}$};
			\node(small-dot) at (12,2.8){$f_{3}$};
			\draw[dashed](x7)--(y7);
			\draw[dashed](x8)--(y8);
			\draw[dashed](x9)--(y9);
			\coordinate (ppp1) at (9.5,2);
			\coordinate (ppp2) at (12,2);
			\coordinate (ppp3) at (12,1);
			\coordinate (ppp4) at (12.5,1);
			\coordinate (ppp5) at (9.5,1);
			\coordinate (ppp6) at (11,1);
			\draw(ppp1)--(ppp2);
			\draw(ppp2)--(ppp3);
			\draw(ppp3)--(ppp4);
			\draw(ppp5)--(ppp6);
			\node(small-dot) at (11,0.2){$F_{z}$};
			\node(small-dot) at (14.5,1.5){$\bar{F}$};
			\node(small-dot) at (14.5,0.8){\tiny{$1\le f_{1}<f_{2}< f_{3}\le n$}};
		\end{tikzpicture}
	\end{figure}	\FloatBarrier
	\begin{proof}[Proof.]
		It is straightforward to verify that $F$ with diagram of the type $\bar{D}$, $\bar{E}$, or $\bar{F}$ is indeed a facet. Let $s_{1}$ be the minimal column index such that $x_{2,s_{1}}\in F_{x}$. In the diagrams $D,F^{*},G$, the point $y$ is not at the $n$th column. Since $F$ is maximal in $\Delta_{0}$, it follows that $y_{2,n}\in F$. Therefore, the second row of $F_{y}$ is nonempty. It is easy to see that $F_{y}$ must be a subset of one set corresponding to the following case:\\
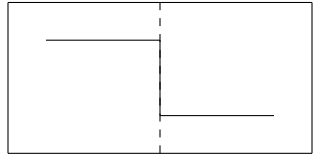
		\begin{figure}[h]
			\centering
			\begin{tikzpicture}
				\coordinate (a) at (0,0.5);
				\coordinate (b) at (4,0.5);
				\coordinate (c) at (0,2.5);
				\coordinate (d) at (4,2.5);
				\draw(a)--(b);
				\draw(b)--(d);
				\draw(d)--(c);
				\draw(c)--(a);
				\coordinate (p1) at (0.5,2);
				\coordinate (p2) at (2,2);
				\coordinate (p3) at (2,1);
				\coordinate (p4) at (3.5,1);
				\coordinate (a1) at (2,0.5);
				\coordinate (b1) at (2,2.5);
				\draw[dashed](a1)--(b1);
				\draw(p1)--(p2);
				\draw(p2)--(p3);
				\draw(p3)--(p4);
			\end{tikzpicture}
		\end{figure}    \FloatBarrier 
		\par There are two cases to consider. One denoted by case (i) is that the first row of $F_{y}$ is nonempty and the other denoted by case (ii) is that the first row of $F_{y}$ is empty.
		\par In case (i), let $s_{2}$ be the maximal column index such that $y_{1,s_{2}}\in F_{y}$ and $t$ be the maximal column index such that $z_{2,t}\in F_{z}$. In the diagrams $D,F^{*},G$, the point $y$ is not at the $n$th column. Then we have $y_{2,n}\in F$ because $F$ is maximal in $\Delta_{0}$. 
		\par If $s_{1}<s_{2}$, then we have $t\le s_{1}$ by restriction on D. Hence any set in $F^{*}, G$ cannot be contained in $F$. We have diagrams for this case as follows: \\
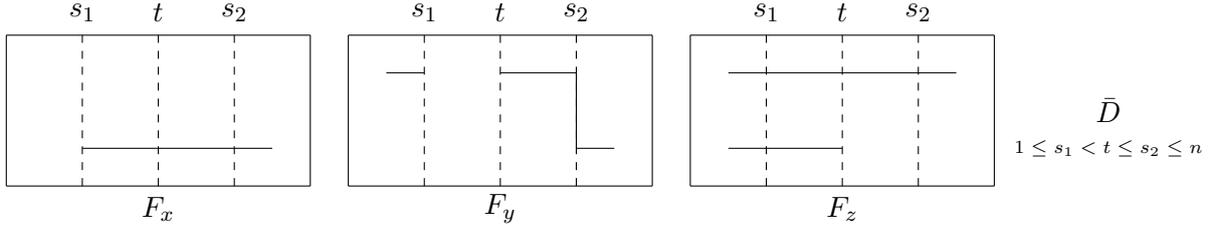
		\begin{figure}[h]
			\centering
			\begin{tikzpicture}
				\coordinate (a) at (0,0.5);
				\coordinate (b) at (4,0.5);
				\coordinate (c) at (0,2.5);
				\coordinate (d) at (4,2.5);
				\draw(a)--(b);
				\draw(b)--(d);
				\draw(d)--(c);
				\draw(c)--(a);
				\coordinate (x1) at (1,0.5);
				\coordinate (y1) at (1,2.5);
				\coordinate (x2) at (2,0.5);
				\coordinate (y2) at (2,2.5);
				\coordinate (x3) at (3,0.5);
				\coordinate (y3) at (3,2.5);
				\node(small-dot) at (1,2.8){$s_{1}$};
				\node(small-dot) at (2,2.8){$t$};
				\node(small-dot) at (3,2.8){$s_{2}$};
				\draw[dashed](x1)--(y1);
				\draw[dashed](x2)--(y2);
				\draw[dashed](x3)--(y3);
				\coordinate (p1) at (1,1);
				\coordinate (p2) at (3.5,1);
				\draw(p1)--(p2);
				\node(small-dot) at (2,0.2){$F_{x}$};
				\coordinate (a1) at (4.5,0.5);
				\coordinate (b1) at (8.5,0.5);
				\coordinate (c1) at (4.5,2.5);
				\coordinate (d1) at (8.5,2.5);
				\draw(a1)--(b1);
				\draw(b1)--(d1);
				\draw(d1)--(c1);
				\draw(c1)--(a1);
				\coordinate (x4) at (5.5,0.5);
				\coordinate (y4) at (5.5,2.5);
				\coordinate (x5) at (6.5,0.5);
				\coordinate (y5) at (6.5,2.5);
				\coordinate (x6) at (7.5,0.5);
				\coordinate (y6) at (7.5,2.5);
				\node(small-dot) at (5.5,2.8){$s_{1}$};
				\node(small-dot) at (6.5,2.8){$t$};
				\node(small-dot) at (7.5,2.8){$s_{2}$};
				\draw[dashed](x4)--(y4);
				\draw[dashed](x5)--(y5);
				\draw[dashed](x6)--(y6);
				\coordinate (pp1) at (5,2);
				\coordinate (pp2) at (5.5,2);
				\coordinate (pp3) at (6.5,2);
				\coordinate (pp4) at (7.5,2);
				\coordinate (pp5) at (7.5,1);
				\coordinate (pp6) at (8,1);
				\draw(pp1)--(pp2);
				\draw(pp3)--(pp4);
				\draw(pp4)--(pp5);
				\draw(pp5)--(pp6);
				\node(small-dot) at (6.5,0.2){$F_{y}$};
				\coordinate (a2) at (9,0.5);
				\coordinate (b2) at (13,0.5);
				\coordinate (c2) at (9,2.5);
				\coordinate (d2) at (13,2.5);
				\draw(a2)--(b2);
				\draw(b2)--(d2);
				\draw(d2)--(c2);
				\draw(c2)--(a2);
				\coordinate (x7) at (10,0.5);
				\coordinate (y7) at (10,2.5);
				\coordinate (x8) at (11,0.5);
				\coordinate (y8) at (11,2.5);
				\coordinate (x9) at (12,0.5);
				\coordinate (y9) at (12,2.5);
				\node(small-dot) at (10,2.8){$s_{1}$};
				\node(small-dot) at (11,2.8){$t$};
				\node(small-dot) at (12,2.8){$s_{2}$};
				\draw[dashed](x7)--(y7);
				\draw[dashed](x8)--(y8);
				\draw[dashed](x9)--(y9);
				\coordinate (ppp1) at (9.5,2);
				\coordinate (ppp2) at (12.5,2);
				\coordinate (ppp3) at (9.5,1);
				\coordinate (ppp4) at (11,1);
				\draw(ppp1)--(ppp2);
				\draw(ppp3)--(ppp4);
				\node(small-dot) at (11,0.2){$F_{z}$};
				\node(small-dot) at (14.5,1.5){$\bar{D}$};
				\node(small-dot) at (14.5,1){\tiny{$1\le s_{1}< t\le s_{2}\le n$}};
			\end{tikzpicture}
		\end{figure}	\FloatBarrier
		\par If $s_{1}\ge s_{2}$, then $F$ cannot contain any set in $D$.\\\\
		\textit{Claim }. The number $t$ is the minimal number such that $t>s_{1}$ and $y_{2,t}\in F_{y}$.
		\begin{proof}[Proof of Claim.]
			If $s_{1}=n$, then any set in $A,\cdots ,G$ cannot be contained in $F$. Suppose $F_{0}$ is a facet for $s_{1},s_{2},t$ in such condition. Let $F^{'}=F_{0}x_{1,n}$. Then $F^{'}\in \Delta_{0}$ and this yields a contradiction. Therefore $s_{1}<n$ and we can consider numbers larger than $s_{1}$. The existence of the minimal number is as $y_{2,n}\in F$. Suppose $t_{1}<t$ is the minimal number such that $t_{1}>s_{1}$ and $y_{2,t_{1}}\in F_{y}$. Then the subset $\{y_{1,s_{2}},x_{2,s_{1}},y_{2,t_{1}},z_{2,t} \}\subset F$ is of the form as diagram $F^{*}$, contradiction.
		\end{proof}
		\par We have diagrams for this case as follows: \\
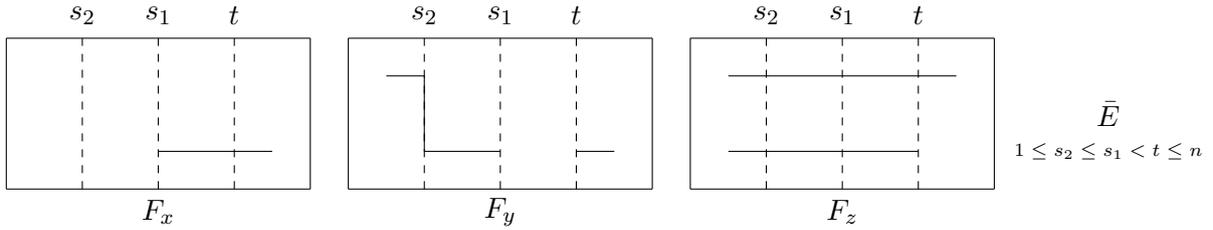
		\begin{figure}[!h]
			\centering
			\begin{tikzpicture}
				\coordinate (a) at (0,0.5);
				\coordinate (b) at (4,0.5);
				\coordinate (c) at (0,2.5);
				\coordinate (d) at (4,2.5);
				\draw(a)--(b);
				\draw(b)--(d);
				\draw(d)--(c);
				\draw(c)--(a);
				\coordinate (x1) at (1,0.5);
				\coordinate (y1) at (1,2.5);
				\coordinate (x2) at (2,0.5);
				\coordinate (y2) at (2,2.5);
				\coordinate (x3) at (3,0.5);
				\coordinate (y3) at (3,2.5);
				\node(small-dot) at (1,2.8){$s_{2}$};
				\node(small-dot) at (2,2.8){$s_{1}$};
				\node(small-dot) at (3,2.8){$t$};
				\draw[dashed](x1)--(y1);
				\draw[dashed](x2)--(y2);
				\draw[dashed](x3)--(y3);
				\coordinate (p1) at (2,1);
				\coordinate (p2) at (3.5,1);
				\draw(p1)--(p2);
				\node(small-dot) at (2,0.2){$F_{x}$};
				\coordinate (a1) at (4.5,0.5);
				\coordinate (b1) at (8.5,0.5);
				\coordinate (c1) at (4.5,2.5);
				\coordinate (d1) at (8.5,2.5);
				\draw(a1)--(b1);
				\draw(b1)--(d1);
				\draw(d1)--(c1);
				\draw(c1)--(a1);
				\coordinate (x4) at (5.5,0.5);
				\coordinate (y4) at (5.5,2.5);
				\coordinate (x5) at (6.5,0.5);
				\coordinate (y5) at (6.5,2.5);
				\coordinate (x6) at (7.5,0.5);
				\coordinate (y6) at (7.5,2.5);
				\node(small-dot) at (5.5,2.8){$s_{2}$};
				\node(small-dot) at (6.5,2.8){$s_{1}$};
				\node(small-dot) at (7.5,2.8){$t$};
				\draw[dashed](x4)--(y4);
				\draw[dashed](x5)--(y5);
				\draw[dashed](x6)--(y6);
				\coordinate (pp1) at (5,2);
				\coordinate (pp2) at (5.5,2);
				\coordinate (pp3) at (5.5,1);
				\coordinate (pp4) at (6.5,1);
				\coordinate (pp5) at (7.5,1);
				\coordinate (pp6) at (8,1);
				\draw(pp1)--(pp2);
				\draw(pp2)--(pp3);
				\draw(pp3)--(pp4);
				\draw(pp5)--(pp6);
				\node(small-dot) at (6.5,0.2){$F_{y}$};
				\coordinate (a2) at (9,0.5);
				\coordinate (b2) at (13,0.5);
				\coordinate (c2) at (9,2.5);
				\coordinate (d2) at (13,2.5);
				\draw(a2)--(b2);
				\draw(b2)--(d2);
				\draw(d2)--(c2);
				\draw(c2)--(a2);
				\coordinate (x7) at (10,0.5);
				\coordinate (y7) at (10,2.5);
				\coordinate (x8) at (11,0.5);
				\coordinate (y8) at (11,2.5);
				\coordinate (x9) at (12,0.5);
				\coordinate (y9) at (12,2.5);
				\node(small-dot) at (10,2.8){$s_{2}$};
				\node(small-dot) at (11,2.8){$s_{1}$};
				\node(small-dot) at (12,2.8){$t$};
				\draw[dashed](x7)--(y7);
				\draw[dashed](x8)--(y8);
				\draw[dashed](x9)--(y9);
				\coordinate (ppp1) at (9.5,2);
				\coordinate (ppp2) at (12.5,2);
				\coordinate (ppp3) at (9.5,1);
				\coordinate (ppp4) at (12,1);
				\draw(ppp1)--(ppp2);
				\draw(ppp3)--(ppp4);
				\node(small-dot) at (11,0.2){$F_{z}$};
				\node(small-dot) at (14.5,1.5){$\bar{E}$};
				\node(small-dot) at (14.5,1){\tiny{$1\le s_{2}\le s_{1}< t\le n$}};
			\end{tikzpicture}
		\end{figure}	\FloatBarrier
		\par In case (ii), any set in $D$ or $F^{*}$ cannot be contained in $F$. Let $s$ be the minimal column index such that $y_{2,s}\in F_{y}$. In the diagram $G$, the point $y$ is not at the first column and the point $x$ is not at the $(n-2)$th, $(n-1)$th, and $n$th columns. Then we have $y_{2,1}\in F$ and $s_{1}\le (n-2)$ because $F$ is maximal in $\Delta_{0}$. Let $s_{2}$ be the minimal column index such that $s_{1}<s_{2}$ and $y_{2,s_{2}}\in F_{y}$. By restriction of $G$, the $\tiny{>}s_{2}$-part of $F_{z}$ must be a path walking from $z_{1,s_{2}+1}$ to $z_{2,n}$ only to the right or down. Let $t$ be the maximal column index such that $z_{1,t}\in F_{z}$. Then we have diagrams for this case as follows:\\
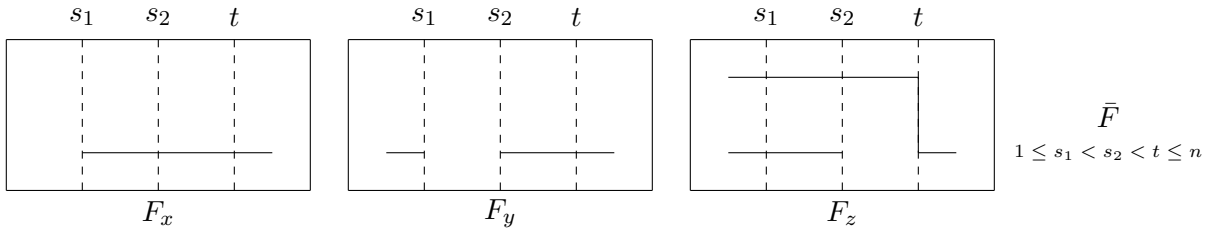
		\begin{figure}[h]
			\centering
			\begin{tikzpicture}
				\coordinate (a) at (0,0.5);
				\coordinate (b) at (4,0.5);
				\coordinate (c) at (0,2.5);
				\coordinate (d) at (4,2.5);
				\draw(a)--(b);
				\draw(b)--(d);
				\draw(d)--(c);
				\draw(c)--(a);
				\coordinate (x1) at (1,0.5);
				\coordinate (y1) at (1,2.5);
				\coordinate (x2) at (2,0.5);
				\coordinate (y2) at (2,2.5);
				\coordinate (x3) at (3,0.5);
				\coordinate (y3) at (3,2.5);
				\node(small-dot) at (1,2.8){$s_{1}$};
				\node(small-dot) at (2,2.8){$s_{2}$};
				\node(small-dot) at (3,2.8){$t$};
				\draw[dashed](x1)--(y1);
				\draw[dashed](x2)--(y2);
				\draw[dashed](x3)--(y3);
				\coordinate (p1) at (1,1);
				\coordinate (p2) at (3.5,1);
				\draw(p1)--(p2);
				\node(small-dot) at (2,0.2){$F_{x}$};
				\coordinate (a1) at (4.5,0.5);
				\coordinate (b1) at (8.5,0.5);
				\coordinate (c1) at (4.5,2.5);
				\coordinate (d1) at (8.5,2.5);
				\draw(a1)--(b1);
				\draw(b1)--(d1);
				\draw(d1)--(c1);
				\draw(c1)--(a1);
				\coordinate (x4) at (5.5,0.5);
				\coordinate (y4) at (5.5,2.5);
				\coordinate (x5) at (6.5,0.5);
				\coordinate (y5) at (6.5,2.5);
				\coordinate (x6) at (7.5,0.5);
				\coordinate (y6) at (7.5,2.5);
				\node(small-dot) at (5.5,2.8){$s_{1}$};
				\node(small-dot) at (6.5,2.8){$s_{2}$};
				\node(small-dot) at (7.5,2.8){$t$};
				\draw[dashed](x4)--(y4);
				\draw[dashed](x5)--(y5);
				\draw[dashed](x6)--(y6);
				\coordinate (pp1) at (5,1);
				\coordinate (pp2) at (5.5,1);
				\coordinate (pp3) at (6.5,1);
				\coordinate (pp4) at (8,1);
				\draw(pp1)--(pp2);
				\draw(pp3)--(pp4);
				\node(small-dot) at (6.5,0.2){$F_{y}$};
				\coordinate (a2) at (9,0.5);
				\coordinate (b2) at (13,0.5);
				\coordinate (c2) at (9,2.5);
				\coordinate (d2) at (13,2.5);
				\draw(a2)--(b2);
				\draw(b2)--(d2);
				\draw(d2)--(c2);
				\draw(c2)--(a2);
				\coordinate (x7) at (10,0.5);
				\coordinate (y7) at (10,2.5);
				\coordinate (x8) at (11,0.5);
				\coordinate (y8) at (11,2.5);
				\coordinate (x9) at (12,0.5);
				\coordinate (y9) at (12,2.5);
				\node(small-dot) at (10,2.8){$s_{1}$};
				\node(small-dot) at (11,2.8){$s_{2}$};
				\node(small-dot) at (12,2.8){$t$};
				\draw[dashed](x7)--(y7);
				\draw[dashed](x8)--(y8);
				\draw[dashed](x9)--(y9);
				\coordinate (ppp1) at (9.5,2);
				\coordinate (ppp2) at (12,2);
				\coordinate (ppp3) at (12,1);
				\coordinate (ppp4) at (12.5,1);
				\coordinate (ppp5) at (9.5,1);
				\coordinate (ppp6) at (11,1);
				\draw(ppp1)--(ppp2);
				\draw(ppp2)--(ppp3);
				\draw(ppp3)--(ppp4);
				\draw(ppp5)--(ppp6);
				\node(small-dot) at (11,0.2){$F_{z}$};
				\node(small-dot) at (14.5,1.5){$\bar{F}$};
				\node(small-dot) at (14.5,1){\tiny{$1\le s_{1}<s_{2}< t\le n$}};
			\end{tikzpicture}
		\end{figure}	\FloatBarrier
	\end{proof}
	\par In conclusion, we have shown that all facets of $\Delta_{0}$ belong to one of the five families: $\bar{A},\bar{C},\bar{D},\bar{E},\bar{F}$. Furthermore, each facet contains exactly $3n+3$ points, so $dim\Delta_{0}=3n+2$.\\\\
	\textbf{Corollary 4.4. }The Stanley-Reisner complex $\Delta_{0}$ is pure of dimension $3n+2$. \\

	\section{Shellability of $\Delta_{0}$}
	\par In this section, we prove that $\Delta_{0}$ is shellable by explicitly constructing a shelling on its facets. 
	\par Let $V=\{ x_{1,1},\cdots x_{2,n}, y_{1,1},\cdots ,y_{2,n},z_{1,1},\cdots ,z_{2,n}\}$ be the vertex set of $\Delta_{0}$. For a facet $F\in \Delta_{0}$, we define a function $\delta _{F}:V\to \{ 0,1\}$ as follows:
	\begin{center}
		$\delta _{F}(v)=\left\{
		\begin{array}{rcl}
			1,\quad v\in F \\
			0,\quad v\notin F
		\end{array}
		\right.$
	\end{center}	 
	\par To define a shelling (just a total ordering of facets at first) on $\Delta_{0}$, we give a new ordering for $V$ as follows:
	\begin{center}
		$x_{1,1}>\cdots >x_{1,n}>x_{2,1}>\cdots >x_{2,n}>y_{1,1}>\cdots >y_{2,n}>z_{1,1}>\cdots >z_{2,n}$.
	\end{center}
	\par We describe the ordering for facets of $\Delta_{0}$ now. The ordering of families of facets is as $\bar{F}<\bar{E}<\bar{D}<\bar{C}<\bar{A}$. For any two facets in two different families, one in the larger family is larger than the other. Let $P$ and $Q$ be two different facets in the same family. We have $P>Q$ if and only if the largest vertice $v\in V$ such that $\delta _{P}(v)\neq \delta _{Q}(v)$ satisfies $\delta _{P}(p)=1$.
	\par For convenience, we call it the $*$-ordering.\\ \\
	\textbf{Example 5.1. }Paths $P$ and $Q$ in the following three examples have the same diagram in $y$-part and $z$-part. We only use diagrams of $P_{x}$ and $Q_{x}$ to show the ordering as defined above.
	
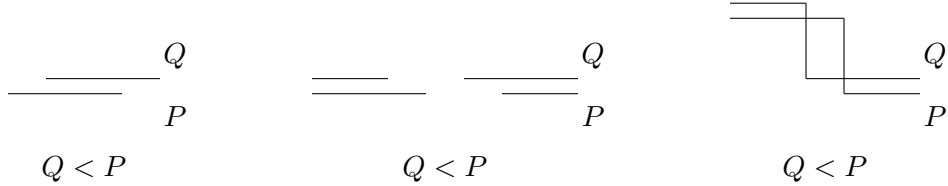
\begin{figure}[h]
		\centering
		\begin{tikzpicture}
			\coordinate (a) at (0.5,0.2);
			\coordinate (b) at (2,0.2);
			\coordinate (c) at (0,0);
			\coordinate (d) at (1.5,0);
			\draw(a)--(b);
			\draw(d)--(c);
			\node(small-dot) at (2.2,0.5){$Q$};
			\node(small-dot) at (2.2,-0.3){$P$};
			\node(small-dot) at (1,-1){$Q<P$};
			\coordinate (aa1) at (4,0.2);
			\coordinate (bb1) at (5,0.2);
			\coordinate (aa2) at (6,0.2);
			\coordinate (bb2) at (7.5,0.2);
			\coordinate (cc1) at (4,0);
			\coordinate (dd1) at (5.5,0);
			\coordinate (cc2) at (6.5,0);
			\coordinate (dd2) at (7.5,0);
			\draw(aa1)--(bb1);
			\draw(aa2)--(bb2);
			\draw(dd1)--(cc1);
			\draw(dd2)--(cc2);
			\node(small-dot) at (7.7,0.5){$Q$};
			\node(small-dot) at (7.7,-0.3){$P$};
			\node(small-dot) at (5.75,-1){$Q<P$};
			\coordinate (a11) at (9.5,1.2);
			\coordinate (b11) at (10.5,1.2);
			\coordinate (a12) at (10.5,0.2);
			\coordinate (b12) at (12,0.2);
			\coordinate (c11) at (9.5,1);
			\coordinate (d11) at (11,1);
			\coordinate (c12) at (11,0);
			\coordinate (d12) at (12,0);
			\draw(a11)--(b11);
			\draw(b11)--(a12);
			\draw(a12)--(b12);
			\draw(d12)--(c12);
			\draw(c12)--(d11);
			\draw(d11)--(c11);
			\node(small-dot) at (12.2,0.5){$Q$};
			\node(small-dot) at (12.2,-0.3){$P$};
			\node(small-dot) at (10.75,-1){$Q<P$};
		\end{tikzpicture}
	\end{figure}    \FloatBarrier
	\par Unless otherwise stated, we will assume $P>Q$ to be two different facets of $\Delta_{0}$. If $P-Q$ contains a unique point $v$, we have $v\in c(P)$ (recall the notation in Definition 2.12). We define two sets $c_{1}(P),c_{2}(P)\subset c(P)$. A vertice $v\in c_{1}(P)$ if and only if there exists $P^{'}<P$ in the same family as $P$ such that $P-P^{'}=\{v\}$. A vertice $v\in c_{2}(P)$ if and only if there exists $P^{''}<P$ in some different family from $P$ such that $P-P^{''}=\{v\}$. Hence we have $c(P)=c_{1}(P)\cup c_{2}(P)$. We say $Q$ covers $c(P)$ if $c(P)\subset Q$.
	\par For a facet $F$ of $\Delta_{0}$, its $F_{x}$-part, $F_{y}$-part, and $F_{z}$-part respectively has three dashed column lines. These dashed lines are denoted by dashed axes of the facet $F$. We give some special names for some vertices on the dashed axes in the diagram of each facet as follows:
	\begin{figure}[!h]
		\centering
		\begin{tikzpicture}
			\coordinate (a1) at (0,5);
			\coordinate (b1) at (13,5);
			\coordinate (a2) at (0,4);
			\coordinate (b2) at (13,4);
			\coordinate (a3) at (0,3);
			\coordinate (b3) at (13,3);
			\coordinate (a4) at (0,2);
			\coordinate (b4) at (13,2);
			\coordinate (a5) at (0,1);
			\coordinate (b5) at (13,1);
			\coordinate (a6) at (0,0);
			\coordinate (b6) at (13,0);
			\coordinate (a7) at (0,6);
			\coordinate (b7) at (13,6);
			\draw(a1)--(b1);
			\draw(a2)--(b2);
			\draw(a3)--(b3);
			\draw(a4)--(b4);
			\draw(a5)--(b5);
			\draw(a6)--(b6);
			\draw(a7)--(b7);
			\coordinate (c1) at (1,6);
			\coordinate (c2) at (4,6);
			\coordinate (c3) at (7,6);
			\coordinate (c4) at (10,6);
			\coordinate (d1) at (1,0);
			\coordinate (d2) at (4,0);
			\coordinate (d3) at (7,0);
			\coordinate (d4) at (10,0);
			\draw(a7)--(a6);
			\draw(b7)--(b6);
			\draw(c3)--(d3);
			\draw(c4)--(d4);
			\draw(c2)--(d2);
			\draw(c1)--(d1);
			\node(small-dot) at (0.5,4.5){$\bar{A}$};
			\node(small-dot) at (0.5,3.5){$\bar{C}$};
			\node(small-dot) at (0.5,2.5){$\bar{D}$};
			\node(small-dot) at (0.5,1.5){$\bar{E}$};
			\node(small-dot) at (0.5,0.5){$\bar{F}$};
			\node(small-dot) at (2.5,4.5){$x_{1,a_{r}}$};
			\node(small-dot) at (2.5,3.5){$x_{1,c_{3}},y_{1,c_{1}}$};
			\node(small-dot) at (2.5,2.5){$y_{1,d_{3}}$};
			\node(small-dot) at (2.5,1.5){$y_{1,e_{1}}$};
			\node(small-dot) at (2.5,0.5){$z_{1,f_{3}}$};
			\node(small-dot) at (5.5,4.5){$x_{2,a_{r}}$};
			\node(small-dot) at (5.5,3.5){$x_{2,c_{3}},y_{1,c_{1}}$};
			\node(small-dot) at (5.5,2.5){$y_{2,d_{3}}$};
			\node(small-dot) at (5.5,1.5){$y_{2,e_{1}}$};
			\node(small-dot) at (5.5,0.5){$z_{2,f_{3}}$};
			\node(small-dot) at (8.5,4.5){$y_{1,a_{2}}$};
			\node(small-dot) at (8.5,3.5){$x_{1,c_{2}}$};
			\node(small-dot) at (8.5,2.5){$y_{1,d_{2}}$};
			\node(small-dot) at (8.5,1.5){$y_{2,e_{3}}$};
			\node(small-dot) at (8.5,0.5){$x_{2,f_{1}},y_{2,f_{2}}$};
			\node(small-dot) at (11.5,4.5){$y_{1,a_{1}}$};
			\node(small-dot) at (11.5,3.5){$y_{2,c_{2}}$};
			\node(small-dot) at (11.5,2.5){$y_{1,d_{1}}$};
			\node(small-dot) at (11.5,1.5){$y_{2,e_{2}}$};
			\node(small-dot) at (11.5,0.5){$y_{2,f_{1}},z_{2,f_{2}}$};
			\node(small-dot) at (2.5,5.5){NE-points};
			\node(small-dot) at (5.5,5.5){SW-points};
			\node(small-dot) at (8.5,5.5){L-points};
			\node(small-dot) at (11.5,5.5){R-points};
		\end{tikzpicture}
	\end{figure}	\FloatBarrier
	\par \noindent \textbf{Example 5.2. }The following diagram is $F_{y}$-part of diagram $\bar{E}$.
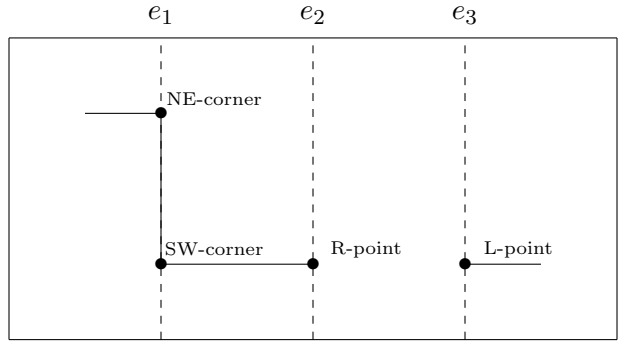
	\begin{figure}[h]
		\centering
		\begin{tikzpicture}
			\coordinate (a) at (0,0.5);
			\coordinate (b) at (8,0.5);
			\coordinate (c) at (0,4.5);
			\coordinate (d) at (8,4.5);
			\draw(a)--(b);
			\draw(b)--(d);
			\draw(d)--(c);
			\draw(c)--(a);
			\coordinate (x1) at (2,0.5);
			\coordinate (y1) at (2,4.5);
			\coordinate (x2) at (4,0.5);
			\coordinate (y2) at (4,4.5);
			\coordinate (x3) at (6,0.5);
			\coordinate (y3) at (6,4.5);
			\draw[dashed](x1)--(y1);
			\draw[dashed](x2)--(y2);
			\draw[dashed](x3)--(y3);
			\coordinate (p1) at (1,3.5);
			\coordinate (p2) at (2,3.5);
			\coordinate (p3) at (2,1.5);
			\coordinate (p4) at (4,1.5);
			\coordinate (p5) at (6,1.5);
			\coordinate (p6) at (7,1.5);
			\draw(p1)--(p2);
			\draw(p2)--(p3);
			\draw(p3)--(p4);
			\draw(p5)--(p6);
			\node(small-dot) at (2,3.5){$\bullet$};
			\node(small-dot) at (2,1.5){$\bullet$};
			\node(small-dot) at (4,1.5){$\bullet$};
			\node(small-dot) at (6,1.5){$\bullet$};
			\node(small-dot) at (2.7,3.7){\tiny{NE-corner}};
			\node(small-dot) at (2.7,1.7){\tiny{SW-corner}};
			\node(small-dot) at (4.7,1.7){\tiny{R-point}};
			\node(small-dot) at (6.7,1.7){\tiny{L-point}};
			\node(small-dot) at (2,4.8){$e_{1}$};
			\node(small-dot) at (4,4.8){$e_{2}$};
			\node(small-dot) at (6,4.8){$e_{3}$};
		\end{tikzpicture}
	\end{figure}	\FloatBarrier
	\par \noindent \textbf{Proposition 5.3. }\textit{If $v\in c_{1}(P)$, then $v$ is either an NE-corner or an R-point or an L-point in the diagram of $P$.}
	\begin{proof}[Proof.]
		Let $P$ and $Q$ be facets and $P>Q$. Suppose $P-Q=\{v\}$ with $v=v_{i_{1},j_{1}}$ and $Q-P=\{w\}$. For convenience, we use integers $p_{i}$ and $q_{i}$ to represent places of dashed axes of $P$ and $Q$ from left to right. We suppose $p_{0}=q_{0}=1$ and $p_{4}=q_{4}=n$
		\par \textit{Step }1. $v$ is on some dashed axes of $P$.
		\par The following discussion is in one of the three parts, $P_{x}$, $P_{y}$, or $P_{z}$, such that $v$ is contained. If $v$ is not on any dashed axis, then there exists $i$ such that $p_{i}<j_{1}<p_{i+1}$ with path $P$ walking from $p_{i}$ to $p_{j}$ at the $i_{1}$th row and path $Q$ empty at $p_{i}<\cdots <p_{i+1}$-part of the $i_{1}$th row. The set $P-Q$ has one element, so $j_{1}=p_{i}+1$ and $p_{i+1}=p_{i}+2$. The $i_{1}$th row of $Q$ contains points at $p_{i}$th column and $(p_{i}+2)$th column. Only $F_{y}$-part of diagrams $\bar{A},\bar{D},\bar{E},\bar{F}$ and $F_{z}$-part of diagram $\bar{F}$ admit such condition. We have $P>Q$ with $*$-ordering, so $q_{i+1}<p_{i+1}$ and $P_{x}$ is the same as $Q_{x}$. 
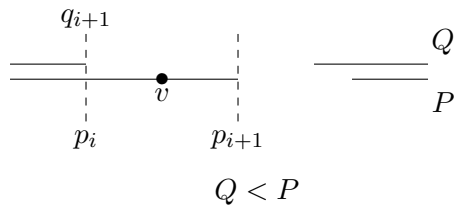
		\begin{figure}[h]
			\centering
			\begin{tikzpicture}
				\coordinate (aa1) at (2,0.2);
				\coordinate (bb1) at (3,0.2);
				\coordinate (aa2) at (6,0.2);
				\coordinate (bb2) at (7.5,0.2);
				\coordinate (cc1) at (2,0);
				\coordinate (dd1) at (5,0);
				\coordinate (cc2) at (6.5,0);
				\coordinate (dd2) at (7.5,0);
				\draw(aa1)--(bb1);
				\draw(aa2)--(bb2);
				\draw(dd1)--(cc1);
				\draw(dd2)--(cc2);
				\node(small-dot) at (7.7,0.5){$Q$};
				\node(small-dot) at (7.7,-0.3){$P$};
				\node(small-dot) at (5.25,-1.5){$Q<P$};
				\node(small-dot) at (3,-0.8){$\tiny{p_{i}}$};
				\node(small-dot) at (5,-0.8){$\tiny{p_{i+1}}$};
				\node(small-dot) at (4,-0.2){$v$};
				\node(small-dot) at (3,0.8){$\tiny{q_{i+1}}$};
				\node(small-dot) at (4,0){$\bullet$};
				\coordinate (a) at (3,0.6);
				\coordinate (b) at (3,-0.6);
				\coordinate (c) at (5,0.6);
				\coordinate (d) at (5,-0.6);
				\draw[dashed](a)--(b);
				\draw[dashed](c)--(d);
			\end{tikzpicture}
		\end{figure}    \FloatBarrier
		\par The four cases of $F_{y}$ as mentioned above are impossible for which $P_{x}$ is not the same as $Q_{x}$. Therefore it suffices to show that the $F_{z}$-part of diagram $\bar{F}$ is also impossible. The $*$-ordering shows that $P_{y}$ is the same as $Q_{y}$ if $v\in P_{z}$. Hence $P_{y}$ and $Q_{y}$ are both the whole second row, otherwise $p_{2}>q_{2}$ induces contradiction. Then we have $p_{1}=p_{2}-1$ and $q_{1}=q_{2}-1$, so $p_{1}>q_{1}$ induces contradiction that $P_{x}$ and $Q_{x}$ are not the same.
		\par \textit{Step }2. SW-corner of $P$ is not in $c_{1}(P)$.
		\par If $v_{i_{1},2}$ is an SW-corner, then the following local diagram of $P$ and $Q$ is the unique condition to make $P-Q=\{v \}$:
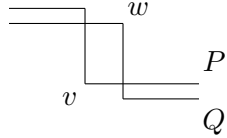
		\begin{figure}[h]
			\centering
			\begin{tikzpicture}
				\coordinate (a11) at (0,1.2);
				\coordinate (b11) at (1,1.2);
				\coordinate (a12) at (1,0.2);
				\coordinate (b12) at (2.5,0.2);
				\coordinate (c11) at (0,1);
				\coordinate (d11) at (1.5,1);
				\coordinate (c12) at (1.5,0);
				\coordinate (d12) at (2.5,0);
				\draw(a11)--(b11);
				\draw(b11)--(a12);
				\draw(a12)--(b12);
				\draw(d12)--(c12);
				\draw(c12)--(d11);
				\draw(d11)--(c11);
				\node(small-dot) at (2.7,0.5){$P$};
				\node(small-dot) at (2.7,-0.3){$Q$};
				\node(small-dot) at (1.7,1.2){$w$};
				\node(small-dot) at (0.8,0){$v$};
			\end{tikzpicture}
		\end{figure}    \FloatBarrier
		\par In this condition, $v$ and $w$ are all fixed so the other parts of $P$ and $Q$ are the same. Then we have $P<Q$, contradiction.
	\end{proof}	
	\par We depict $c_{1}(P)$ for every facet $P\in \Delta_{0}$ following. Proposition 5.3 shows it enough to consider NE-corners, R-points, and L-points of $P$. If $v\in c_{1}(P)$ is generated by $P-Q$ in the same family, then the local diagram is as follows:
	
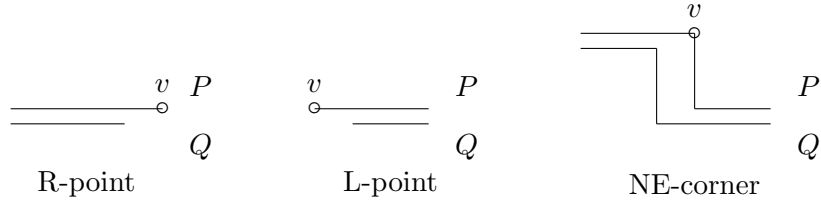
\begin{figure}[h]
		\centering
		\begin{tikzpicture}
			\coordinate (a) at (0,0.2);
			\coordinate (b) at (2,0.2);
			\coordinate (c) at (0,0);
			\coordinate (d) at (1.5,0);
			\draw(a)--(b);
			\draw(d)--(c);
			\node(small-dot) at (2.5,0.5){$P$};
			\node(small-dot) at (2.5,-0.3){$Q$};
			\node(small-dot) at (2,0.2){$\circ$};
			\coordinate (aa1) at (4,0.2);
			\coordinate (bb1) at (5.5,0.2);
			\coordinate (cc1) at (4.5,0);
			\coordinate (dd1) at (5.5,0);
			\draw(aa1)--(bb1);
			\draw(dd1)--(cc1);
			\node(small-dot) at (6,0.5){$P$};
			\node(small-dot) at (6,-0.3){$Q$};
			\node(small-dot) at (4,0.2){$\circ$};
			\coordinate (a11) at (7.5,1.2);
			\coordinate (b11) at (9,1.2);
			\coordinate (a12) at (9,0.2);
			\coordinate (b12) at (10,0.2);
			\coordinate (c11) at (7.5,1);
			\coordinate (d11) at (8.5,1);
			\coordinate (c12) at (8.5,0);
			\coordinate (d12) at (10,0);
			\draw(a11)--(b11);
			\draw(b11)--(a12);
			\draw(a12)--(b12);
			\draw(d12)--(c12);
			\draw(c12)--(d11);
			\draw(d11)--(c11);
			\node(small-dot) at (10.5,0.5){$P$};
			\node(small-dot) at (10.5,-0.3){$Q$};
			\node(small-dot) at (9,1.2){$\circ$};
			\node(small-dot) at (2,0.5){$v$};
			\node(small-dot) at (4,0.5){$v$};
			\node(small-dot) at (9,1.5){$v$};
			\node(small-dot) at (1,-0.8){R-point};
			\node(small-dot) at (5,-0.8){L-point};
			\node(small-dot) at (9,-0.8){NE-corner};
		\end{tikzpicture}
	\end{figure}    \FloatBarrier
	\par We consider generic condition, that is, any two adjacent dashed axes are neither coincident nor in the adjacent columns and the 1st-column and nth-column are not axis. We label points in $c_{1}(P)$ from left to right and draw a hollow circle to represent them, details as follows:
	
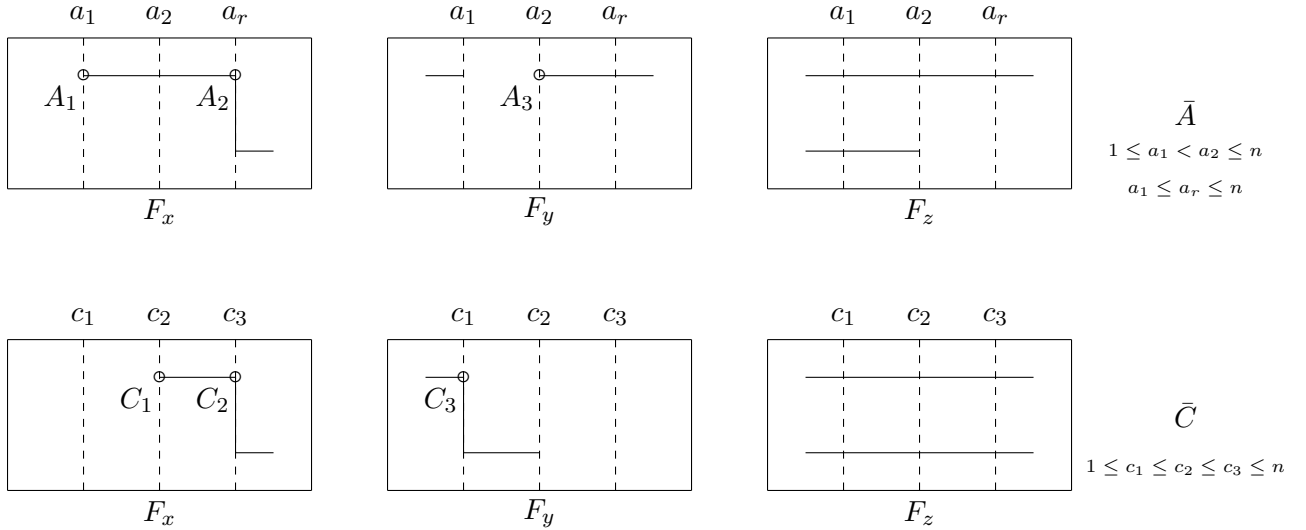
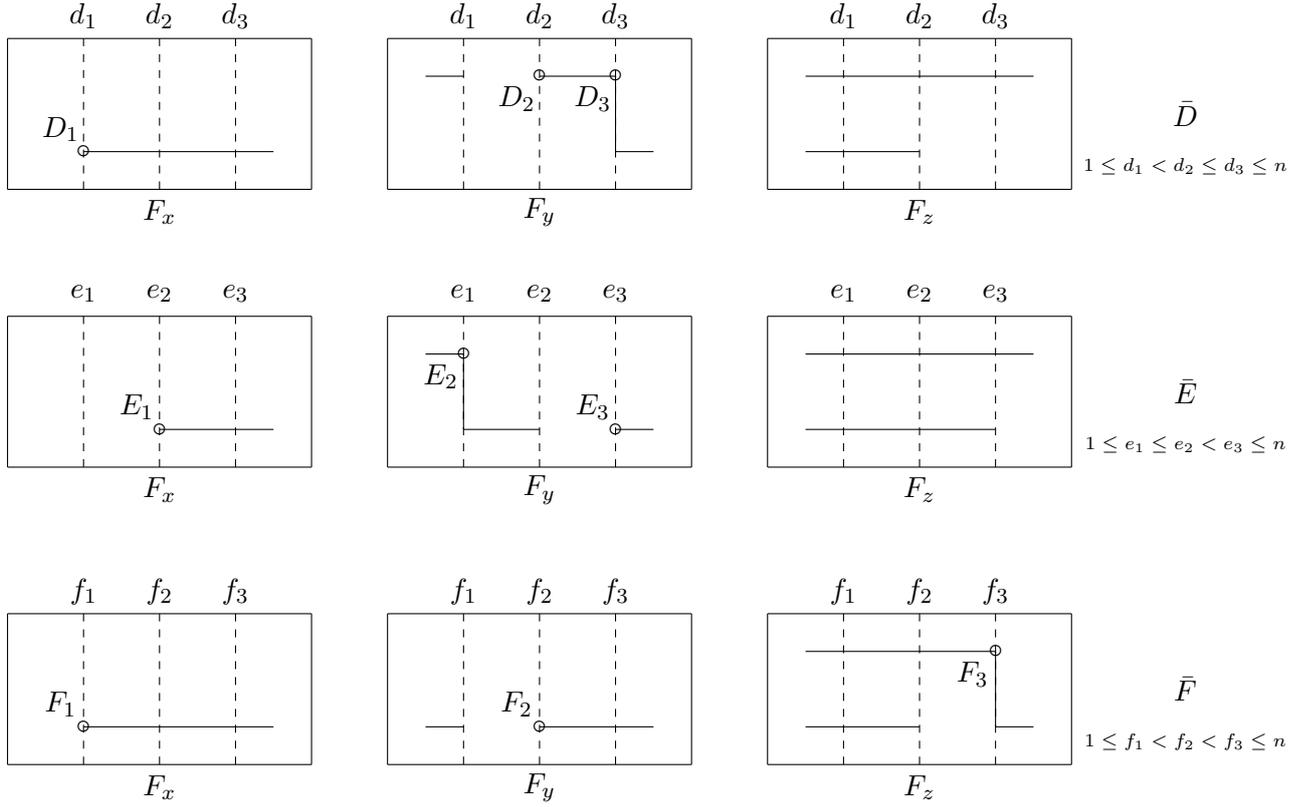
\begin{figure}[!h]
		\centering
		\begin{tikzpicture}
			\coordinate (a) at (0,0.5);
			\coordinate (b) at (4,0.5);
			\coordinate (c) at (0,2.5);
			\coordinate (d) at (4,2.5);
			\draw(a)--(b);
			\draw(b)--(d);
			\draw(d)--(c);
			\draw(c)--(a);
			\coordinate (x1) at (1,0.5);
			\coordinate (y1) at (1,2.5);
			\coordinate (x2) at (2,0.5);
			\coordinate (y2) at (2,2.5);
			\coordinate (x3) at (3,0.5);
			\coordinate (y3) at (3,2.5);
			\node(small-dot) at (1,2.8){$a_{1}$};
			\node(small-dot) at (2,2.8){$a_{2}$};
			\node(small-dot) at (3,2.8){$a_{r}$};
			\draw[dashed](x1)--(y1);
			\draw[dashed](x2)--(y2);
			\draw[dashed](x3)--(y3);
			\coordinate (p1) at (1,2);
			\coordinate (p2) at (3,2);
			\coordinate (p3) at (3,1);
			\coordinate (p4) at (3.5,1);
			\draw(p1)--(p2);
			\draw(p2)--(p3);
			\draw(p3)--(p4);
			\node(small-dot) at (2,0.2){$F_{x}$};
			\coordinate (a1) at (5,0.5);
			\coordinate (b1) at (9,0.5);
			\coordinate (c1) at (5,2.5);
			\coordinate (d1) at (9,2.5);
			\draw(a1)--(b1);
			\draw(b1)--(d1);
			\draw(d1)--(c1);
			\draw(c1)--(a1);
			\coordinate (x4) at (6,0.5);
			\coordinate (y4) at (6,2.5);
			\coordinate (x5) at (7,0.5);
			\coordinate (y5) at (7,2.5);
			\coordinate (x6) at (8,0.5);
			\coordinate (y6) at (8,2.5);
			\node(small-dot) at (6,2.8){$a_{1}$};
			\node(small-dot) at (7,2.8){$a_{2}$};
			\node(small-dot) at (8,2.8){$a_{r}$};
			\draw[dashed](x4)--(y4);
			\draw[dashed](x5)--(y5);
			\draw[dashed](x6)--(y6);
			\coordinate (pp1) at (5.5,2);
			\coordinate (pp2) at (6,2);
			\coordinate (pp3) at (7,2);
			\coordinate (pp4) at (8.5,2);
			\draw(pp1)--(pp2);
			\draw(pp3)--(pp4);
			\node(small-dot) at (7,0.2){$F_{y}$};
			\coordinate (a2) at (10,0.5);
			\coordinate (b2) at (14,0.5);
			\coordinate (c2) at (10,2.5);
			\coordinate (d2) at (14,2.5);
			\draw(a2)--(b2);
			\draw(b2)--(d2);
			\draw(d2)--(c2);
			\draw(c2)--(a2);
			\coordinate (x7) at (11,0.5);
			\coordinate (y7) at (11,2.5);
			\coordinate (x8) at (12,0.5);
			\coordinate (y8) at (12,2.5);
			\coordinate (x9) at (13,0.5);
			\coordinate (y9) at (13,2.5);
			\node(small-dot) at (11,2.8){$a_{1}$};
			\node(small-dot) at (12,2.8){$a_{2}$};
			\node(small-dot) at (13,2.8){$a_{r}$};
			\draw[dashed](x7)--(y7);
			\draw[dashed](x8)--(y8);
			\draw[dashed](x9)--(y9);
			\coordinate (ppp1) at (10.5,2);
			\coordinate (ppp2) at (13.5,2);
			\coordinate (ppp3) at (10.5,1);
			\coordinate (ppp4) at (12,1);
			\draw(ppp1)--(ppp2);
			\draw(ppp3)--(ppp4);
			\node(small-dot) at (12,0.2){$F_{z}$};
			\node(small-dot) at (15.5,1.5){$\bar{A}$};
			\node(small-dot) at (1,2){$\circ$};
			\node(small-dot) at (0.7,1.7){$A_{1}$};
			\node(small-dot) at (3,2){$\circ$};
			\node(small-dot) at (2.7,1.7){$A_{2}$};
			\node(small-dot) at (7,2){$\circ$};
			\node(small-dot) at (6.7,1.7){$A_{3}$};
			\node(small-dot) at (15.5,1){\tiny{$1\le a_{1}< a_{2}\le n$}};
			\node(small-dot) at (15.5,0.5){\tiny{$a_{1}\le a_{r}\le n$}};
		\end{tikzpicture}
	\end{figure}	\FloatBarrier
	\begin{figure}[!h]
		\centering
		\begin{tikzpicture}
			\coordinate (a) at (0,0.5);
			\coordinate (b) at (4,0.5);
			\coordinate (c) at (0,2.5);
			\coordinate (d) at (4,2.5);
			\draw(a)--(b);
			\draw(b)--(d);
			\draw(d)--(c);
			\draw(c)--(a);
			\coordinate (x1) at (1,0.5);
			\coordinate (y1) at (1,2.5);
			\coordinate (x2) at (2,0.5);
			\coordinate (y2) at (2,2.5);
			\coordinate (x3) at (3,0.5);
			\coordinate (y3) at (3,2.5);
			\node(small-dot) at (1,2.8){$c_{1}$};
			\node(small-dot) at (2,2.8){$c_{2}$};
			\node(small-dot) at (3,2.8){$c_{3}$};
			\draw[dashed](x1)--(y1);
			\draw[dashed](x2)--(y2);
			\draw[dashed](x3)--(y3);
			\coordinate (p1) at (2,2);
			\coordinate (p2) at (3,2);
			\coordinate (p3) at (3,1);
			\coordinate (p4) at (3.5,1);
			\draw(p1)--(p2);
			\draw(p2)--(p3);
			\draw(p3)--(p4);
			\node(small-dot) at (2,0.2){$F_{x}$};
			\coordinate (a1) at (5,0.5);
			\coordinate (b1) at (9,0.5);
			\coordinate (c1) at (5,2.5);
			\coordinate (d1) at (9,2.5);
			\draw(a1)--(b1);
			\draw(b1)--(d1);
			\draw(d1)--(c1);
			\draw(c1)--(a1);
			\coordinate (x4) at (6,0.5);
			\coordinate (y4) at (6,2.5);
			\coordinate (x5) at (7,0.5);
			\coordinate (y5) at (7,2.5);
			\coordinate (x6) at (8,0.5);
			\coordinate (y6) at (8,2.5);
			\node(small-dot) at (6,2.8){$c_{1}$};
			\node(small-dot) at (7,2.8){$c_{2}$};
			\node(small-dot) at (8,2.8){$c_{3}$};
			\draw[dashed](x4)--(y4);
			\draw[dashed](x5)--(y5);
			\draw[dashed](x6)--(y6);
			\coordinate (pp1) at (5.5,2);
			\coordinate (pp2) at (6,2);
			\coordinate (pp3) at (6,1);
			\coordinate (pp4) at (7,1);
			\draw(pp1)--(pp2);
			\draw(pp2)--(pp3);
			\draw(pp3)--(pp4);
			\node(small-dot) at (7,0.2){$F_{y}$};
			\coordinate (a2) at (10,0.5);
			\coordinate (b2) at (14,0.5);
			\coordinate (c2) at (10,2.5);
			\coordinate (d2) at (14,2.5);
			\draw(a2)--(b2);
			\draw(b2)--(d2);
			\draw(d2)--(c2);
			\draw(c2)--(a2);
			\coordinate (x7) at (11,0.5);
			\coordinate (y7) at (11,2.5);
			\coordinate (x8) at (12,0.5);
			\coordinate (y8) at (12,2.5);
			\coordinate (x9) at (13,0.5);
			\coordinate (y9) at (13,2.5);
			\node(small-dot) at (11,2.8){$c_{1}$};
			\node(small-dot) at (12,2.8){$c_{2}$};
			\node(small-dot) at (13,2.8){$c_{3}$};
			\draw[dashed](x7)--(y7);
			\draw[dashed](x8)--(y8);
			\draw[dashed](x9)--(y9);
			\coordinate (ppp1) at (10.5,2);
			\coordinate (ppp2) at (13.5,2);
			\coordinate (ppp3) at (10.5,1);
			\coordinate (ppp4) at (13.5,1);
			\draw(ppp1)--(ppp2);
			\draw(ppp3)--(ppp4);
			\node(small-dot) at (12,0.2){$F_{z}$};
			\node(small-dot) at (15.5,1.5){$\bar{C}$};
			\node(small-dot) at (2,2){$\circ$};
			\node(small-dot) at (1.7,1.7){$C_{1}$};
			\node(small-dot) at (3,2){$\circ$};
			\node(small-dot) at (2.7,1.7){$C_{2}$};
			\node(small-dot) at (6,2){$\circ$};
			\node(small-dot) at (5.7,1.7){$C_{3}$};
			\node(small-dot) at (15.5,0.8){\tiny{$1\le c_{1}\le c_{2}\le c_{3}\le n$}};
		\end{tikzpicture}
	\end{figure}	\FloatBarrier
	\begin{figure}[!h]
		\centering
		\begin{tikzpicture}
			\coordinate (a) at (0,0.5);
			\coordinate (b) at (4,0.5);
			\coordinate (c) at (0,2.5);
			\coordinate (d) at (4,2.5);
			\draw(a)--(b);
			\draw(b)--(d);
			\draw(d)--(c);
			\draw(c)--(a);
			\coordinate (x1) at (1,0.5);
			\coordinate (y1) at (1,2.5);
			\coordinate (x2) at (2,0.5);
			\coordinate (y2) at (2,2.5);
			\coordinate (x3) at (3,0.5);
			\coordinate (y3) at (3,2.5);
			\node(small-dot) at (1,2.8){$d_{1}$};
			\node(small-dot) at (2,2.8){$d_{2}$};
			\node(small-dot) at (3,2.8){$d_{3}$};
			\draw[dashed](x1)--(y1);
			\draw[dashed](x2)--(y2);
			\draw[dashed](x3)--(y3);
			\coordinate (p1) at (1,1);
			\coordinate (p2) at (3.5,1);
			\draw(p1)--(p2);
			\node(small-dot) at (2,0.2){$F_{x}$};
			\coordinate (a1) at (5,0.5);
			\coordinate (b1) at (9,0.5);
			\coordinate (c1) at (5,2.5);
			\coordinate (d1) at (9,2.5);
			\draw(a1)--(b1);
			\draw(b1)--(d1);
			\draw(d1)--(c1);
			\draw(c1)--(a1);
			\coordinate (x4) at (6,0.5);
			\coordinate (y4) at (6,2.5);
			\coordinate (x5) at (7,0.5);
			\coordinate (y5) at (7,2.5);
			\coordinate (x6) at (8,0.5);
			\coordinate (y6) at (8,2.5);
			\node(small-dot) at (6,2.8){$d_{1}$};
			\node(small-dot) at (7,2.8){$d_{2}$};
			\node(small-dot) at (8,2.8){$d_{3}$};
			\draw[dashed](x4)--(y4);
			\draw[dashed](x5)--(y5);
			\draw[dashed](x6)--(y6);
			\coordinate (pp1) at (5.5,2);
			\coordinate (pp2) at (6,2);
			\coordinate (pp3) at (7,2);
			\coordinate (pp4) at (8,2);
			\coordinate (pp5) at (8,1);
			\coordinate (pp6) at (8.5,1);
			\draw(pp1)--(pp2);
			\draw(pp3)--(pp4);
			\draw(pp4)--(pp5);
			\draw(pp5)--(pp6);
			\node(small-dot) at (7,0.2){$F_{y}$};
			\coordinate (a2) at (10,0.5);
			\coordinate (b2) at (14,0.5);
			\coordinate (c2) at (10,2.5);
			\coordinate (d2) at (14,2.5);
			\draw(a2)--(b2);
			\draw(b2)--(d2);
			\draw(d2)--(c2);
			\draw(c2)--(a2);
			\coordinate (x7) at (11,0.5);
			\coordinate (y7) at (11,2.5);
			\coordinate (x8) at (12,0.5);
			\coordinate (y8) at (12,2.5);
			\coordinate (x9) at (13,0.5);
			\coordinate (y9) at (13,2.5);
			\node(small-dot) at (11,2.8){$d_{1}$};
			\node(small-dot) at (12,2.8){$d_{2}$};
			\node(small-dot) at (13,2.8){$d_{3}$};
			\draw[dashed](x7)--(y7);
			\draw[dashed](x8)--(y8);
			\draw[dashed](x9)--(y9);
			\coordinate (ppp1) at (10.5,2);
			\coordinate (ppp2) at (13.5,2);
			\coordinate (ppp3) at (10.5,1);
			\coordinate (ppp4) at (12,1);
			\draw(ppp1)--(ppp2);
			\draw(ppp3)--(ppp4);
			\node(small-dot) at (12,0.2){$F_{z}$};
			\node(small-dot) at (15.5,1.5){$\bar{D}$};
			\node(small-dot) at (1,1){$\circ$};
			\node(small-dot) at (0.7,1.3){$D_{1}$};
			\node(small-dot) at (7,2){$\circ$};
			\node(small-dot) at (6.7,1.7){$D_{2}$};
			\node(small-dot) at (8,2){$\circ$};
			\node(small-dot) at (7.7,1.7){$D_{3}$};
			\node(small-dot) at (15.5,0.8){\tiny{$1\le d_{1}< d_{2}\le d_{3}\le n$}};
		\end{tikzpicture}
	\end{figure}	\FloatBarrier
	\begin{figure}[!h]
		\centering
		\begin{tikzpicture}
			\coordinate (a) at (0,0.5);
			\coordinate (b) at (4,0.5);
			\coordinate (c) at (0,2.5);
			\coordinate (d) at (4,2.5);
			\draw(a)--(b);
			\draw(b)--(d);
			\draw(d)--(c);
			\draw(c)--(a);
			\coordinate (x1) at (1,0.5);
			\coordinate (y1) at (1,2.5);
			\coordinate (x2) at (2,0.5);
			\coordinate (y2) at (2,2.5);
			\coordinate (x3) at (3,0.5);
			\coordinate (y3) at (3,2.5);
			\node(small-dot) at (1,2.8){$e_{1}$};
			\node(small-dot) at (2,2.8){$e_{2}$};
			\node(small-dot) at (3,2.8){$e_{3}$};
			\draw[dashed](x1)--(y1);
			\draw[dashed](x2)--(y2);
			\draw[dashed](x3)--(y3);
			\coordinate (p1) at (2,1);
			\coordinate (p2) at (3.5,1);
			\draw(p1)--(p2);
			\node(small-dot) at (2,0.2){$F_{x}$};
			\coordinate (a1) at (5,0.5);
			\coordinate (b1) at (9,0.5);
			\coordinate (c1) at (5,2.5);
			\coordinate (d1) at (9,2.5);
			\draw(a1)--(b1);
			\draw(b1)--(d1);
			\draw(d1)--(c1);
			\draw(c1)--(a1);
			\coordinate (x4) at (6,0.5);
			\coordinate (y4) at (6,2.5);
			\coordinate (x5) at (7,0.5);
			\coordinate (y5) at (7,2.5);
			\coordinate (x6) at (8,0.5);
			\coordinate (y6) at (8,2.5);
			\node(small-dot) at (6,2.8){$e_{1}$};
			\node(small-dot) at (7,2.8){$e_{2}$};
			\node(small-dot) at (8,2.8){$e_{3}$};
			\draw[dashed](x4)--(y4);
			\draw[dashed](x5)--(y5);
			\draw[dashed](x6)--(y6);
			\coordinate (pp1) at (5.5,2);
			\coordinate (pp2) at (6,2);
			\coordinate (pp3) at (6,1);
			\coordinate (pp4) at (7,1);
			\coordinate (pp5) at (8,1);
			\coordinate (pp6) at (8.5,1);
			\draw(pp1)--(pp2);
			\draw(pp2)--(pp3);
			\draw(pp3)--(pp4);
			\draw(pp5)--(pp6);
			\node(small-dot) at (7,0.2){$F_{y}$};
			\coordinate (a2) at (10,0.5);
			\coordinate (b2) at (14,0.5);
			\coordinate (c2) at (10,2.5);
			\coordinate (d2) at (14,2.5);
			\draw(a2)--(b2);
			\draw(b2)--(d2);
			\draw(d2)--(c2);
			\draw(c2)--(a2);
			\coordinate (x7) at (11,0.5);
			\coordinate (y7) at (11,2.5);
			\coordinate (x8) at (12,0.5);
			\coordinate (y8) at (12,2.5);
			\coordinate (x9) at (13,0.5);
			\coordinate (y9) at (13,2.5);
			\node(small-dot) at (11,2.8){$e_{1}$};
			\node(small-dot) at (12,2.8){$e_{2}$};
			\node(small-dot) at (13,2.8){$e_{3}$};
			\draw[dashed](x7)--(y7);
			\draw[dashed](x8)--(y8);
			\draw[dashed](x9)--(y9);
			\coordinate (ppp1) at (10.5,2);
			\coordinate (ppp2) at (13.5,2);
			\coordinate (ppp3) at (10.5,1);
			\coordinate (ppp4) at (13,1);
			\draw(ppp1)--(ppp2);
			\draw(ppp3)--(ppp4);
			\node(small-dot) at (12,0.2){$F_{z}$};
			\node(small-dot) at (15.5,1.5){$\bar{E}$};
			\node(small-dot) at (2,1){$\circ$};
			\node(small-dot) at (1.7,1.3){$E_{1}$};
			\node(small-dot) at (6,2){$\circ$};
			\node(small-dot) at (5.7,1.7){$E_{2}$};
			\node(small-dot) at (8,1){$\circ$};
			\node(small-dot) at (7.7,1.3){$E_{3}$};
			\node(small-dot) at (15.5,0.8){\tiny{$1\le e_{1}\le e_{2}< e_{3}\le n$}};
		\end{tikzpicture}
	\end{figure}	\FloatBarrier
	\begin{figure}[!h]
		\centering
		\begin{tikzpicture}
			\coordinate (a) at (0,0.5);
			\coordinate (b) at (4,0.5);
			\coordinate (c) at (0,2.5);
			\coordinate (d) at (4,2.5);
			\draw(a)--(b);
			\draw(b)--(d);
			\draw(d)--(c);
			\draw(c)--(a);
			\coordinate (x1) at (1,0.5);
			\coordinate (y1) at (1,2.5);
			\coordinate (x2) at (2,0.5);
			\coordinate (y2) at (2,2.5);
			\coordinate (x3) at (3,0.5);
			\coordinate (y3) at (3,2.5);
			\node(small-dot) at (1,2.8){$f_{1}$};
			\node(small-dot) at (2,2.8){$f_{2}$};
			\node(small-dot) at (3,2.8){$f_{3}$};
			\draw[dashed](x1)--(y1);
			\draw[dashed](x2)--(y2);
			\draw[dashed](x3)--(y3);
			\coordinate (p1) at (1,1);
			\coordinate (p2) at (3.5,1);
			\draw(p1)--(p2);
			\node(small-dot) at (2,0.2){$F_{x}$};
			\coordinate (a1) at (5,0.5);
			\coordinate (b1) at (9,0.5);
			\coordinate (c1) at (5,2.5);
			\coordinate (d1) at (9,2.5);
			\draw(a1)--(b1);
			\draw(b1)--(d1);
			\draw(d1)--(c1);
			\draw(c1)--(a1);
			\coordinate (x4) at (6,0.5);
			\coordinate (y4) at (6,2.5);
			\coordinate (x5) at (7,0.5);
			\coordinate (y5) at (7,2.5);
			\coordinate (x6) at (8,0.5);
			\coordinate (y6) at (8,2.5);
			\node(small-dot) at (6,2.8){$f_{1}$};
			\node(small-dot) at (7,2.8){$f_{2}$};
			\node(small-dot) at (8,2.8){$f_{3}$};
			\draw[dashed](x4)--(y4);
			\draw[dashed](x5)--(y5);
			\draw[dashed](x6)--(y6);
			\coordinate (pp1) at (5.5,1);
			\coordinate (pp2) at (6,1);
			\coordinate (pp3) at (7,1);
			\coordinate (pp4) at (8.5,1);
			\draw(pp1)--(pp2);
			\draw(pp3)--(pp4);
			\node(small-dot) at (7,0.2){$F_{y}$};
			\coordinate (a2) at (10,0.5);
			\coordinate (b2) at (14,0.5);
			\coordinate (c2) at (10,2.5);
			\coordinate (d2) at (14,2.5);
			\draw(a2)--(b2);
			\draw(b2)--(d2);
			\draw(d2)--(c2);
			\draw(c2)--(a2);
			\coordinate (x7) at (11,0.5);
			\coordinate (y7) at (11,2.5);
			\coordinate (x8) at (12,0.5);
			\coordinate (y8) at (12,2.5);
			\coordinate (x9) at (13,0.5);
			\coordinate (y9) at (13,2.5);
			\node(small-dot) at (11,2.8){$f_{1}$};
			\node(small-dot) at (12,2.8){$f_{2}$};
			\node(small-dot) at (13,2.8){$f_{3}$};
			\draw[dashed](x7)--(y7);
			\draw[dashed](x8)--(y8);
			\draw[dashed](x9)--(y9);
			\coordinate (ppp1) at (10.5,2);
			\coordinate (ppp2) at (13,2);
			\coordinate (ppp3) at (13,1);
			\coordinate (ppp4) at (13.5,1);
			\coordinate (ppp5) at (10.5,1);
			\coordinate (ppp6) at (12,1);
			\draw(ppp1)--(ppp2);
			\draw(ppp2)--(ppp3);
			\draw(ppp3)--(ppp4);
			\draw(ppp5)--(ppp6);
			\node(small-dot) at (12,0.2){$F_{z}$};
			\node(small-dot) at (15.5,1.5){$\bar{F}$};
			\node(small-dot) at (1,1){$\circ$};
			\node(small-dot) at (0.7,1.3){$F_{1}$};
			\node(small-dot) at (7,1){$\circ$};
			\node(small-dot) at (6.7,1.3){$F_{2}$};
			\node(small-dot) at (13,2){$\circ$};
			\node(small-dot) at (12.7,1.7){$F_{3}$};
			\node(small-dot) at (15.5,0.8){\tiny{$1\le f_{1}<f_{2}< f_{3}\le n$}};
		\end{tikzpicture}
	\end{figure}	\FloatBarrier
	\par These points will vanish in some non-generic conditions due to boundary effects, which is detailed in the following table:
	\begin{figure}[!h]
		\centering
		\begin{tikzpicture}
			\coordinate (a1) at (0,5);
			\coordinate (b1) at (15,5);
			\coordinate (a2) at (0,4);
			\coordinate (b2) at (15,4);
			\coordinate (a3) at (0,3);
			\coordinate (b3) at (15,3);
			\coordinate (a4) at (0,2);
			\coordinate (b4) at (15,2);
			\coordinate (a5) at (0,1);
			\coordinate (b5) at (15,1);
			\coordinate (a6) at (0,0);
			\coordinate (b6) at (15,0);
			\draw(a1)--(b1);
			\draw(a2)--(b2);
			\draw(a3)--(b3);
			\draw(a4)--(b4);
			\draw(a5)--(b5);
			\draw(a6)--(b6);
			\coordinate (c1) at (1,5);
			\coordinate (c2) at (5,5);
			\coordinate (c3) at (6,5);
			\coordinate (c4) at (10,5);
			\coordinate (c5) at (11,5);
			\coordinate (d1) at (1,0);
			\coordinate (d2) at (5,0);
			\coordinate (d3) at (6,0);
			\coordinate (d4) at (10,0);
			\coordinate (d5) at (11,0);
			\draw(a1)--(a6);
			\draw(b1)--(b6);
			\draw(c3)--(d3);
			\draw(c4)--(d4);
			\draw(c5)--(d5);
			\draw(c2)--(d2);
			\draw(c1)--(d1);
			\node(small-dot) at (0.5,4.5){$A_{1}$};
			\node(small-dot) at (5.5,4.5){$A_{2}$};
			\node(small-dot) at (10.5,4.5){$A_{3}$};
			\node(small-dot) at (2.3,4.75){\ding{172}$a_{1}=a_{r}$ or}; \node(small-dot) at (3,4.25){\ding{173}$a_{1}=n-1$, $a_{2}=n$};
			\node(small-dot) at (7,4.5){$a_{1}=a_{r}$};
			\node(small-dot) at (11.9,4.5){$a_{2}=n$};
			\node(small-dot) at (0.5,3.5){$C_{1}$};
			\node(small-dot) at (5.5,3.5){$C_{2}$};
			\node(small-dot) at (10.5,3.5){$C_{3}$};
			\node(small-dot) at (1.8,3.5){$c_{2}=c_{3}$};
			\node(small-dot) at (7,3.5){$c_{2}=c_{3}$};
			\node(small-dot) at (11.93,3.5){$c_{1}=1$};
			\node(small-dot) at (0.5,2.5){$D_{1}$};
			\node(small-dot) at (5.5,2.5){$D_{2}$};
			\node(small-dot) at (10.5,2.5){$D_{3}$};
			\node(small-dot) at (3,2.5){$d_{1}=d_{2}-1$, $d_{2}=d_{3}$};
			\node(small-dot) at (7,2.5){$d_{2}=d_{3}$};
			\node(small-dot) at (12,2.5){$d_{2}=d_{3}$};
			\node(small-dot) at (0.5,1.5){$E_{1}$};
			\node(small-dot) at (5.5,1.5){$E_{2}$};
			\node(small-dot) at (10.5,1.5){$E_{3}$};
			\node(small-dot) at (2.85,1.5){$e_{2}=n-1$, $e_{3}=n$};
			\node(small-dot) at (6.93,1.5){$e_{1}=1$};
			\node(small-dot) at (11.95,1.5){$e_{3}=n$};
			\node(small-dot) at (0.5,0.5){$F_{1}$};
			\node(small-dot) at (5.5,0.5){$F_{2}$};
			\node(small-dot) at (10.5,0.5){$F_{3}$};
			\node(small-dot) at (2.2,0.75){$f_{1}=n-2$,};
			\node(small-dot) at (2.85,0.25){$f_{2}=n-1$, $f_{3}=n$};
			\node(small-dot) at (8,0.5){$f_{2}=n-1$, $f_{3}=n$};
			\node(small-dot) at (12.38,0.5){$f_{2}=f_{3}-1$};
			\node(small-dot) at (7.5,5.5){Vanishing conditions};
		\end{tikzpicture}
	\end{figure}	\FloatBarrier
	\par We consider $c_{2}(P)$ for every facet $P\in \Delta_{0}$ following by checking all possible cases. We will use the same label if the position of a point in $c_{2}(P)$ is at the same position as some point in $c_{1}(P)$. Following the upper path represents $P$ and the lower path represents $Q$. Sets $P-Q$ and $Q-P$ both only have one element. We use a hollow circle to represent points in $c_{2}(P)$, which is just the unique element in $P-Q$.
	\par \noindent \textit{Case $P\in \bar{E}$, $Q\in \bar{F}$.} $y_{1,1}\in P$ and $y_{1,1}\notin Q$ imply $P-Q=\{y_{1,1}\}$. Hence $e_{1}=1$ and $f_{3}=n$. If $e_{3}<n$, then $Q-P=\{z_{2,n}\}$. We have $e_{3}=f_{2}$ and $e_{2}=f_{1}$.
	
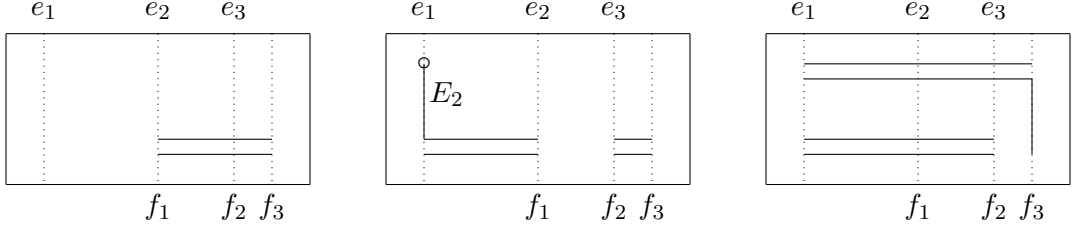
\begin{figure}[h]
		\centering
		\begin{tikzpicture}
			\coordinate (a) at (0,0.5);
			\coordinate (b) at (4,0.5);
			\coordinate (c) at (0,2.5);
			\coordinate (d) at (4,2.5);
			\draw(a)--(b);
			\draw(b)--(d);
			\draw(d)--(c);
			\draw(c)--(a);
			\node(small-dot) at (0.5,2.8){$e_{1}$};
			\node(small-dot) at (2,2.8){$e_{2}$};
			\node(small-dot) at (3,2.8){$e_{3}$};
			\node(small-dot) at (2,0.2){$f_{1}$};
			\node(small-dot) at (3,0.2){$f_{2}$};
			\node(small-dot) at (3.5,0.2){$f_{3}$};
			\coordinate (p1) at (2,1.1);
			\coordinate (p2) at (3.5,1.1);
			\coordinate (q1) at (2,0.9);
			\coordinate (q2) at (3.5,0.9);
			\draw(p1)--(p2);
			\draw(q1)--(q2);
			\coordinate (a1) at (5,0.5);
			\coordinate (b1) at (9,0.5);
			\coordinate (c1) at (5,2.5);
			\coordinate (d1) at (9,2.5);
			\draw(a1)--(b1);
			\draw(b1)--(d1);
			\draw(d1)--(c1);
			\draw(c1)--(a1);
			\node(small-dot) at (5.5,2.8){$e_{1}$};
			\node(small-dot) at (7,2.8){$e_{2}$};
			\node(small-dot) at (8,2.8){$e_{3}$};
			\node(small-dot) at (7,0.2){$f_{1}$};
			\node(small-dot) at (8,0.2){$f_{2}$};
			\node(small-dot) at (8.5,0.2){$f_{3}$};
			\coordinate (pp1) at (5.5,1.1);
			\coordinate (pp2) at (7,1.1);
			\coordinate (pp3) at (8,1.1);
			\coordinate (pp4) at (8.5,1.1);
			\coordinate (pp5) at (5.5,2.1);
			\coordinate (qq1) at (5.5,0.9);
			\coordinate (qq2) at (7,0.9);
			\coordinate (qq3) at (8,0.9);
			\coordinate (qq4) at (8.5,0.9);
			\draw(pp1)--(pp2);
			\draw(pp3)--(pp4);
			\draw(pp1)--(pp5);
			\draw(qq1)--(qq2);
			\draw(qq3)--(qq4);
			\coordinate (a2) at (10,0.5);
			\coordinate (b2) at (14,0.5);
			\coordinate (c2) at (10,2.5);
			\coordinate (d2) at (14,2.5);
			\draw(a2)--(b2);
			\draw(b2)--(d2);
			\draw(d2)--(c2);
			\draw(c2)--(a2);
			\node(small-dot) at (10.5,2.8){$e_{1}$};
			\node(small-dot) at (12,2.8){$e_{2}$};
			\node(small-dot) at (13,2.8){$e_{3}$};
			\node(small-dot) at (12,0.2){$f_{1}$};
			\node(small-dot) at (13,0.2){$f_{2}$};
			\node(small-dot) at (13.5,0.2){$f_{3}$};
			\coordinate (ppp1) at (10.5,2.1);
			\coordinate (ppp2) at (13.5,2.1);
			\coordinate (ppp3) at (10.5,1.1);
			\coordinate (ppp4) at (13,1.1);
			\coordinate (qqq1) at (10.5,1.9);
			\coordinate (qqq2) at (13.5,1.9);
			\coordinate (qqq3) at (10.5,0.9);
			\coordinate (qqq4) at (13,0.9);
			\coordinate (qqq5) at (13.5,0.9);
			\draw(ppp1)--(ppp2);
			\draw(ppp3)--(ppp4);
			\draw(qqq1)--(qqq2);
			\draw(qqq3)--(qqq4);
			\draw(qqq2)--(qqq5);
			\node(small-dot) at (5.5,2.1){$\circ$};
			\node(small-dot) at (5.8,1.7){$E_{2}$};
			\coordinate (s1) at (0.5,2.5);
			\coordinate (t1) at (0.5,0.5);
			\draw[dotted](s1)--(t1);
			\coordinate (s2) at (2,2.5);
			\coordinate (t2) at (2,0.5);
			\draw[dotted](s2)--(t2);
			\coordinate (s3) at (3,2.5);
			\coordinate (t3) at (3,0.5);
			\draw[dotted](s3)--(t3);
			\coordinate (s4) at (3.5,2.5);
			\coordinate (t4) at (3.5,0.5);
			\draw[dotted](s4)--(t4);
			\coordinate (ss1) at (5.5,2.5);
			\coordinate (tt1) at (5.5,0.5);
			\draw[dotted](ss1)--(tt1);
			\coordinate (ss2) at (7,2.5);
			\coordinate (tt2) at (7,0.5);
			\draw[dotted](ss2)--(tt2);
			\coordinate (ss3) at (8,2.5);
			\coordinate (tt3) at (8,0.5);
			\draw[dotted](ss3)--(tt3);
			\coordinate (ss4) at (8.5,2.5);
			\coordinate (tt4) at (8.5,0.5);
			\draw[dotted](ss4)--(tt4);
			\coordinate (sss1) at (10.5,2.5);
			\coordinate (ttt1) at (10.5,0.5);
			\draw[dotted](sss1)--(ttt1);
			\coordinate (sss2) at (12,2.5);
			\coordinate (ttt2) at (12,0.5);
			\draw[dotted](sss2)--(ttt2);
			\coordinate (sss3) at (13,2.5);
			\coordinate (ttt3) at (13,0.5);
			\draw[dotted](sss3)--(ttt3);
			\coordinate (sss4) at (13.5,2.5);
			\coordinate (ttt4) at (13.5,0.5);
			\draw[dotted](sss4)--(ttt4);
		\end{tikzpicture}
	\end{figure}	\FloatBarrier 
	\par If $e_{3}=n$, then $f_{2}=n-1$. We have $Q-P=\{ y_{2,n-1}\}$ and $e_{2}=f_{1}$ in this case. 
	
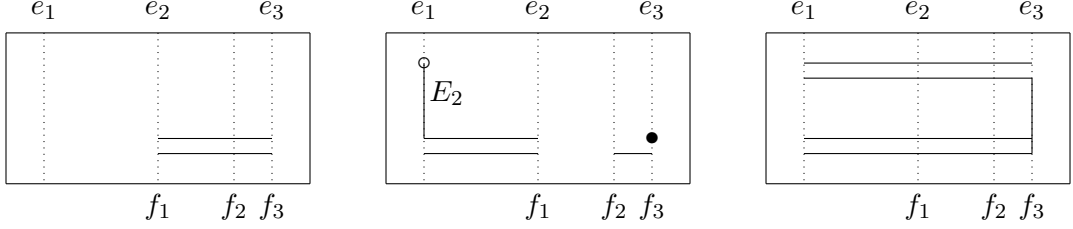
\begin{figure}[h]
		\centering
		\begin{tikzpicture}
			\coordinate (a) at (0,0.5);
			\coordinate (b) at (4,0.5);
			\coordinate (c) at (0,2.5);
			\coordinate (d) at (4,2.5);
			\draw(a)--(b);
			\draw(b)--(d);
			\draw(d)--(c);
			\draw(c)--(a);
			\node(small-dot) at (0.5,2.8){$e_{1}$};
			\node(small-dot) at (2,2.8){$e_{2}$};
			\node(small-dot) at (3.5,2.8){$e_{3}$};
			\node(small-dot) at (2,0.2){$f_{1}$};
			\node(small-dot) at (3,0.2){$f_{2}$};
			\node(small-dot) at (3.5,0.2){$f_{3}$};
			\coordinate (p1) at (2,1.1);
			\coordinate (p2) at (3.5,1.1);
			\coordinate (q1) at (2,0.9);
			\coordinate (q2) at (3.5,0.9);
			\draw(p1)--(p2);
			\draw(q1)--(q2);
			\coordinate (a1) at (5,0.5);
			\coordinate (b1) at (9,0.5);
			\coordinate (c1) at (5,2.5);
			\coordinate (d1) at (9,2.5);
			\draw(a1)--(b1);
			\draw(b1)--(d1);
			\draw(d1)--(c1);
			\draw(c1)--(a1);
			\node(small-dot) at (5.5,2.8){$e_{1}$};
			\node(small-dot) at (7,2.8){$e_{2}$};
			\node(small-dot) at (8.5,2.8){$e_{3}$};
			\node(small-dot) at (7,0.2){$f_{1}$};
			\node(small-dot) at (8,0.2){$f_{2}$};
			\node(small-dot) at (8.5,0.2){$f_{3}$};
			\coordinate (pp1) at (5.5,1.1);
			\coordinate (pp2) at (7,1.1);
			\coordinate (pp5) at (5.5,2.1);
			\coordinate (qq1) at (5.5,0.9);
			\coordinate (qq2) at (7,0.9);
			\coordinate (qq3) at (8,0.9);
			\coordinate (qq4) at (8.5,0.9);
			\draw(pp1)--(pp2);
			\draw(pp1)--(pp5);
			\draw(qq1)--(qq2);
			\draw(qq3)--(qq4);
			\coordinate (a2) at (10,0.5);
			\coordinate (b2) at (14,0.5);
			\coordinate (c2) at (10,2.5);
			\coordinate (d2) at (14,2.5);
			\draw(a2)--(b2);
			\draw(b2)--(d2);
			\draw(d2)--(c2);
			\draw(c2)--(a2);
			\node(small-dot) at (10.5,2.8){$e_{1}$};
			\node(small-dot) at (12,2.8){$e_{2}$};
			\node(small-dot) at (13.5,2.8){$e_{3}$};
			\node(small-dot) at (12,0.2){$f_{1}$};
			\node(small-dot) at (13,0.2){$f_{2}$};
			\node(small-dot) at (13.5,0.2){$f_{3}$};
			\coordinate (ppp1) at (10.5,2.1);
			\coordinate (ppp2) at (13.5,2.1);
			\coordinate (ppp3) at (10.5,1.1);
			\coordinate (ppp4) at (13.5,1.1);
			\coordinate (qqq1) at (10.5,1.9);
			\coordinate (qqq2) at (13.5,1.9);
			\coordinate (qqq3) at (10.5,0.9);
			\coordinate (qqq4) at (13.5,0.9);
			\coordinate (qqq5) at (13.5,0.9);
			\draw(ppp1)--(ppp2);
			\draw(ppp3)--(ppp4);
			\draw(qqq1)--(qqq2);
			\draw(qqq3)--(qqq4);
			\draw(qqq2)--(qqq5);
			\node(small-dot) at (5.5,2.1){$\circ$};
			\node(small-dot) at (5.8,1.7){$E_{2}$};
			\node(small-dot) at (8.5,1.1){$\bullet$};
			\coordinate (s1) at (0.5,2.5);
			\coordinate (t1) at (0.5,0.5);
			\draw[dotted](s1)--(t1);
			\coordinate (s2) at (2,2.5);
			\coordinate (t2) at (2,0.5);
			\draw[dotted](s2)--(t2);
			\coordinate (s3) at (3,2.5);
			\coordinate (t3) at (3,0.5);
			\draw[dotted](s3)--(t3);
			\coordinate (s4) at (3.5,2.5);
			\coordinate (t4) at (3.5,0.5);
			\draw[dotted](s4)--(t4);
			\coordinate (ss1) at (5.5,2.5);
			\coordinate (tt1) at (5.5,0.5);
			\draw[dotted](ss1)--(tt1);
			\coordinate (ss2) at (7,2.5);
			\coordinate (tt2) at (7,0.5);
			\draw[dotted](ss2)--(tt2);
			\coordinate (ss3) at (8,2.5);
			\coordinate (tt3) at (8,0.5);
			\draw[dotted](ss3)--(tt3);
			\coordinate (ss4) at (8.5,2.5);
			\coordinate (tt4) at (8.5,0.5);
			\draw[dotted](ss4)--(tt4);
			\coordinate (sss1) at (10.5,2.5);
			\coordinate (ttt1) at (10.5,0.5);
			\draw[dotted](sss1)--(ttt1);
			\coordinate (sss2) at (12,2.5);
			\coordinate (ttt2) at (12,0.5);
			\draw[dotted](sss2)--(ttt2);
			\coordinate (sss3) at (13,2.5);
			\coordinate (ttt3) at (13,0.5);
			\draw[dotted](sss3)--(ttt3);
			\coordinate (sss4) at (13.5,2.5);
			\coordinate (ttt4) at (13.5,0.5);
			\draw[dotted](sss4)--(ttt4);
		\end{tikzpicture}
	\end{figure}	\FloatBarrier
	\par \noindent \textit{Case $P\in \bar{D}$, $Q\in \bar{F}$.} $P-Q$ contains at least two points $y_{1,1}$ and $y_{1,d_{3}}$, contradiction.
	\par \noindent \textit{Case $P\in \bar{D}$, $Q\in \bar{E}$.} This case is a little complicated. There are three conditions $v\in P_{x}$, $v\in P_{y}$, and $v\in P_{z}$ with supposition $P-Q=\{v\}$. 
	\par If $v\in P_{x}$, then $d_{1}=e_{2}-1$ and $d_{3}=e_{1}$. We have $d_{1}\le d_{3}-1=e_{1}-1\le e_{2}-1=d_{1}$, so $d_{1}=d_{2}-1=d_{3}-1$ and $e_{1}=e_{2}$. There is no point in $P_{y}-Q_{y}$, so $e_{3}=e_{2}+1$.
	
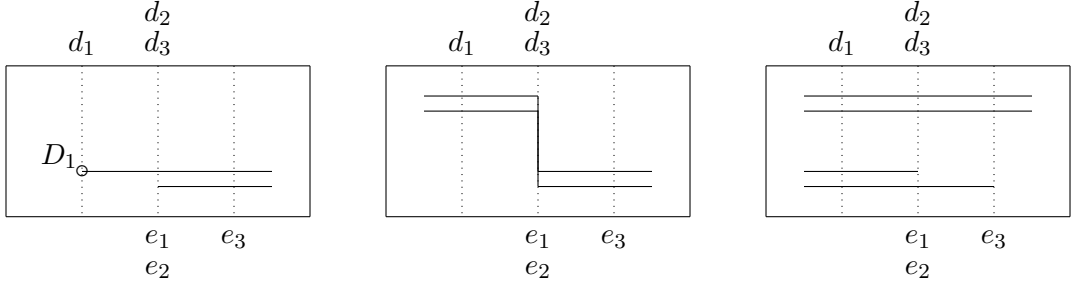
\begin{figure}[h]
		\centering
		\begin{tikzpicture}
			\coordinate (a) at (0,0.5);
			\coordinate (b) at (4,0.5);
			\coordinate (c) at (0,2.5);
			\coordinate (d) at (4,2.5);
			\draw(a)--(b);
			\draw(b)--(d);
			\draw(d)--(c);
			\draw(c)--(a);
			\node(small-dot) at (1,2.8){$d_{1}$};
			\node(small-dot) at (2,3.2){$d_{2}$};
			\node(small-dot) at (2,2.8){$d_{3}$};
			\node(small-dot) at (2,0.2){$e_{1}$};
			\node(small-dot) at (2,-0.2){$e_{2}$};
			\node(small-dot) at (3,0.2){$e_{3}$};
			\coordinate (p1) at (1,1.1);
			\coordinate (p2) at (3.5,1.1);
			\coordinate (q1) at (2,0.9);
			\coordinate (q2) at (3.5,0.9);
			\draw(p1)--(p2);
			\draw(q1)--(q2);
			\coordinate (a1) at (5,0.5);
			\coordinate (b1) at (9,0.5);
			\coordinate (c1) at (5,2.5);
			\coordinate (d1) at (9,2.5);
			\draw(a1)--(b1);
			\draw(b1)--(d1);
			\draw(d1)--(c1);
			\draw(c1)--(a1);
			\node(small-dot) at (6,2.8){$d_{1}$};
			\node(small-dot) at (7,3.2){$d_{2}$};
			\node(small-dot) at (7,2.8){$d_{3}$};
			\node(small-dot) at (7,0.2){$e_{1}$};
			\node(small-dot) at (7,-0.2){$e_{2}$};
			\node(small-dot) at (8,0.2){$e_{3}$};
			\coordinate (pp1) at (5.5,2.1);
			\coordinate (pp2) at (6,2.1);
			\coordinate (pp3) at (7,2.1);
			\coordinate (pp4) at (7,1.1);
			\coordinate (pp5) at (8.5,1.1);
			\coordinate (qq1) at (5.5,1.9);
			\coordinate (qq2) at (7,1.9);
			\coordinate (qq3) at (7,0.9);
			\coordinate (qq4) at (8,0.9);
			\coordinate (qq5) at (8.5,0.9);
			\draw(qq1)--(qq2);
			\draw(qq3)--(qq4);
			\draw(qq2)--(qq3);
			\draw(qq5)--(qq4);
			\draw(pp1)--(pp2);
			\draw(pp3)--(pp2);
			\draw(pp4)--(pp5);
			\draw(pp3)--(pp4);
			\coordinate (a2) at (10,0.5);
			\coordinate (b2) at (14,0.5);
			\coordinate (c2) at (10,2.5);
			\coordinate (d2) at (14,2.5);
			\draw(a2)--(b2);
			\draw(b2)--(d2);
			\draw(d2)--(c2);
			\draw(c2)--(a2);
			\node(small-dot) at (11,2.8){$d_{1}$};
			\node(small-dot) at (12,3.2){$d_{2}$};
			\node(small-dot) at (12,2.8){$d_{3}$};
			\node(small-dot) at (12,0.2){$e_{1}$};
			\node(small-dot) at (12,-0.2){$e_{2}$};
			\node(small-dot) at (13,0.2){$e_{3}$};
			\coordinate (ppp1) at (10.5,2.1);
			\coordinate (ppp2) at (13.5,2.1);
			\coordinate (ppp3) at (10.5,1.1);
			\coordinate (ppp4) at (12,1.1);
			\coordinate (qqq1) at (10.5,1.9);
			\coordinate (qqq2) at (13.5,1.9);
			\coordinate (qqq3) at (10.5,0.9);
			\coordinate (qqq4) at (13,0.9);
			\draw(ppp1)--(ppp2);
			\draw(ppp3)--(ppp4);
			\draw(qqq1)--(qqq2);
			\draw(qqq3)--(qqq4);
			\node(small-dot) at (1,1.1){$\circ$};
			\node(small-dot) at (0.7,1.3){$D_{1}$};
			\coordinate (s1) at (1,2.5);
			\coordinate (t1) at (1,0.5);
			\draw[dotted](s1)--(t1);
			\coordinate (s2) at (2,2.5);
			\coordinate (t2) at (2,0.5);
			\draw[dotted](s2)--(t2);
			\coordinate (s3) at (3,2.5);
			\coordinate (t3) at (3,0.5);
			\draw[dotted](s3)--(t3);
			\coordinate (ss1) at (6,2.5);
			\coordinate (tt1) at (6,0.5);
			\draw[dotted](ss1)--(tt1);
			\coordinate (ss2) at (7,2.5);
			\coordinate (tt2) at (7,0.5);
			\draw[dotted](ss2)--(tt2);
			\coordinate (ss3) at (8,2.5);
			\coordinate (tt3) at (8,0.5);
			\draw[dotted](ss3)--(tt3);
			\coordinate (sss1) at (11,2.5);
			\coordinate (ttt1) at (11,0.5);
			\draw[dotted](sss1)--(ttt1);
			\coordinate (sss2) at (12,2.5);
			\coordinate (ttt2) at (12,0.5);
			\draw[dotted](sss2)--(ttt2);
			\coordinate (sss3) at (13,2.5);
			\coordinate (ttt3) at (13,0.5);
			\draw[dotted](sss3)--(ttt3);
		\end{tikzpicture}
	\end{figure}	\FloatBarrier
	\par If $v\in P_{y}$, then we have $d_{1}\ge e_{2}$. The point $y_{2,e_{1}}$ is in $Q-P$ and the point $y_{1,d_{3}}$ is in $P-Q$. The diagrams of $P_{x}$ and $Q_{x}$ are the same, so $d_{1}=e_{2}$.
	The diagrams of $P_{z}$ and $Q_{z}$ are the same, so $d_{2}=e_{3}$. Hence $e_{1}=e_{2}$ and $d_{3}=d_{2}$ according to $y$-parts.
	
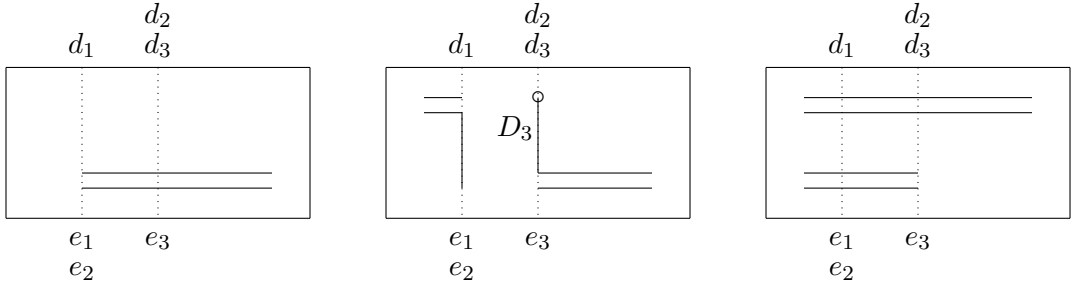
\begin{figure}[h]
		\centering
		\begin{tikzpicture}
			\coordinate (a) at (0,0.5);
			\coordinate (b) at (4,0.5);
			\coordinate (c) at (0,2.5);
			\coordinate (d) at (4,2.5);
			\draw(a)--(b);
			\draw(b)--(d);
			\draw(d)--(c);
			\draw(c)--(a);
			\node(small-dot) at (1,2.8){$d_{1}$};
			\node(small-dot) at (2,3.2){$d_{2}$};
			\node(small-dot) at (2,2.8){$d_{3}$};
			\node(small-dot) at (1,0.2){$e_{1}$};
			\node(small-dot) at (1,-0.2){$e_{2}$};
			\node(small-dot) at (2,0.2){$e_{3}$};
			\coordinate (p1) at (1,1.1);
			\coordinate (p2) at (3.5,1.1);
			\coordinate (q1) at (1,0.9);
			\coordinate (q2) at (3.5,0.9);
			\draw(p1)--(p2);
			\draw(q1)--(q2);
			\coordinate (a1) at (5,0.5);
			\coordinate (b1) at (9,0.5);
			\coordinate (c1) at (5,2.5);
			\coordinate (d1) at (9,2.5);
			\draw(a1)--(b1);
			\draw(b1)--(d1);
			\draw(d1)--(c1);
			\draw(c1)--(a1);
			\node(small-dot) at (6,2.8){$d_{1}$};
			\node(small-dot) at (7,3.2){$d_{2}$};
			\node(small-dot) at (7,2.8){$d_{3}$};
			\node(small-dot) at (6,0.2){$e_{1}$};
			\node(small-dot) at (6,-0.2){$e_{2}$};
			\node(small-dot) at (7,0.2){$e_{3}$};
			\coordinate (pp1) at (5.5,2.1);
			\coordinate (pp2) at (6,2.1);
			\coordinate (pp3) at (7,2.1);
			\coordinate (pp4) at (7,1.1);
			\coordinate (pp5) at (8.5,1.1);
			\coordinate (qq1) at (5.5,1.9);
			\coordinate (qq2) at (6,1.9);
			\coordinate (qq3) at (6,0.9);
			\coordinate (qq4) at (7,0.9);
			\coordinate (qq5) at (8.5,0.9);
			\draw(pp1)--(pp2);
			\draw(pp3)--(pp4);
			\draw(pp4)--(pp5);
			\draw(qq1)--(qq2);
			\draw(qq2)--(qq3);
			\draw(qq4)--(qq5);
			\coordinate (a2) at (10,0.5);
			\coordinate (b2) at (14,0.5);
			\coordinate (c2) at (10,2.5);
			\coordinate (d2) at (14,2.5);
			\draw(a2)--(b2);
			\draw(b2)--(d2);
			\draw(d2)--(c2);
			\draw(c2)--(a2);
			\node(small-dot) at (11,2.8){$d_{1}$};
			\node(small-dot) at (12,3.2){$d_{2}$};
			\node(small-dot) at (12,2.8){$d_{3}$};
			\node(small-dot) at (11,0.2){$e_{1}$};
			\node(small-dot) at (11,-0.2){$e_{2}$};
			\node(small-dot) at (12,0.2){$e_{3}$};
			\coordinate (ppp1) at (10.5,2.1);
			\coordinate (ppp2) at (13.5,2.1);
			\coordinate (ppp3) at (10.5,1.1);
			\coordinate (ppp4) at (12,1.1);
			\coordinate (qqq1) at (10.5,1.9);
			\coordinate (qqq2) at (13.5,1.9);
			\coordinate (qqq3) at (10.5,0.9);
			\coordinate (qqq4) at (12,0.9);
			\draw(ppp1)--(ppp2);
			\draw(ppp3)--(ppp4);
			\draw(qqq1)--(qqq2);
			\draw(qqq3)--(qqq4);
			\node(small-dot) at (7,2.1){$\circ$};
			\node(small-dot) at (6.7,1.7){$D_{3}$};
			\coordinate (s1) at (1,2.5);
			\coordinate (t1) at (1,0.5);
			\draw[dotted](s1)--(t1);
			\coordinate (s2) at (2,2.5);
			\coordinate (t2) at (2,0.5);
			\draw[dotted](s2)--(t2);
			\coordinate (s3) at (3,2.5);
			\coordinate (ss1) at (6,2.5);
			\coordinate (tt1) at (6,0.5);
			\draw[dotted](ss1)--(tt1);
			\coordinate (ss2) at (7,2.5);
			\coordinate (tt2) at (7,0.5);
			\draw[dotted](ss2)--(tt2);
			\coordinate (sss1) at (11,2.5);
			\coordinate (ttt1) at (11,0.5);
			\draw[dotted](sss1)--(ttt1);
			\coordinate (sss2) at (12,2.5);
			\coordinate (ttt2) at (12,0.5);
			\draw[dotted](sss2)--(ttt2);
		\end{tikzpicture}
	\end{figure}	\FloatBarrier
	\par If $v\in P_{z}$, then $d_{2}=e_{3}+1$ and $d_{3}=e_{1}$. We have $d_{2}\le d_{3}=e_{1}\le e_{3}<d_{2}$, contradiction.
	\par \noindent \textit{Case $P\in \bar{C}$, $Q\in \bar{F}$.} $P-Q$ contains at least two points $x_{1,c_{2}}$ and $y_{1,c_{1}}$, contradiction.
	\par \noindent \textit{Case $P\in \bar{C}$, $Q\in \bar{E}$.} $x_{1,c_{2}}\in P-Q$ implies $c_{2}=c_{3}$ and $e_{3}=n$. We have $c_{1}=e_{1}$ and $c_{2}=e_{2}$ according to $y$-parts.
	
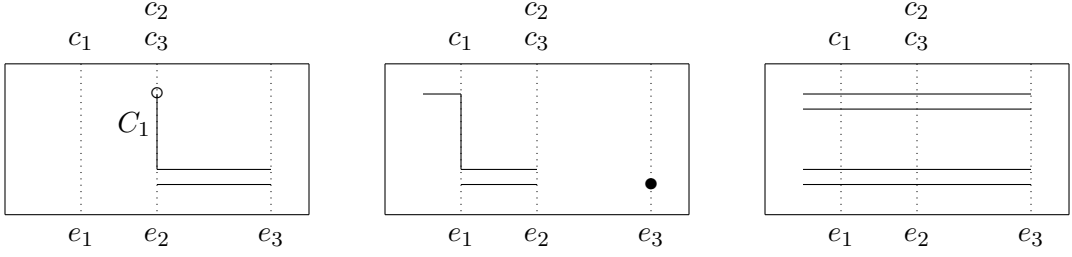
\begin{figure}[h]
		\centering
		\begin{tikzpicture}
			\coordinate (a) at (0,0.5);
			\coordinate (b) at (4,0.5);
			\coordinate (c) at (0,2.5);
			\coordinate (d) at (4,2.5);
			\draw(a)--(b);
			\draw(b)--(d);
			\draw(d)--(c);
			\draw(c)--(a);
			\node(small-dot) at (1,2.8){$c_{1}$};
			\node(small-dot) at (2,3.2){$c_{2}$};
			\node(small-dot) at (2,2.8){$c_{3}$};
			\node(small-dot) at (1,0.2){$e_{1}$};
			\node(small-dot) at (2,0.2){$e_{2}$};
			\node(small-dot) at (3.5,0.2){$e_{3}$};
			\coordinate (p1) at (2,2.1);
			\coordinate (p2) at (2,1.1);
			\coordinate (p3) at (3.5,1.1);
			\coordinate (q1) at (2,0.9);
			\coordinate (q2) at (3.5,0.9);
			\draw(p1)--(p2);
			\draw(p2)--(p3);
			\draw(q1)--(q2);
			\coordinate (a1) at (5,0.5);
			\coordinate (b1) at (9,0.5);
			\coordinate (c1) at (5,2.5);
			\coordinate (d1) at (9,2.5);
			\draw(a1)--(b1);
			\draw(b1)--(d1);
			\draw(d1)--(c1);
			\draw(c1)--(a1);
			\node(small-dot) at (6,2.8){$c_{1}$};
			\node(small-dot) at (7,3.2){$c_{2}$};
			\node(small-dot) at (7,2.8){$c_{3}$};
			\node(small-dot) at (6,0.2){$e_{1}$};
			\node(small-dot) at (7,0.2){$e_{2}$};
			\node(small-dot) at (8.5,0.2){$e_{3}$};
			\coordinate (pp1) at (5.5,2.1);
			\coordinate (pp2) at (6,2.1);
			\coordinate (pp3) at (6,1.1);
			\coordinate (pp4) at (7,1.1);
			\coordinate (qq1) at (5.5,1.9);
			\coordinate (qq2) at (6,1.9);
			\coordinate (qq3) at (6,0.9);
			\coordinate (qq2)--(qq3);
			\draw(qq3)--(qq4);
			\draw(pp1)--(pp2);
			\draw(pp2)--(pp3);
			\draw(pp3)--(pp4);
			\coordinate (a2) at (10,0.5);
			\coordinate (b2) at (14,0.5);
			\coordinate (c2) at (10,2.5);
			\coordinate (d2) at (14,2.5);
			\draw(a2)--(b2);
			\draw(b2)--(d2);
			\draw(d2)--(c2);
			\draw(c2)--(a2);
			\node(small-dot) at (11,2.8){$c_{1}$};
			\node(small-dot) at (12,3.2){$c_{2}$};
			\node(small-dot) at (12,2.8){$c_{3}$};
			\node(small-dot) at (11,0.2){$e_{1}$};
			\node(small-dot) at (12,0.2){$e_{2}$};
			\node(small-dot) at (13.5,0.2){$e_{3}$};
			\coordinate (ppp1) at (10.5,2.1);
			\coordinate (ppp2) at (13.5,2.1);
			\coordinate (ppp3) at (10.5,1.1);
			\coordinate (ppp4) at (13.5,1.1);
			\coordinate (qqq1) at (10.5,1.9);
			\coordinate (qqq2) at (13.5,1.9);
			\coordinate (qqq3) at (10.5,0.9);
			\coordinate (qqq4) at (13.5,0.9);
			\draw(ppp1)--(ppp2);
			\draw(ppp3)--(ppp4);
			\draw(qqq1)--(qqq2);
			\draw(qqq3)--(qqq4);
			\node(small-dot) at (8.5,0.9){$\bullet$};
			\node(small-dot) at (2,2.1){$\circ$};
			\node(small-dot) at (1.7,1.7){$C_{1}$};
			\coordinate (s1) at (1,2.5);
			\coordinate (t1) at (1,0.5);
			\draw[dotted](s1)--(t1);
			\coordinate (s2) at (2,2.5);
			\coordinate (t2) at (2,0.5);
			\draw[dotted](s2)--(t2);
			\coordinate (s3) at (3.5,2.5);
			\coordinate (t3) at (3.5,0.5);
			\draw[dotted](s3)--(t3);
			\coordinate (ss1) at (6,2.5);
			\coordinate (tt1) at (6,0.5);
			\draw[dotted](ss1)--(tt1);
			\coordinate (ss2) at (7,2.5);
			\coordinate (tt2) at (7,0.5);
			\draw[dotted](ss2)--(tt2);
			\coordinate (ss3) at (8.5,2.5);
			\coordinate (tt3) at (8.5,0.5);
			\draw[dotted](ss3)--(tt3);
			\coordinate (sss1) at (11,2.5);
			\coordinate (ttt1) at (11,0.5);
			\draw[dotted](sss1)--(ttt1);
			\coordinate (sss2) at (12,2.5);
			\coordinate (ttt2) at (12,0.5);
			\draw[dotted](sss2)--(ttt2);
			\coordinate (sss3) at (13.5,2.5);
			\coordinate (ttt3) at (13.5,0.5);
			\draw[dotted](sss3)--(ttt3);
		\end{tikzpicture}
	\end{figure}	\FloatBarrier
	\par \noindent \textit{Case $P\in \bar{C}$, $Q\in \bar{D}$.} $x_{1,c_{3}}\in P-Q$ implies $c_{2}=c_{3}$. According to $y$-parts and $z$-parts, we have $c_{1}=d_{3}$ and $d_{2}=n$, so $c_{1}=c_{2}=c_{3}=d_{2}=d_{3}=n$. The unique point in $Q-P$ is $x_{2,d_{1}}$ since $d_{1}<d_{2}$. Hence we have $d_{1}=n-1$.
	
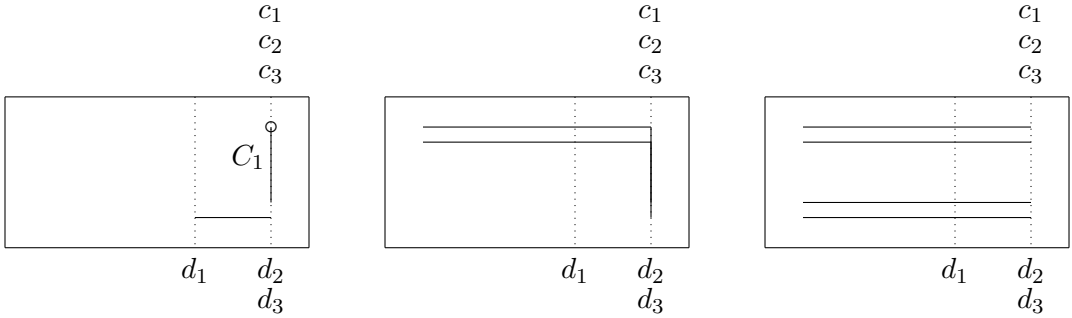
\begin{figure}[h]
		\centering
		\begin{tikzpicture}
			\coordinate (a) at (0,0.5);
			\coordinate (b) at (4,0.5);
			\coordinate (c) at (0,2.5);
			\coordinate (d) at (4,2.5);
			\draw(a)--(b);
			\draw(b)--(d);
			\draw(d)--(c);
			\draw(c)--(a);
			\node(small-dot) at (3.5,3.6){$c_{1}$};
			\node(small-dot) at (3.5,3.2){$c_{2}$};
			\node(small-dot) at (3.5,2.8){$c_{3}$};
			\node(small-dot) at (2.5,0.2){$d_{1}$};
			\node(small-dot) at (3.5,0.2){$d_{2}$};
			\node(small-dot) at (3.5,-0.2){$d_{3}$};
			\coordinate (p1) at (3.5,2.1);
			\coordinate (p2) at (3.5,1.1);
			\coordinate (q1) at (2.5,0.9);
			\coordinate (q2) at (3.5,0.9);
			\draw(p1)--(p2);
			\draw(q1)--(q2);
			\coordinate (a1) at (5,0.5);
			\coordinate (b1) at (9,0.5);
			\coordinate (c1) at (5,2.5);
			\coordinate (d1) at (9,2.5);
			\draw(a1)--(b1);
			\draw(b1)--(d1);
			\draw(d1)--(c1);
			\draw(c1)--(a1);
			\node(small-dot) at (8.5,3.6){$c_{1}$};
			\node(small-dot) at (8.5,3.2){$c_{2}$};
			\node(small-dot) at (8.5,2.8){$c_{3}$};
			\node(small-dot) at (7.5,0.2){$d_{1}$};
			\node(small-dot) at (8.5,0.2){$d_{2}$};
			\node(small-dot) at (8.5,-0.2){$d_{3}$};
			\coordinate (pp1) at (5.5,2.1);
			\coordinate (pp2) at (8.5,2.1);
			\coordinate (pp3) at (8.5,1.1);
			\coordinate (qq1) at (5.5,1.9);
			\coordinate (qq3) at (8.5,1.9);
			\coordinate (qq4) at (8.5,0.9);
			\draw(qq1)--(qq3);
			\draw(qq3)--(qq4);
			\draw(pp1)--(pp2);
			\draw(pp2)--(pp3);
			\coordinate (a2) at (10,0.5);
			\coordinate (b2) at (14,0.5);
			\coordinate (c2) at (10,2.5);
			\coordinate (d2) at (14,2.5);
			\draw(a2)--(b2);
			\draw(b2)--(d2);
			\draw(d2)--(c2);
			\draw(c2)--(a2);
			\node(small-dot) at (13.5,3.6){$c_{1}$};
			\node(small-dot) at (13.5,3.2){$c_{2}$};
			\node(small-dot) at (13.5,2.8){$c_{3}$};
			\node(small-dot) at (12.5,0.2){$d_{1}$};
			\node(small-dot) at (13.5,0.2){$d_{2}$};
			\node(small-dot) at (13.5,-0.2){$d_{3}$};
			\coordinate (ppp1) at (10.5,2.1);
			\coordinate (ppp2) at (13.5,2.1);
			\coordinate (ppp3) at (10.5,1.1);
			\coordinate (ppp4) at (13.5,1.1);
			\coordinate (qqq1) at (10.5,1.9);
			\coordinate (qqq2) at (13.5,1.9);
			\coordinate (qqq3) at (10.5,0.9);
			\coordinate (qqq4) at (13.5,0.9);
			\draw(ppp1)--(ppp2);
			\draw(ppp3)--(ppp4);
			\draw(qqq1)--(qqq2);
			\draw(qqq3)--(qqq4);
			\node(small-dot) at (3.5,2.1){$\circ$};
			\node(small-dot) at (3.2,1.7){$C_{1}$};
			\coordinate (s1) at (2.5,2.5);
			\coordinate (t1) at (2.5,0.5);
			\draw[dotted](s1)--(t1);
			\coordinate (s2) at (3.5,2.5);
			\coordinate (t2) at (3.5,0.5);
			\draw[dotted](s2)--(t2);
			\coordinate (ss1) at (7.5,2.5);
			\coordinate (tt1) at (7.5,0.5);
			\draw[dotted](ss1)--(tt1);
			\coordinate (ss2) at (8.5,2.5);
			\coordinate (tt2) at (8.5,0.5);
			\draw[dotted](ss2)--(tt2);
			\coordinate (sss1) at (12.5,2.5);
			\coordinate (ttt1) at (12.5,0.5);
			\draw[dotted](sss1)--(ttt1);
			\coordinate (sss2) at (13.5,2.5);
			\coordinate (ttt2) at (13.5,0.5);
			\draw[dotted](sss2)--(ttt2);
		\end{tikzpicture}
	\end{figure}	\FloatBarrier
	\par \noindent \textit{Case $P\in \bar{A}$, $Q\in \bar{F}$ or $\bar{E}$ .} $P-Q$ contains at least two points $x_{1,a_{1}}$ and $y_{1,n}$, contradiction.
	\par \noindent \textit{Case $P\in \bar{A}$, $Q\in \bar{D}$.} $x_{1,a_{1}}\in P-Q$ implies $a_{1}=a_{r}$ and $d_{3}=n$, so the unique point in $Q-P$ is $y_{2,d_{3}}$. Hence we have $a_{1}=d_{1}$ and $a_{2}=d_{2}$.
	
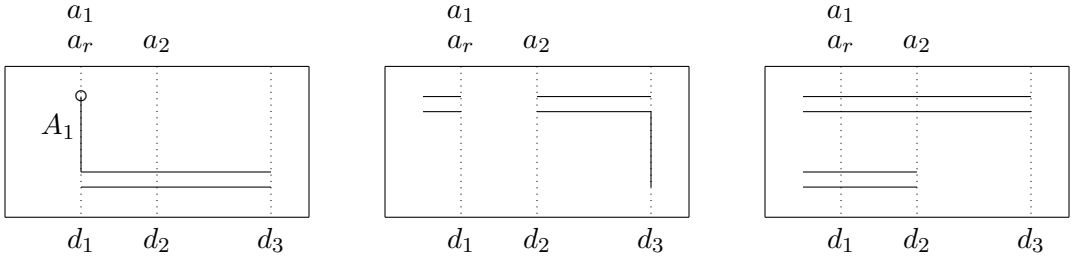
\begin{figure}[h]
		\centering
		\begin{tikzpicture}
			\coordinate (a) at (0,0.5);
			\coordinate (b) at (4,0.5);
			\coordinate (c) at (0,2.5);
			\coordinate (d) at (4,2.5);
			\draw(a)--(b);
			\draw(b)--(d);
			\draw(d)--(c);
			\draw(c)--(a);
			\node(small-dot) at (1,3.2){$a_{1}$};
			\node(small-dot) at (2,2.8){$a_{2}$};
			\node(small-dot) at (1,2.8){$a_{r}$};
			\node(small-dot) at (1,0.2){$d_{1}$};
			\node(small-dot) at (2,0.2){$d_{2}$};
			\node(small-dot) at (3.5,0.2){$d_{3}$};
			\coordinate (p1) at (1,2.1);
			\coordinate (p2) at (1,1.1);
			\coordinate (p3) at (3.5,1.1);
			\coordinate (q1) at (1,0.9);
			\coordinate (q2) at (3.5,0.9);
			\draw(p1)--(p2);
			\draw(p2)--(p3);
			\draw(q1)--(q2);
			\coordinate (a1) at (5,0.5);
			\coordinate (b1) at (9,0.5);
			\coordinate (c1) at (5,2.5);
			\coordinate (d1) at (9,2.5);
			\draw(a1)--(b1);
			\draw(b1)--(d1);
			\draw(d1)--(c1);
			\draw(c1)--(a1);
			\node(small-dot) at (6,3.2){$a_{1}$};
			\node(small-dot) at (7,2.8){$a_{2}$};
			\node(small-dot) at (6,2.8){$a_{r}$};
			\node(small-dot) at (6,0.2){$d_{1}$};
			\node(small-dot) at (7,0.2){$d_{2}$};
			\node(small-dot) at (8.5,0.2){$d_{3}$};
			\coordinate (pp1) at (5.5,2.1);
			\coordinate (pp2) at (6,2.1);
			\coordinate (pp3) at (7,2.1);
			\coordinate (pp4) at (8.5,2.1);
			\coordinate (qq1) at (5.5,1.9);
			\coordinate (qq2) at (6,1.9);
			\coordinate (qq3) at (7,1.9);
			\coordinate (qq4) at (8.5,1.9);
			\coordinate (qq5) at (8.5,0.9);
			\draw(pp1)--(pp2);
			\draw(pp3)--(pp4);
			\draw(qq1)--(qq2);
			\draw(qq3)--(qq4);
			\draw(qq4)--(qq5);
			\coordinate (a2) at (10,0.5);
			\coordinate (b2) at (14,0.5);
			\coordinate (c2) at (10,2.5);
			\coordinate (d2) at (14,2.5);
			\draw(a2)--(b2);
			\draw(b2)--(d2);
			\draw(d2)--(c2);
			\draw(c2)--(a2);
			\node(small-dot) at (11,3.2){$a_{1}$};
			\node(small-dot) at (12,2.8){$a_{2}$};
			\node(small-dot) at (11,2.8){$a_{r}$};
			\node(small-dot) at (11,0.2){$d_{1}$};
			\node(small-dot) at (12,0.2){$d_{2}$};
			\node(small-dot) at (13.5,0.2){$d_{3}$};
			\coordinate (ppp1) at (10.5,2.1);
			\coordinate (ppp2) at (13.5,2.1);
			\coordinate (ppp3) at (10.5,1.1);
			\coordinate (ppp4) at (12,1.1);
			\coordinate (qqq1) at (10.5,1.9);
			\coordinate (qqq2) at (13.5,1.9);
			\coordinate (qqq3) at (10.5,0.9);
			\coordinate (qqq4) at (12,0.9);
			\draw(ppp1)--(ppp2);
			\draw(ppp3)--(ppp4);
			\draw(qqq1)--(qqq2);
			\draw(qqq3)--(qqq4);
			\node(small-dot) at (1,2.1){$\circ$};
			\node(small-dot) at (0.7,1.7){$A_{1}$};
			\coordinate (s1) at (1,2.5);
			\coordinate (t1) at (1,0.5);
			\draw[dotted](s1)--(t1);
			\coordinate (s2) at (2,2.5);
			\coordinate (t2) at (2,0.5);
			\draw[dotted](s2)--(t2);
			\coordinate (s3) at (3.5,2.5);
			\coordinate (t3) at (3.5,0.5);
			\draw[dotted](s3)--(t3);
			\coordinate (ss1) at (6,2.5);
			\coordinate (tt1) at (6,0.5);
			\draw[dotted](ss1)--(tt1);
			\coordinate (ss2) at (7,2.5);
			\coordinate (tt2) at (7,0.5);
			\draw[dotted](ss2)--(tt2);
			\coordinate (ss3) at (8.5,2.5);
			\coordinate (tt3) at (8.5,0.5);
			\draw[dotted](ss3)--(tt3);
			\coordinate (sss1) at (11,2.5);
			\coordinate (ttt1) at (11,0.5);
			\draw[dotted](sss1)--(ttt1);
			\coordinate (sss2) at (12,2.5);
			\coordinate (ttt2) at (12,0.5);
			\draw[dotted](sss2)--(ttt2);
			\coordinate (sss3) at (13.5,2.5);
			\coordinate (ttt3) at (13.5,0.5);
			\draw[dotted](sss3)--(ttt3);
		\end{tikzpicture}
	\end{figure}	\FloatBarrier
	\par \noindent \textit{Case $P\in \bar{A}$, $Q\in \bar{C}$.} The second row of $P_{y}$ is empty so $c_{1}=c_{2}$ and $Q-P=\{y_{2,c_{1}}\}$. Hence we have $a_{2}=n$ according to $z$-parts. Suppose $P-Q=\{v\}$ and there is two possible cases $v\in P_{x}$ or $P_{y}$.
	\par If $v\in P_{x}$, then $c_{1}=n$ since $y_{1,n}\in P_{y}$. We have $c_{1}=c_{2}=c_{3}=n$, $a_{1}=n-1$, and $a_{r}=n$. 
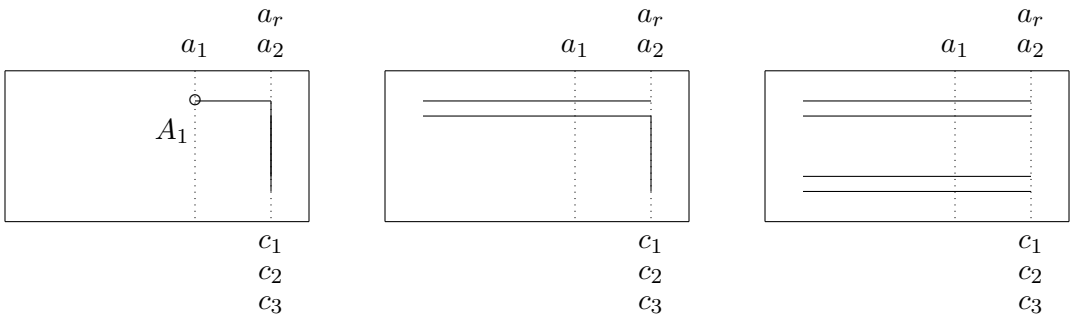
	\begin{figure}[!h]
		\centering
		\begin{tikzpicture}
			\coordinate (a) at (0,0.5);
			\coordinate (b) at (4,0.5);
			\coordinate (c) at (0,2.5);
			\coordinate (d) at (4,2.5);
			\draw(a)--(b);
			\draw(b)--(d);
			\draw(d)--(c);
			\draw(c)--(a);
			\node(small-dot) at (2.5,2.8){$a_{1}$};
			\node(small-dot) at (3.5,2.8){$a_{2}$};
			\node(small-dot) at (3.5,3.2){$a_{r}$};
			\node(small-dot) at (3.5,0.2){$c_{1}$};
			\node(small-dot) at (3.5,-0.2){$c_{2}$};
			\node(small-dot) at (3.5,-0.6){$c_{3}$};
			\coordinate (p1) at (2.5,2.1);
			\coordinate (p2) at (3.5,2.1);
			\coordinate (p3) at (3.5,1.1);
			\coordinate (q1) at (3.5,1.9);
			\coordinate (q2) at (3.5,0.9);
			\draw(p1)--(p2);
			\draw(p2)--(p3);
			\draw(q1)--(q2);
			\coordinate (a1) at (5,0.5);
			\coordinate (b1) at (9,0.5);
			\coordinate (c1) at (5,2.5);
			\coordinate (d1) at (9,2.5);
			\draw(a1)--(b1);
			\draw(b1)--(d1);
			\draw(d1)--(c1);
			\draw(c1)--(a1);
			\node(small-dot) at (7.5,2.8){$a_{1}$};
			\node(small-dot) at (8.5,2.8){$a_{2}$};
			\node(small-dot) at (8.5,3.2){$a_{r}$};
			\node(small-dot) at (8.5,0.2){$c_{1}$};
			\node(small-dot) at (8.5,-0.2){$c_{2}$};
			\node(small-dot) at (8.5,-0.6){$c_{3}$};
			\coordinate (pp1) at (5.5,2.1);
			\coordinate (pp2) at (8,2.1);
			\coordinate (pp3) at (8.5,2.1);
			\coordinate (qq1) at (5.5,1.9);
			\coordinate (qq2) at (8.5,1.9);
			\coordinate (qq3) at (8.5,0.9);
			\draw(pp1)--(pp3);
			\draw(qq1)--(qq2);
			\draw(qq2)--(qq3);
			\coordinate (a2) at (10,0.5);
			\coordinate (b2) at (14,0.5);
			\coordinate (c2) at (10,2.5);
			\coordinate (d2) at (14,2.5);
			\draw(a2)--(b2);
			\draw(b2)--(d2);
			\draw(d2)--(c2);
			\draw(c2)--(a2);
			\node(small-dot) at (12.5,2.8){$a_{1}$};
			\node(small-dot) at (13.5,2.8){$a_{2}$};
			\node(small-dot) at (13.5,3.2){$a_{r}$};
			\node(small-dot) at (13.5,0.2){$c_{1}$};
			\node(small-dot) at (13.5,-0.2){$c_{2}$};
			\node(small-dot) at (13.5,-0.6){$c_{3}$};
			\coordinate (ppp1) at (10.5,2.1);
			\coordinate (ppp2) at (13.5,2.1);
			\coordinate (ppp3) at (10.5,1.1);
			\coordinate (ppp4) at (13.5,1.1);
			\coordinate (qqq1) at (10.5,1.9);
			\coordinate (qqq2) at (13.5,1.9);
			\coordinate (qqq3) at (10.5,0.9);
			\coordinate (qqq4) at (13.5,0.9);
			\draw(ppp1)--(ppp2);
			\draw(ppp3)--(ppp4);
			\draw(qqq1)--(qqq2);
			\draw(qqq3)--(qqq4);
			\node(small-dot) at (2.5,2.1){$\circ$};
			\node(small-dot) at (2.2,1.7){$A_{1}$};
			\coordinate (s1) at (2.5,2.5);
			\coordinate (t1) at (2.5,0.5);
			\draw[dotted](s1)--(t1);
			\coordinate (s2) at (3.5,2.5);
			\coordinate (t2) at (3.5,0.5);
			\draw[dotted](s2)--(t2);
			\coordinate (ss1) at (7.5,2.5);
			\coordinate (tt1) at (7.5,0.5);
			\draw[dotted](ss1)--(tt1);
			\coordinate (ss2) at (8.5,2.5);
			\coordinate (tt2) at (8.5,0.5);
			\draw[dotted](ss2)--(tt2);
			\coordinate (sss1) at (12.5,2.5);
			\coordinate (ttt1) at (12.5,0.5);
			\draw[dotted](sss1)--(ttt1);
			\coordinate (sss2) at (13.5,2.5);
			\coordinate (ttt2) at (13.5,0.5);
			\draw[dotted](sss2)--(ttt2);
		\end{tikzpicture}
	\end{figure}	\FloatBarrier
	\par If $v\in P_{y}$, then $a_{1}=c_{2}$ since $y_{2,c_{2}}\in Q-P$ and $a_{r}=c_{3}$ since $x_{2,a_{r}}\in P_{x}$.
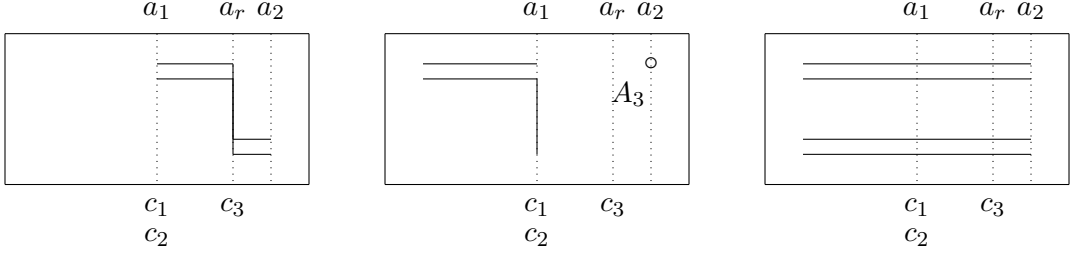
	\begin{figure}[h]
		\centering
		\begin{tikzpicture}
			\coordinate (a) at (0,0.5);
			\coordinate (b) at (4,0.5);
			\coordinate (c) at (0,2.5);
			\coordinate (d) at (4,2.5);
			\draw(a)--(b);
			\draw(b)--(d);
			\draw(d)--(c);
			\draw(c)--(a);
			\node(small-dot) at (2,2.8){$a_{1}$};
			\node(small-dot) at (3.5,2.8){$a_{2}$};
			\node(small-dot) at (3,2.8){$a_{r}$};
			\node(small-dot) at (2,0.2){$c_{1}$};
			\node(small-dot) at (2,-0.2){$c_{2}$};
			\node(small-dot) at (3,0.2){$c_{3}$};
			\coordinate (p1) at (2,2.1);
			\coordinate (p2) at (3,2.1);
			\coordinate (p3) at (3,1.1);
			\coordinate (p4) at (3.5,1.1);
			\coordinate (q1) at (2,1.9);
			\coordinate (q2) at (3,1.9);
			\coordinate (q3) at (3,0.9);
			\coordinate (q4) at (3.5,0.9);
			\draw(p1)--(p2);
			\draw(p2)--(p3);
			\draw(p3)--(p4);
			\draw(q1)--(q2);
			\draw(q2)--(q3);
			\draw(q3)--(q4);
			\coordinate (a1) at (5,0.5);
			\coordinate (b1) at (9,0.5);
			\coordinate (c1) at (5,2.5);
			\coordinate (d1) at (9,2.5);
			\draw(a1)--(b1);
			\draw(b1)--(d1);
			\draw(d1)--(c1);
			\draw(c1)--(a1);
			\node(small-dot) at (7,2.8){$a_{1}$};
			\node(small-dot) at (8.5,2.8){$a_{2}$};
			\node(small-dot) at (8,2.8){$a_{r}$};
			\node(small-dot) at (7,0.2){$c_{1}$};
			\node(small-dot) at (7,-0.2){$c_{2}$};
			\node(small-dot) at (8,0.2){$c_{3}$};
			\coordinate (pp1) at (5.5,2.1);
			\coordinate (pp2) at (7,2.1);
			\coordinate (qq1) at (5.5,1.9);
			\coordinate (qq2) at (7,1.9);
			\coordinate (qq3) at (7,0.9);
			\draw(qq1)--(qq2);
			\draw(qq2)--(qq3);
			\draw(pp1)--(pp2);
			\coordinate (a2) at (10,0.5);
			\coordinate (b2) at (14,0.5);
			\coordinate (c2) at (10,2.5);
			\coordinate (d2) at (14,2.5);
			\draw(a2)--(b2);
			\draw(b2)--(d2);
			\draw(d2)--(c2);
			\draw(c2)--(a2);
			\node(small-dot) at (12,2.8){$a_{1}$};
			\node(small-dot) at (13.5,2.8){$a_{2}$};
			\node(small-dot) at (13,2.8){$a_{r}$};
			\node(small-dot) at (12,0.2){$c_{1}$};
			\node(small-dot) at (12,-0.2){$c_{2}$};
			\node(small-dot) at (13,0.2){$c_{3}$};
			\coordinate (ppp1) at (10.5,2.1);
			\coordinate (ppp2) at (13.5,2.1);
			\coordinate (ppp3) at (10.5,1.1);
			\coordinate (ppp4) at (13.5,1.1);
			\coordinate (qqq1) at (10.5,1.9);
			\coordinate (qqq2) at (13.5,1.9);
			\coordinate (qqq3) at (10.5,0.9);
			\coordinate (qqq4) at (13.5,0.9);
			\draw(ppp1)--(ppp2);
			\draw(ppp3)--(ppp4);
			\draw(qqq1)--(qqq2);
			\draw(qqq3)--(qqq4);
			\node(small-dot) at (8.5,2.1){$\circ$};
			\node(small-dot) at (8.2,1.7){$A_{3}$};
			\coordinate (s1) at (2,2.5);
			\coordinate (t1) at (2,0.5);
			\draw[dotted](s1)--(t1);
			\coordinate (s2) at (3,2.5);
			\coordinate (t2) at (3,0.5);
			\draw[dotted](s2)--(t2);
			\coordinate (s3) at (3.5,2.5);
			\coordinate (t3) at (3.5,0.5);
			\draw[dotted](s3)--(t3);
			\coordinate (ss1) at (7,2.5);
			\coordinate (tt1) at (7,0.5);
			\draw[dotted](ss1)--(tt1);
			\coordinate (ss2) at (8,2.5);
			\coordinate (tt2) at (8,0.5);
			\draw[dotted](ss2)--(tt2);
			\coordinate (ss3) at (8.5,2.5);
			\coordinate (tt3) at (8.5,0.5);
			\draw[dotted](ss3)--(tt3);
			\coordinate (sss1) at (12,2.5);
			\coordinate (ttt1) at (12,0.5);
			\draw[dotted](sss1)--(ttt1);
			\coordinate (sss2) at (13,2.5);
			\coordinate (ttt2) at (13,0.5);
			\draw[dotted](sss2)--(ttt2);
			\coordinate (sss3) at (13.5,2.5);
			\coordinate (ttt3) at (13.5,0.5);
			\draw[dotted](sss3)--(ttt3);
		\end{tikzpicture}
	\end{figure}	\FloatBarrier
	\par According to above discussion, $c_{2}(P)\subset c_{1}(P)$ for facet $P$ in generic condition. However, for non-generic condition, some points vanishing in $c_{1}(P)$ will appear in $c_{2}(P)$. We summarize these additive conditions into the following table:
	\begin{figure}[h]
		\centering
		\begin{tikzpicture}
			\coordinate (a2) at (0,4.5);
			\coordinate (b2) at (10,4.5);
			\coordinate (a3) at (0,3);
			\coordinate (b3) at (10,3);
			\coordinate (a4) at (0,2);
			\coordinate (b4) at (10,2);
			\coordinate (a5) at (0,1);
			\coordinate (b5) at (10,1);
			\coordinate (a6) at (0,0);
			\coordinate (b6) at (10,0);
			\draw(a2)--(b2);
			\draw(a3)--(b3);
			\draw(a4)--(b4);
			\draw(a5)--(b5);
			\draw(a6)--(b6);
			\coordinate (c1) at (1,4.5);
			\coordinate (c2) at (5,4.5);
			\coordinate (c3) at (6,4.5);
			\coordinate (c4) at (10,4.5);
			\coordinate (d1) at (1,0);
			\coordinate (d2) at (5,3);
			\coordinate (d3) at (6,3);
			\coordinate (d4) at (10,0);
			\draw(a2)--(a6);
			\draw(b2)--(b6);
			\draw(c3)--(d3);
			\draw(c4)--(d4);
			\draw(c2)--(d2);
			\draw(c1)--(d1);
			\coordinate (t1) at (5,2);
			\coordinate (t2) at (5,1);
			\coordinate (t3) at (6,2);
			\coordinate (t4) at (6,1);
			\draw(t1)--(t2);
			\draw(t3)--(t4);
			\node(small-dot) at (0.5,3.75){$A_{1}$};
			\node(small-dot) at (5.5,3.75){$A_{3}$};
			\node(small-dot) at (2.05,4.25){\ding{172}$a_{1}=a_{r}$}; 
			\node(small-dot) at (2.4,3.75){\ding{173}$a_{1}=n-1$,}; \node(small-dot) at (2.3,3.25){$a_{2}=a_{r}=n$};
			\node(small-dot) at (6.9,3.75){$a_{2}=n$};
			\node(small-dot) at (0.5,2.5){$C_{1}$};
			\node(small-dot) at (1.88,2.5){$c_{2}=c_{3}$};
			\node(small-dot) at (0.5,1.5){$D_{1}$};
			\node(small-dot) at (5.5,1.5){$D_{2}$};
			\node(small-dot) at (2.34,1.75){$d_{1}=d_{2}-1$,}; \node(small-dot) at (2.7,1.25){$d_{2}=d_{3}\le n-1$};
			\node(small-dot) at (7,1.5){$d_{2}=d_{3}$};
			\node(small-dot) at (0.5,0.5){$E_{2}$};
			\node(small-dot) at (1.85,0.5){$e_{1}=1$};
			\node(small-dot) at (5,5){Additive conditions};
		\end{tikzpicture}
	\end{figure}	\FloatBarrier
	\par We prove the first part of our main theorem as the end of this section.\\\\
	\textbf{Theorem 5.4. }\textit{$\Delta_{0}$ is shellable and the $*$-ordering of facets gives a shelling of $\Delta_{0}$.}\\\\
	\textit{Proof.} Suppose $P>Q$ are two facets of $\Delta_{0}$. By Corollary 4.4, it suffices to show that $Q$ cannot cover $c(P)$. Then $*$-ordering is a shelling of $\Delta_{0}$.
	\par We first deal with cases that $P$ and $Q$ are from different families of facets of $\Delta_{0}$. If $P\in \bar{A}$ and $Q\in \bar{D}$, $\bar{E}$, or $\bar{F}$, then $A_{1}$ exists and cannot be covered by $Q$. If $P\in \bar{A}$ and $Q\in \bar{C}$, then $A_{1}$ and $A_{3}$ exist and cannot be covered by $Q$, otherwise $c_{2}\le a_{1}<a_{2}\le c_{1}$ leads to contribution. If $P\in \bar{C}$, then $C_{1}$ exists and cannot be covered by $Q$. If $P\in \bar{D}$ and $Q\in \bar{E}$, then $D_{2}$ exists. If $Q$ covers $D_{2}$, then we have $d_{1}<d_{2}\le e_{1}\le e_{2}\le n-1$. Hence $D_{1}$ exists and cannot be covered by $Q$ since $d_{1}<e_{2}$. If $P\in \bar{D}$ or $\bar{E}$ and $Q\in \bar{F}$, then $D_{2}$ or $E_{2}$ exist and cannot be covered by $Q$.
	\par Then we deal with cases that $P$ and $Q$ are from the same families of facets of $\Delta_{0}$. The idea of proof is to find contradiction under the assumption $Q$ cover $c(P)$. For convenience, we use $a_{i}(P)$ and $a_{i}(Q)$ to represent axes of $P$, $Q\in \bar{A}$, which is similar for $\bar{C}$, $\bar{D}$, $\bar{E}$, $\bar{F}$.
	\par \textit{Case $P$, $Q\in \bar{A}$. }$A_{1}$ always exists, so $a_{1}(P)\ge a_{1}(Q)$. Since $P>Q$, we have $a_{1}(P)=a_{1}(Q)$. If $a_{1}(P)\neq a_{r}(P)$, then $A_{2}$ exists and $a_{r}(P)=a_{r}(Q)$. If $a_{1}(P)=a_{r}(P)$, then $a_{r}(Q)=a_{r}(P)$ since $P>Q$. $A_{3}$ always exists so $a_{2}(P)=a_{2}(Q)$. All the three axes are at the same position, which is contradictory to $P>Q$. 
	\par \textit{Case $P$, $Q\in \bar{C}$. }$C_{1}$ always exists, so $c_{2}(P)\ge c_{2}(Q)$. Since $P>Q$, we have $c_{2}(P)=c_{2}(Q)$. Similarly to the discussion in $\bar{A}$, we have $c_{3}(P)=c_{3}(Q)$. If $c_{1}(P)\neq1$, then $C_{3}$ exists and $c_{1}(P)=c_{1}(Q)$. If $c_{1}(P)=1$, then $c_{1}(Q)=1$ since $P>Q$. All the three axes are at the same position, which is contradictory to $P>Q$.
	\par \textit{Case $P$, $Q\in \bar{D}$. }$D_{1}$ vanishes only when $d_{1}(P)=n-1$ and $d_{2}(P)=d_{3}(P)=n$. In this condition, $P$ is the minimal element in $\bar{D}$, so we only have to consider cases that $D_{1}$ exists. We have $e_{1}(P)=e_{1}(Q)$ since $P>Q$. $D_{2}$ exists always, so $e_{2}(P)=e_{2}(Q)$ and $e_{3}(P)=e_{3}(Q)$ similarly to the discussion in $\bar{A}$, which comes to contradiction.
	\par \textit{Case $P$, $Q\in \bar{E}$. }If $e_{2}(P)=n-1$, then $e_{2}(Q)=n-1$ since $P>Q$. If $e_{2}(P)\neq n-1$, then $E_{1}$ exists, which implies $e_{2}(P)=e_{2}(Q)$. $E_{2}$ always exists so $e_{1}(P)=e_{1}(Q)$. Similarly to the discussion of $e_{2}$, we have $e_{3}(P)=e_{3}(Q)$. All the three axes are at the same position, which is contradictory to $P>Q$.
	\par \textit{Case $P$, $Q\in \bar{F}$. }Similarly to discussion in $\bar{E}$, we have $f_{1}(P)=f_{1}(Q)$ and $f_{2}(P)=f_{2}(Q)$. If $f_{2}(P)\neq f_{3}(P)-1$, then $F_{3}$ exists and $f_{3}(P)=f_{3}(Q)$. If $f_{2}(P)=f_{3}(P)-1$, then we have $f_{3}(P)=f_{3}(Q)$ since $P>Q$. The two conditions both lead to contradiction.   \hfill{$\square$} 
	
	\section{Hilbert Series of $\mathcal{L}^{2,n}_{2,2}$}
	\par We compute $h_{j}(T)$ (recall the notation in Definition 2.12) for $T=\bar{A},\bar{C},\bar{D},\bar{E},\bar{F}$ and $j=1,2,3$ according to the vanishing condition table and additive condition table. The results are as follows:
	\par $h_{1}(\bar{A})=0$,
	\par $h_{1}(\bar{C})=n$, for condition $c_{1}=1,1\le c_{2}=c_{3}\le n$,
	\par $h_{1}(\bar{D})=1$, for condition $d_{1}=n-1,d_{2}=d_{3}=n$,
	\par $h_{1}(\bar{E})=n-1$, for condition $1\le e_{1}\le n-1,e_{2}=n-1,e_{3}=n$,
	\par $h_{1}(\bar{F})=n-3$, for condition $1\le f_{1}\le n-3,f_{2}=n-1,f_{3}=n$,
	\par $h_{2}(\bar{A})=\binom{n}{2}=\frac{n(n-1)}{2}$, for condition $1\le a_{1}=a_{r}<a_{2}\le n$,
	\par $h_{2}(\bar{C})=\binom{n}{2}+\binom{n-1}{2}+\binom{n-1}{1}=n(n-1)$, for conditions $1=c_{1}\le c_{2}<c_{3}\le n$ and $1<c_{1}\le c_{2}=c_{3}\le n$,
	\par $h_{2}(\bar{D})=\binom{n}{2}-\binom{n-1}{1}+n-2=\frac{(n+1)(n-2)}{2}$, for conditions $1\le d_{1}<d_{2}-1=d_{3}-1\le n-1$ and $1\le d_{1}=d_{2}-1=d_{3}-1\le n-1$,
	\par $h_{2}(\bar{E})=\binom{n-2}{2}+\binom{n-2}{1}=\frac{(n-1)(n-2)}{2}$, for condition $1\le e_{1}\le e_{2}<n-1,e_{3}=n$,
	\par $h_{2}(\bar{F})=\binom{n-2}{2}=\frac{(n-2)(n-3)}{2}$, for condition $1\le f_{1}<f_{2}=f_{3}-1<n-1$,
	\par $h_{3}(\bar{A})=\underset{i=1}{\overset{n-1}{\sum}}i^{2}=\frac{n(n-1)(2n-1)}{6}$, for condition $1\le a_{1}<a_{r}\le n,\le a_{1}<a_{2}\le n$,
	\par $h_{3}(\bar{C})=\binom{n-1}{3}+\binom{n-1}{2}=\frac{n(n-1)(n-2)}{6}$, for condition $1<c_{1}\le c_{2}<c_{3}\le n$,
	\par $h_{3}(\bar{D})=\binom{n}{3}=\frac{n(n-1)(n-2)}{6}$, for condition $1\le d_{1}<d_{2}<d_{3}\le n$,
	\par $h_{3}(\bar{E})=\binom{n-1}{3}+\binom{n-1}{2}=\frac{n(n-1)(n-2)}{6}$, for condition $1\le e_{1}\le e_{2}<e_{3}<n$,
	\par $h_{3}(\bar{F})=\binom{n}{3}-\binom{n-1}{2}=\frac{(n-1)(n-2)(n-3)}{6}$, for condition $1\le f_{1}<f_{2}<f_{3}-1\le n-1$.
	\par With above results, we can finish the second part of our main theorem.\\\\
	\textbf{Theorem 6.1. }\textit{The Hilbert series of $\mathcal{L}^{2,n}_{2,2}$ is  $H_{\mathcal{L}^{2,n}_{2,2}}(z)=\big( \frac{1+(n-1)z}{(1-z)^{n+1}}\big)^{3}$.}
	\begin{proof}[Proof.]
		The number $h_{i}$ is just the sum of $h_{i}(\bar{A})$ to $h_{i}(\bar{F})$. As a result of above computation, we have $h_{0}=1$, $h_{1}=3(n-1)$, $h_{2}=3(n-1)^{2}$, $h_{3}=(n-1)^{3}$. Then we have 
		\begin{center}
			$H_{\mathcal{L}^{2,n}_{2,2}}(z)=\frac{h_{0}+h_{1}z+h_{2}z^{2}+h_{3}z^{3}}{(1-z)^{3n+3}}=\big( \frac{1+(n-1)z}{(1-z)^{n+1}}\big)^{3}$. 
		\end{center} 
		\par \noindent by Corollary 4.4, Theorem 5.4, and formula (5). 
	\end{proof} 
	\par It is natural to find some relation between the Hilbert series of determinantal varieties and their jets. For the complicacy of the Hilbert series of arbitrary order jets, we can only deal with some easy cases here. By formula (1), we have 
	\begin{center}
		$H_{\mathcal{L}^{m,n}_{2}}(z)=\frac{\underset{k=0}{\overset{m-1}{\sum}}\binom{m-1}{k} \binom{n-1}{k} z^{k}}{(1-z)^{m+n-1}}$
	\end{center}	
	\par \noindent and it is $H_{\mathcal{L}^{2,n}_{2}}(z)=\frac{1+(n-1)z}{(1-z)^{n+1}}$ for $m=2$. By Theorem 6.1 and (2), it is easy to see that $H_{\mathcal{L}^{m,n}_{2,1}}(z)=H_{\mathcal{L}^{m,n}_{2}}(z)^{2}$ and $H_{\mathcal{L}^{2,n}_{2,2}}(z)=H_{\mathcal{L}^{2,n}_{2}}(z)^{3}$. 
	\par For general condition, we have a conjecture about the Hilbert series:\\\\
	\textbf{Conjecture 6.2. }The Hilbert series of of $\mathcal{L}^{m,n}_{r,k}$ is the $(k+1)$-th power of the Hilbert series of $\mathcal{L}^{m,n}_{r}$. \\\\

	\end{document}